\documentclass[10pt]{article}
\usepackage[psamsfonts]{amssymb}
\usepackage{bm}
\usepackage{theorem}
\theorembodyfont{\upshape}

\usepackage{amsmath}
\usepackage{latexsym}

\begin{document}

\newcommand\qed{\hfill$\square$}

\newcounter{prop}

\addtocounter{prop}{1}

\newenvironment{prop}{\smallskip {\bf {\textsf{Proposition}}}\
{\bf{\textsf{\arabic{section}.\arabic{prop}.}}}\ \it }{\smallskip\addtocounter{prop}{1}}

\newenvironment{propbt}{\smallskip {\bf {\textsf{Proposition}}}\
{\bf {\textsf {\arabic{section}.\arabic{prop}}}}\ \it }{\smallskip\addtocounter{prop}{1}}

\newenvironment{lem}{\smallskip {\bf {\textsf{Lemma}}}\
{\bf {\textsf{\arabic{section}.\arabic{prop}.}}}\ \it }{\smallskip
\addtocounter{prop}{1}}

\newenvironment{thm}{\smallskip {\bf {\textsf{Theorem}}}\
{\bf {\textsf{\arabic{section}.\arabic{prop}.}}}\ \it }{\smallskip
\addtocounter{prop}{1}}

\newenvironment{cor}{\smallskip {\bf {\textsf{Corollary}}}\
{\bf {\textsf{\arabic{section}.\arabic{prop}.}}}\ \rm }{\smallskip
\addtocounter{prop}{1}}

\newenvironment{lembt}{\smallskip {\bf Lemma}\
{\bf \arabic{section}.\arabic{prop}}\ \it }{\smallskip }

\newenvironment{rem}{\smallskip {\bf {\textsf{Remark}}}\
{\bf {\textsf{\arabic{section}.\arabic{prop}.}}}\ \rm }{\smallskip
\addtocounter{prop}{1}}

\newenvironment{exmp}{\smallskip {\bf {\textsf{Example}}}\
{\bf {\textsf{\arabic{section}.\arabic{prop}.}}}\ \rm }{\smallskip
\addtocounter{prop}{1}}

\newenvironment{proof}{\smallskip {\it Proof.}\ \rm }{\medskip } 

\renewcommand{\thesection }{\normalsize{\textsf {\arabic{section}.}}}

\renewcommand{\thetheorem}{\normalsize{\textsf {\arabic{theorem}.}}}

\newcounter{sct}

\addtocounter{sct}{1}

\renewcommand{\thesubsection }{\normalsize{\textsf
{\arabic{section}.\arabic{subsection}.}}}

\renewcommand{\theequation}{\arabic{sct}.\arabic{equation}}

\newcommand{\vien}{\mbox{\rm 1\hspace{-.2em}I}\,}
\newcommand{\e}{\,\mbox{\rm e}}
\newcommand{\dd}{\,\mbox{\rm d}}
\newcommand{\dr}{{\rm d}\,}
\newcommand{\dek}{\mathop{\mbox{\large {\rm X}}\,}}
\newcommand{\lsup}{\mathop{\mbox{\rm limsup}\,}}
\newcommand{\linf}{\mathop{\mbox {\rm liminf}\,}}
\newcommand{\diam}{\mathop{\mbox {\rm diam}\,}}
\newcommand{\dist}{\mathop{\mbox {\rm dist}\,}}
\newcommand{\iz}[1]{_{{}_{#1}}}
\newcommand{\izz}[1]{\raisebox{-.6\height}{\scriptsize $#1$}}
\newcommand{\cit}[1]{~[#1]}
\newcommand{\ed}{\end{document}}
\newcommand{\br}[1]{\vert #1\vert }

\newcommand{\cF}{\mathcal {F}}
\newcommand{\cH}{\mathcal {H}}
\newcommand{\cG}{\mathcal {G}}
\newcommand{\cN}{\mathcal {N}}
\newcommand{\cM}{\mathcal {M}}
\newcommand{\cT}{\mathcal {T}}
\newcommand{\cP}{\mathcal {P}}
\newcommand{\cE}{\mathcal {E}}
\newcommand{\rra }{\alpha }
\newcommand{\rre }{\varepsilon }
\newcommand{\rrg }{\gamma }
\newcommand{\rrl }{\lambda }
\newcommand{\rrL }{\Lambda }
\newcommand{\rrk }{\kappa }
\newcommand{\rro }{\omega }
\newcommand{\rrO }{\Omega }
\newcommand{\rrr }{\rho }
\newcommand{\uxi }{\underline{\xi }}
\newcommand{\thh }{\theta }
\newcommand{\la  }{\leftarrow }

\newcommand{\ttD }{\mathbb{D}}
\newcommand{\ttA }{\mathbb{A}}
\newcommand{\ttE }{\mathbb{E}}
\newcommand{\ttN }{\mathbb{N}}
\newcommand{\ttR }{\mathbb{R}}
\newcommand{\ttZ }{\mathbb{Z}}
\newcommand{\ttB }{\mathbb{B}}
\newcommand{\ttS }{\mathbb{S}}
\newcommand{\ttV }{\mathbb{V}}
\newcommand{\ttG }{\mathfrak{G}}
\newcommand{\bP }{{\mathbb P}}
\newcommand{\bE }{{\mathbb E}}
\newcommand{\mro }{\mbox{\rm o}}
\newcommand{\mo }{\rm o}
\newcommand{\rO }{\mbox{\rm O}}
\newcommand{\spect }{\mbox{\rm Spect}\,}
\newcommand{\Spect }{\mbox{\rm Spect}\,}
\newcommand{\const }{\mbox{\rm const}\,}
\newcommand{\opt }{\mbox{\rm opt}\,}
\newcommand{\mconst }{\mbox{\small {\rm const}}}
\newcommand{\supp }{\mbox{\rm supp}\,}
\newcommand{\dly }{\geqslant }
\newcommand{\mly }{\leqslant }
\newcommand{\lk}{``}
\newcommand{\rk}{''}
\newcommand{\lf}{\left }
\newcommand{\rg }{\right }
\newcommand{\lv}{\left \vert }
\newcommand{\ve }{\vert }
\newcommand{\tus }{\varnothing }
\newcommand{\pr }{^\prime }
\newcommand{\zv }{^\ast }
\newcommand{\wdt }{\widetilde }
\newcommand{\wdh }{\widehat }
\newcommand{\ds }{\displaystyle }

\newcommand{\un}{\underline}

\newcommand{\bunderline}[1]{\underline{#1\mkern-4mu}\mkern4mu }
\newcommand{\bun}[1]{\underline{#1\mkern-7mu}\mkern7mu }

\newcommand{\wt}{\widetilde}
\newcommand{\rrf}{\varphi}
\newcommand{\n}{\nonumber }

\newcommand{\gs }{\hfill$\square$\medskip }

\newcommand{\nn}{\nonumber }
\newcommand{\beq}{\begin{eqnarray}}
\newcommand{\eeq}{\end{eqnarray}}
\newcommand{\beqz}{\begin{eqnarray*}}
\newcommand{\eeqz}{\end{eqnarray*}}
\newcommand{\rrG }{\Gamma }
\newcommand{\rrd}{\partial}
\newcommand{\mdd}{\mbox{\rm {\small d}}}

\newcommand\ntimes{{\mkern-2mu{\times}\mkern-2mu}}
\newcommand{\ttmes}{\mathop{\mathchoice{\raisebox{-.75ex}{\Huge$\ntimes$}}{\raisebox{-.35ex}{\LARGE$\ntimes$}}{\textstyle\ntimes}{\scriptstyle\ntimes}}}

\noindent {\bf {\Large From Extreme Values of I.I.D.~Random Fields\\
[5pt] to Extreme Eigenvalues of Finite-volume Anderson
Hamiltonian}}

\vspace{.7cm}

\noindent {\bf A.~Astrauskas}\footnote{ VU Institute of
Mathematics
and Informatics, Akademijos str. 4, LT-08663 Vilnius, Lithuania;\\
\mbox{\qquad }{e-mail}: arvydas.astrauskas@mii.vu.lt}

\vspace{0.7cm}

{\footnotesize
{\leftskip=1.6cm

\noindent \hrulefill

\smallskip

\noindent The aim of this paper is to study asymptotic geometric
properties almost surely or/and in probability of extreme order
statistics of an i.i.d.~random field (potential) indexed by sites
of multidimensional lattice cube, the volume of which unboundedly
increases. We discuss the following topics: (I) high level
exceedances, in particular, clustering of exceedances; (II) decay
rate of spacings in comparison with increasing rate of extreme
order statistics; (III) minimum of spacings of successive order
statistics; (IV) asymptotic behavior of values neighboring to
extremes and so on. The conditions of the results are formulated
in terms of regular variation (RV) of the cumulative hazard
function and its inverse. A relationship between RV classes of the
present paper as well as their links to the well-known RV classes
(including domains of attraction of max-stable distributions) are
discussed.

The asymptotic behavior of functionals (I)--(IV) determines the
asymptotic structure of the top eigenvalues and the corresponding
eigenfunctions of the large-volume discrete Schr\" odinger
operators with an i.i.d.~potential (Anderson Hamiltonian). Thus,
another aim of the present paper is to review and comment a recent
progress on the extreme value theory for eigenvalues of random Schr\"
odinger operators as well as to provide a clear and rigorous
understanding of the relationship between the top eigenvalues and
extreme values of i.i.d.~random potentials. We also discuss their
links to the long-time intermittent behavior of the
parabolic problems associated with the Anderson Hamiltonian via spectral
representation of solutions.

\noindent \hrulefill

\smallskip
\noindent \textsf {\textbf{KEY WORDS:}} extreme value theory;
Poisson limit theorems; extreme order statistics; high-level
exceedances; spacings; regular variation; Weibull distribution;
discrete Schr\" odinger operator; Anderson Hamiltonian; random
potential; largest eigenvalues; principal eigenvalues;
localisation; parabolic Anderson model; intermittency.

\medskip

\noindent {\sc msc {\footnotesize 2010} subject classification:
primary -  {\footnotesize 60}G{\footnotesize 70}, {\footnotesize
60}H{\footnotesize 25}, {\footnotesize 82}B{\footnotesize 44},
{\footnotesize 35}P{\footnotesize 05} se\-con\-dary -
{\footnotesize 60}F{\footnotesize 05}, {\footnotesize
60}F{\footnotesize 15}, {\footnotesize 60}G{\footnotesize 60},
{\footnotesize 82}C{\footnotesize 44}, {\footnotesize
35}P{\footnotesize 15}, {\footnotesize 15}B{\footnotesize 52},
{\footnotesize 26}A{\footnotesize 12}}

\par}

}

\tableofcontents

\section*{\normalsize{\textsf {REFERENCES}}\hbox{\hbox to 9,7cm{ \dotfill}\hfill \quad\ 78}}

\newpage

\hspace{1cm}

\section{\normalsize{\textsf {INTRODUCTION}}}

\subsection{\normalsize{\textsf {Extremes of  i.i.d.~random fields } }}

In this paper, we assume that $\xi(x)$, $x\in \ttZ^\nu $, are 
independent identically distributed (i.i.d.) random variables on a
probability space $(\rrO ,\cF ,\bP )$, indexed by sites of the
$\nu $-dimensional integer lattice $\ttZ^\nu$, with a distribution
function $ \bP(\xi(0)\mly t)=:1-\e^{-Q(t)}$, $t\in \ttR$; here $Q$
denotes the cumulative hazard function of distribution. Define
$V=[-n;n]^\nu\cap\ttZ^{\nu}$, the cubes in $\ttZ^\nu$. Let $\ve V
\ve$ denote the number of sites in $V$. We write $ \vert
x\vert=\sum^\nu_{i=1} \ve x^i\ve$ for the lattice $l^1$-distance
between $x=(x^1,\ldots ,x^\nu)\in \ttZ^\nu $ and $0\in \ttZ^\nu $.

We consider the variational series (order statistics) \begin{equation} \xi\iz{1,
V}:=\xi(z\iz{1, V})\dly \xi\iz{2, V}:=\xi(z\iz{2, V})\dly \ldots
\dly \xi\iz{\vert V\vert , V}:=\xi(z\iz{\vert V\vert , V})
\end{equation}
based on the sample $\xi_V :=\{ \xi(x)\colon\ x\in V\}$; here $V =
\big \{ z_{k,V} \colon\ 1 \mly k \mly \ve V \ve \big \}$. The
first $\ve V\ve^\rre$ ($0<\rre <1$) terms of the variational
series (1.1) are referred to as $\xi_V $-extremes or
$\xi_V$-peaks. The coordinate $z_{k,V} \in V$ stands for a
location of the $k$th extreme value of $\xi_V$; $1 \mly k \mly \ve
V\ve$.

In this paper, letting $\ve V\ve\to\infty$, we study the
asymptotic geometric properties of $\xi_V $-extremes almost surely
and/or in probability. We are interested in the following
functionals of order statistics (1.1):

\textbf{(EX)}  {\it Exceedances} of the sample $\xi_V$ over high
levels $L_V$, in particular, clustering of exceedances (Theorem
3.1).

\textbf{(SP)} { \it The decay rate of the spacings}
$\xi\iz{K,V}-\xi\iz{K+1,V}$ and $\xi\iz{[\ve
V\ve^\rre],V}-\xi\iz{[\ve V\ve^\theta],V}$ in comparison with
increasing rate of $\xi\iz{K,V}$ for fixed natural $K\in \ttN$ and
$0\mly \rre< \theta <1$ (Theorems 4.3--4.7).

\textbf{(MIN)} { \it Minimum of the spacings}
$\xi\iz{l,V}-\xi\iz{l+1,V}$, $1\mly l\mly \ve V\ve^\rre$, for each
$0<\rre <1$ (Theorems 4.8 and 4.9).

\textbf{(NEI)} {\it $\xi_V$-values neighboring to $\xi_V$-extremes},
in particular, $\xi(z\iz{l, V}+y)$ for $1\mly l\mly \ve V\ve^\rre$
and for fixed $y\neq 0$ (Lemma 5.1 and Theorems 5.3, 5.4).

The conditions of the asymptotic results for (EX), (SP), (MIN) and
(NEI) are given in terms of regular variation (RV) of the inverse
function of $Q$. In Appendix A and Section 6, we discuss a relationship between
RV classes of the present paper as well as their links to the
well-known RV classes including domains of attraction of
max-stable distributions.

The asymptotic results for (EX), (SP), (MIN) and
(NEI) and related RV classes were announced without the proof in 
(Astrauskas 2007; 2008; 2012; 2013). In this survey,  these results 
are given in the most general setting with the detailed proof; therefore, they 
present self-contained  topics of probability theory and may be considered of independent interest. 

\medskip

\subsection{\normalsize{\textsf {Extreme value theory for eigenvalues of large-volume Anderson Hamilto-\ nians} }}

Let us consider the finite-volume Schr\" odinger operators
$\cH_V=\rrk\Delta_V +\xi_V$ on $l^2(V)$ with periodic boundary
conditions (Anderson Hamiltonian); here $\rrk >0$ is a diffusion
constant; $\Delta \psi(x):=\sum _{\vert y-x\vert =1}\psi(y)$ is the
lattice Laplacian, and the i.i.d.~random field $\xi_V :=\{
\xi(x)\colon\ x\in V\}$ is the multiplication operator
(potential). Denote by $\rrl\iz{K,V}$ the $K$th largest eigenvalue
of the operators $\cH_V$, and let $\psi(x;\rrl\iz{K,V})$ ($x\in
V$) be the corresponding eigenfunction normalized to have unit
$l^2$-norm, $\sum_{x\in V}\psi(x;\rrl\iz{K,V})^2=1$. Another aim
of this paper is to show in what manner the asymptotic behavior of
functionals (EX), (SP), (MIN) and (NEI) determines the asymptotic
structure of the top eigenvalues $\rrl\iz{K,V}$ and the
corresponding eigenfunctions, as $V\uparrow \ttZ^\nu$ and $K \dly
1$ fixed. In Section 2, we give an overview of rigorous statements
on this relationship which are proved in the papers by
G\"artner and Molchanov (1998) and Astrauskas~(2007; 2008; 2012).
In Section~6, we review and comment results on the asymptotic expansion
formulas and Poisson limit theorems for the largest eigenvalues
$\rrl\iz{K,V}$ as well as localization properties of the
corresponding eigenfunctions. These results are proved by
Astrauskas and Molchanov~(1992), G\"artner and Molchanov (1998),
Astrauskas~(2007; 2008; 2012; 2013), Germinet and Klopp (2013),
Biskup and K\"onig (2016) and other mathematicians. These papers
are complemented by the present survey on the asymptotic geometric
properties of $\xi_V$-extremes and related RV classes of
distributions. We here give proof sketches of the results on the extreme value theory 
for eigenvalues $\rrl_{K,V}$ (Sections 2 and 6) demonstrating 
their connections to asymptotic properties of $\xi_V$-extremes.

Thus, in this survey, we discuss in detail the following important
branches of probability theory: (i) extreme value theory for
eigenvalues of the Anderson Hamiltonian $\cH_V=\rrk\Delta_V
+\xi_V$ which is a particular model of random matrices (Sections
1.4, 2 and 6); (ii) asymptotic geometric properties of random
i.i.d.~fields (Sections 3--5) and (iii) regular variation of
distribution functions (Appendix A). On the other hand, we briefly
comment the links of the extreme value theory for eigenvalues to the
following important topics of statistical physics: (iv) Anderson
localization for the random Schr\" odinger operators
$\cH=\rrk\Delta +\xi(\cdot)$ in the whole lattice $\ttZ^\nu$
(Section~1.3), and (v)~long-time intermittent behavior of
solutions $u$ of the parabolic problems associated with the
Anderson Hamiltonian (PAM) via spectral representation of $u$
(Section~7).

Of course, asymptotic results for the Anderson models
(time-dependent or -indepen-\ dent) are heavily necessitated by the
asymptotic structure of high $\xi(\cdot)$-values, which in turn is
determined by conditions on the regularity and tail decay of the
distribution $ \bP(\xi(0) > t)=\e^{-Q(t)}$ at its right endpoint
$t_Q:={\it esssup}\ \xi(0)$. In the present survey, we focus on
the case of unbounded from above i.i.d.~potentials, i.e.
$t_Q:=\infty$, with distributional tails heavier than double
exponential, i.e., $\bP(\xi(0) > t)=\exp \{-\e^{\mo (t)}\}$ as
$t\to \infty$; cf.~Sections 1.2, 2.2--2.3 and 6.1. For such
distributions satisfying additional RV and continuity conditions,
we will show that with probability one the $\xi_V$-peaks are
spatially separated and differ in height as $V\uparrow \ttZ^\nu$;
therefore, the top eigenvalue $\rrl\iz{K,V}$ of the operator
$\cH_V=\rrk\Delta_V +\xi_V$ is approximated by an isolated
$\xi_V$-peak, say $\xi(z\iz{\tau(K),V})$, plus some corrections of
order $\mro(1)$ involving neighboring $\xi_V$-values. Moreover,
the $K$th eigenfunction $\psi(\cdotp; \rrl_{K,V})$ is
asymptotically delta like function at the site $z\iz{\tau(K),V}
\in V$, the localization center. (In this case, we will say that
the $K$th eigenvalue is associated with the site
$z\iz{\tau(K),V}$, viz. $\rrl_{K,V}\leftrightarrow
z\iz{\tau(K),V}$). Therefore, we are able to apply the standard
extreme value theory to prove Poisson limit theorems for the
normalized extreme eigenvalues and their localization centers.
From these Poisson limit theorems one obtains the limiting joint
(max-stable) distribution for the normalized largest eigenvalues
and their spacings, limiting uniform distribution for the
normalized localization centers and other important limiting
distributions for eigenvalue statistics (Section 6.1). Eigenvalue
statistics in turn play a crucial role in studying the intermittent
behavior of the parabolic Anderson model, PAM (Section~7).

For the lighter upper tails including the double exponential $\bP(\xi(0)
> t)=\exp \{-\e^t\}$ and bounded tails i.e. $t_Q < \infty$, we will prove
the rough asymptotic expansion formulas for the largest eigenvalues. For such distributional
tails, it will turn out that all $\xi_V$-extremes are of comparable amplitude; therefore,
the $K$th largest eigenvalue is associated with a large island of higher
$\xi_V$-values of a particular preferred shape, rather than an isolated $\xi_V$-peak.

To illustrate the relationship between
$\xi_V$-extremes and the largest eigenvalues $\rrl_{K,V}$
($V\uparrow \ttZ^\nu$) more precisely, we now formulate Propositions 1.1--1.3
which are "typical" examples of the statements given in Sections
2--6. The first proposition tells us that, if the peaks of {\it
deterministic (nonrandom) functions} $\xi_V=:\xi_V(\cdotp)$ are
extremely sharp and widely spaced, then the $K$th largest
eigenvalue $\rrl\iz{K,V}$ is approximated by the $K$th largest
value of $\xi_V$ with sufficiently small error.

\begin{propbt}{\rm (see Theorem 2.2(ii) in Section 2.2).}
Fix constants $K\in \ttN$ and $0<\theta <1/2$, and assume that the
deterministic functions $\xi_V$ satisfy the following conditions
as $V\uparrow \ttZ^\nu$: \beq \min_{1 \mly l \mly
K}\xi_{l+1,V}(\xi\iz{l,V}-\xi\iz{l+1,V}) \to \infty \quad
\hbox{(distinct height of peaks)},
\eeq \beq \frac{1}{\log\ve V\ve}\min_{1\mly k<n \mly\ve
V\ve^{\theta}} \big|z\iz{k, V}-z\iz{n, V}\big| \to \infty \quad
\hbox{(sparseness of peaks)}
\eeq and, finally, \beq \xi\iz{[\ve V\ve^\theta],V}\break
/\xi\iz{K,V} < \const(\theta)
\eeq for some $0<\const(\theta)<1$ ({\it negligibility of the lower peaks}). Then
$$
\rrl\iz{l,V}=\xi\iz{l,V}+\rO(1/\xi\iz{l,V}) \quad \hbox{for
all}\quad 1\mly l \mly K.
$$
\end{propbt}

We now give an example of i.i.d.~random field $\xi(\cdot)$ with sufficiently
\lk heavy tails\rk\ possessing extremes like those in Proposition
1.1.

\begin{propbt}{\rm (see Theorem 4.3(i) with $p=1$, Theorem 3.1
with $R=0$ and Theorem 4.5).} If $\xi(0)$ has the Weibull distribution
\beq \bP(\xi(0)> t)=\e^{-Q(t)}=\e^{-t^\alpha} \quad
(t\dly 0)
\eeq
with $\rra <2$, then the i.i.d.~sample $\xi_V$ ($V\uparrow \ttZ^\nu$) satisfies {\rm (1.2)--(1.4)}
with probability $1+\mro(1)$.
\end{propbt}

Since, by Propositions 1.1 and 1.2, the eigenvalues $\rrl\iz{K,V}$
are very close to $\xi\iz{K,V}$ as $V\uparrow \ttZ^\nu$, it turns
out that Poisson limit theorems (and the corresponding
renormalization constants) for the largest eigenvalues are the
same as that for $\xi_V$-extremes according to the following
proposition.

\begin{propbt}{\rm (see Theorem 6.9 and (Astrauskas 2012)).} Assume that
$Q(t)=t^\alpha $ with $\rra <2$, and write $b_V:=(\log \ve
V\ve)^{1/\rra}$. Define the point process $\cN_V^\rrl$ on
$[-1/2;1/2]^\nu \times \ttR$ by
$$
\cN_V^\rrl :=\sum_{k=1}^{\ve V\ve} \delta_{\Lambda \iz{V}(k)}
\quad \hbox {with} \quad \Lambda \iz{V}(k):=\bigg(\frac
{z\iz{k,V}}{|V|^{1/\nu}},\frac{\rrl\iz{k,V}-b_V}{\rra^{-1} b_V^{1-\rra}} \bigg),
$$
where $\delta_X$ denotes the Dirac measure at $X \in [-1/2;1/2]^\nu
\times \ttR$. Then $\cN_V^\rrl$ converges weakly to the Poisson
process on $[-1/2;1/2]^\nu \times \ttR$ with the intensity measure
$\dr x\times \e^{-t} \dr t$.

Moreover, for fixed $K \dly 1$, the eigenfunction $\psi(\cdot;\rrl\iz{K,V})$
is exponentially localised at the site $z\iz{K,V}$:
\beq \lsup_V\max_{x \neq z_{K,V}}\frac
{\log\big|\psi(x;\rrl\iz{K,V})\big|}
{|x-z\iz{K,V}|\log b_V} \mly -1 \eeq
in probability.
\end{propbt}

According to Proposition 1.3 the $K$th largest eigenvalue
$\rrl_{K,V}$ is associated with the $K$th largest value of
$\xi_V$, viz., $\rrl_{K,V}\leftrightarrow z_{K,V}$. For the
lighter tails, say, Weibull distributions (1.5) with $\rra \dly 2$, the
landscape of $\xi_V$ gets \lk smoother\rk, in particular, (1.2)
fails. Therefore, $\rrl_{K,V}$ is associated with a lower and \lk
slightly supported\rk \ $\xi_V$-peak, viz.,
$\rrl_{K,V}\leftrightarrow z_{\tau(K),V}$, where for $\rra > 3$,
the index $\tau(K)=\tau_V(K)$ tends to infinity as $\ve V\ve \to
\infty$.  This in turn implies that further terms in expansion for
$\rrl\iz{K,V}$ become essential; see (2.22) and Examples
6.12--6.13. Let us distinguish three classes (J)--(JJJ) of light
tailed distributions (i.e., {\it universality classes}),
which ensure a different asymptotic behavior of
the eigenvalues $\rrl_{K,V}$.

\textbf{(J)} {\it Distribution tails heavier than the double
exponential function.} Assume that \beq \log Q(t)=\mro(t)
\eeq and $Q$ satisfies additional regularity and continuity
conditions as $t\to \infty$. This class is presented by Weibull distributions
(1.5) for arbitrary $\rra >0$ and those with fractional-double-exponential tails
\beq \bP(\xi(0)> t)=\e^{-Q(t)}=\exp \{-\e^{t^\rrg}\} \quad (t \dly t_0)
\eeq
for $\rrg < 1$. For such distributions, $\xi_V$-extremes possess a {\it
strongly pronounced geometric structure} which can be described as
follows:

For arbitrary sufficiently small constants $0 < \rre <\theta$,
there exist constants $c_1 >c_2 >0$ and (large) $C>0$ such that
almost surely \beq \min_{1\mly l<n \mly\ve
V\ve^{\theta}}\big(\xi\iz{l,V}-\xi\iz{n,V}\big) \dly \e^{-\ve
V\ve^{c_2}} \quad \hbox{(distinct height of peaks)},
\eeq \beq \min_{1\mly l<n \mly\ve V\ve^{\theta}} \big|z\iz{l,
V}-z\iz{n, V}\big| \dly \ve V\ve^{c_1} \quad \hbox{(sparseness of
peaks)}
\eeq and, finally, \beq \xi\iz{[\ve V\ve^\rre],V}-\xi\iz{[\ve
V\ve^\theta],V}\dly C \quad \hbox{(negligibility of the lower peaks)}
\eeq for each large $V$; see Theorem 4.8 with $\rrk =0$, Theorem
3.1 with $R=0$ and Theorem~4.6 with $\rho= \infty$. By the
standard finite-rank perturbation arguments in (Astrauskas and
Molchanov 1992) and (Astrauskas 2008), these properties of $\xi_V$
yield that there is no resonance between $\xi\iz{V}$-peaks in the
Anderson model for large $V$; therefore, the eigenvalues
associated with a block of peaks can be determined by the local
eigenvalues associated with separate peaks (i.e.,{\it \lk relevant
single peak\rk approximation}). More precisely, for fixed natural
$K$ and $V\uparrow \ttZ^\nu$, almost surely the eigenvalue
$\rrl_{K,V}$ of $\cH_V=\rrk\Delta_V +\xi_V$ is approximated by the
principal (i.e., the first largest) eigenvalue of the \lk single
peak\rk Hamiltonian $\rrk \Delta \iz{V}+\wdt \xi(\cdotp
)+\xi(z_{\tau(K),V})\delta _{z_{\tau(K),V}}$ where $\log
\tau(K)=\mro(\log \ve V\ve)$. Here $\wdt \xi(\cdotp )$ is the \lk
noise\rk\ potential; the site $z_{\tau(K),V}\in V$ is a
localization center of the $K$th eigenfunction
$\psi(\cdot;\rrl\iz{K,V})$ of $\cH_V$. Thus, Poisson limit
theorems for the eigenvalues $\rrl_{K,V}$ of $\cH\iz{V}$ are
reduced to those for the principal eigenvalues of the \lk single
peak\rk \ Hamiltonians, which in turn are expanded into certain
(nonlinear) series in $\xi(x)$ ($x\in V$); cf.
formulas~(2.20)--(2.22), Theorems~2.3,~6.2 and discussions in
Section~6.4. We finally notice that the $K$th eigenfunction
$\psi(\cdot;\rrl\iz{K,V})$ is exponentially well localized, i.e., there exist non-random constants (decay rates)
$C>0$ and $0<M_V \to \infty$ such that with probability one \beq \vert \psi(x; \rrl
\iz{K,V})\vert \mly C \exp \{ -M_V\vert x-z\iz{\tau (K), V}\vert
\}\quad (x\in V)
\eeq for all $V$ large enough. Consequently,
$\psi(\cdot;\rrl\iz{K,V})$ is asymptotically delta-like function
at $z_{\tau(K),V}$ (Astrauskas 2008; 2013). This refers to the
correspondence $\rrl_{K,V}\leftrightarrow z_{\tau(K),V}$. Under
assumption (1.7), asymptotic expansion formulas, Poisson limit
theorems and localization theorems for the largest eigenvalues are
derived by Astrauskas and Molchanov (1992), Astrauskas (2007;
2008; 2012; 2013). See also Grenkova et al.~(1983) and Grenkova et
al.~(1990) for the case of Weibull distribution (1.5) with $\rra <
2$.

\textbf{(JJ)} {\it Distribution tails lighter than the double
exponential function.} Assume that \beq t^{-1}\log Q(t) \to \infty
\eeq and $Q$ satisfies additional regularity conditions as $t$
tends to $t_Q$ (= the right endpoint of $Q$). This class of
potentials contains the important case of $\xi(\cdot)$ which is
bounded from above ($t_Q <\infty$) and those with
fractional-double-exponential tails (1.8) for $\rrg > 1$. For such
$\xi(\cdot)$, it turns out that $\xi_V$-peaks possess a {\it
weakly pronounced geometric structure}. In particular, almost
surely $\xi\iz{[\ve V\ve^\rre],V}-\xi\iz{[\ve V\ve^\theta],V} \to
0$ as $\ve V\ve \to \infty$, for all $0 \mly \rre <\theta <1$, so
that the height of all $\xi_V$-extremes is of the same order
$\xi_{1,V} +\mro(1)$ (see Theorem~4.6 with $\rho =0$ and Theorem
3.1(i) with arbitrary $R \dly 1$ and $\theta(\cdot)\equiv
\theta=\const$).  In this case, the eigenvalue $\rrl \iz{K,V}$ ($K
\dly 1$ fixed) does not longer correspond to an isolated potential
peak, but to a flat extremely large \lk relevant island\rk \ of
high $\xi_V$-values. More precisely, almost surely the top
eigenvalue $\rrl\iz{K,V}$ of the Hamiltonian $\cH_V$ approaches
the local principal eigenvalue of the Hamiltonian restricted to a random
connected region $\ttA_{V; \rm opt}^K \subset V$ with the
following features: The diameter of $\ttA_{V; \rm opt}^K $
unboundedly increases, and $\xi(\cdot)$ possesses in $\ttA_{V; \rm
opt}^K $ values of the order $\xi_{1,V} +\mro(1)$ as $V\uparrow
\ttZ^\nu$, i.e., relevant island of potential values. Moreover,
the $K$th eigenfunction is expected to be highly concentrated in
the neighborhood of the region $\ttA_{V; \rm opt}^K$. In this
case, we will say that the $K$th eigenvalue is associated with
$\ttA_{V; \rm opt}^K$, viz. $\rrl_{K,V}\leftrightarrow \ttA_{V;
\rm opt}^K$. For more explanations, see Theorem~6.14 and the
proof of Theorem 2.6, where the second order expansion formula for
$\rrl \iz{K,V}$ is obtained.

For the Bernoulli i.i.d.~random variables $\xi(x)$ with $t_Q=1$
(which is the particular case of (1.13)), Bishop and Wehr (2012)
derived a more accurate expansion formula for the principal
eigenvalue $\rrl \iz{1,V}$ of the one-dimensional Hamiltonian
$\cH_V$ ($\nu =1$). They have showed that $\rrl \iz{1,V}$ is
associated with the longest consecutive sequence of sites $x \in
V$ with $\xi(x)=1$, i.e., the \lk relevant island \rk of
$\xi_V$-extremes, the length of which unboundedly increases as $V
\uparrow \ttZ$. See also (Sznitman 1998) for similar asymptotic
results in the case of spatially continuous Schr\" odinger
operators with a bounded Poisson potential of obstacles. Cf.~Section
6.2 below.

To the best of our knowledge, Poisson limit theorems for the
(unfolded) largest eigenvalues were proved only in the case
$\nu=1$ and bounded $\xi(0)$, provided the distribution
$1-\e^{-Q}$ satisfies additional continuity and tail decay
conditions (Germinet and Klopp 2013); see also Section 6.2 below.
For $\nu \dly 2$ or general RV conditions on $Q$ satisfying
(1.13), the Poissonian convergence of the top eigenvalues still
remains an open problem.

\textbf{(JJJ)} {\it Double exponential type tails.} Finally,
assume that $t^{-1}\log Q(t)$ tends to a positive finite constant $ \rrr^{-1}$
as $t \to \infty$, i.e., the double exponential
tails
\beq \bP(\xi(0)> t)=\e^{-Q(t)}=\exp \{-\e^{t(\rrr^{-1}+\mo(1))}\}
\eeq satisfying additional RV and continuity conditions. This
class of distributions presents the intermediate case between (J)
and (JJ). For such $\e^{-Q}$, it turns out that  all
$\xi_V$-extremes are of comparable amplitude; i.e., almost surely
$\xi\iz{[\ve V\ve^\rre],V}-\xi\iz{[\ve V\ve^\theta],V}=\rO(1)$ as
$\ve V\ve \to \infty$, for any $0 \mly \rre <\theta <1$.
Therefore, with probability one the top eigenvalue $\rrl\iz{K,V}$
($K \dly 1$ fixed) of the Hamiltonian $\cH_V$ is approximated by
the local principal eigenvalue of the Hamiltonian restricted to a random
connected region $\ttA_{V; \rm opt}^K \subset V$ of bounded
diameter, where $\xi_V$ possesses high values of the optimal shape (so
that  $\rrl_{K,V}\leftrightarrow \ttA_{V; \rm opt}^K$). The
optimal shape of $\xi_V$-values in $\ttA_{V; \rm opt}^K \subset V$
is specified by deterministic variational principles. These
considerations are referred to as the {\it \lk relevant island\rk\
approximation}; see Theorems 2.7 and 6.19 for the second order
expansion formulas for the first largest eigenvalue $\rrl
\iz{1,V}$, which have been originally derived by G\"artner and
Molchanov~(1998).  Moreover, the $K$th eigenfunction $\psi(\cdotp;
\rrl_{K,V})$ is highly concentrated in the neighborhood of the
region $\ttA_{V; \rm opt}^K$, as proved by Astrauskas (2008; 2013)
for $\rrr$ large enough, and by Biskup and K\"onig (2016) for
arbitrary $\rrr$; see also Section 6.3 of the present survey.

Rigorous results on Poisson limit theorems and further
localization properties for the largest eigenvalues (in the case
of double exponential tails) have been proved by Astrauskas
(2007; 2008; 2013) for $\rrr$ large enough, and by Biskup and
K\"onig (2016) for arbitrary $\rrr$; see also the review paper by
K\"onig (2016) and Sections 6.3--6.4 of the present survey for the
discussions on their results and the proofs.

Let us finally summarize the above observations: As the upper
tails of potential distribution get lighter, the $\xi_V$-extremes
($V\uparrow \ttZ^\nu$) get less expressed; therefore, the number
of higher $\xi_V$-values contributing to the asymptotic amount of the
top eigenvalues gets larger and concentration properties of the
corresponding eigenfunctions become weaker. On the other hand, the
general theory of Anderson localization (cf.~Section 1.3 below)
suggests that almost surely the eigenfunctions associated with the
upper spectral edge of the Hamiltonian $\cH_V$ decay exponentially
like in (1.17), for arbitrary potential distribution satisfying
certain continuity conditions.

\medskip

\subsection{\normalsize{\textsf {Relations to infinite-volume Anderson Hamiltonians} }}

{\bf (I)} {\it Anderson localization}. The Anderson model on the
whole lattice $\ttZ^\nu$ is given by the Hamiltonian
$$
\cH=\rrk\Delta +\xi(\cdotp )
$$
acting on $l^2(\ttZ^\nu )$. Here, as above, $\Delta $ is the
lattice Laplacian, $\rrk >0$ is a diffusion constant, and $\xi(x)$
($x\in \ttZ^\nu $) are i.i.d.~random variables  with a common
distribution function $1-\e^{-Q}$. This is a basic model of
disordered quantum systems introduced to describe the regions of
energy levels (spectrum) of the electron in the random potential
modelling electrical conductance regimes of alloys, crystals with
impurities and so on (Anderson 1958). The energy spectrum $\spect
(\cH)$ of the Hamiltonian $\cH=\rrk\Delta +\xi(\cdotp )$ is almost
surely nonrandom:
$$
\spect (\cH)=\spect(\rrk\Delta)+\spect (\xi(\cdotp ))=[-2\nu \rrk;2\nu \rrk]+\supp(1-\\e^{-Q}),
$$
where $\supp(F)$ is the support of a probability measure generated
by the distribution function $F$, and "+" denotes the algebraic
sum of subsets of real line. Therefore, with probability one, the
spectrum consists of spectral bands situated in the interval
$[L_{\min};L_{\max}]$, where $L_{\min}$ and $L_{\max}$ are
respectively the infimum (i.e.~bottom) and the supremum (i.e.~upper edge) of
the spectrum. The most important property of the Hamiltonian $\cH$
on $l^2(\ttZ^\nu)$ is the presence of pure point spectrum in the
neighborhood of edges of spectral bands for any $\rrk>0$ and any
$\nu\dly 1$. In particular, there exist (nonrandom) real constants
$L_i=L_i(\rrk, \nu, Q)$, $L_1 <L_2$, such that with probability
one
$$
\hbox {the\ spectrum\ in}\ \ (L_{\min}; L_1)\cup(L_2; L_{\max} )\
\ \hbox {is\ dense\ purely\ point},
$$
say $\{\rrl_k\}$, and the corresponding eigenfunctions
$\psi(\cdot;\rrl\iz{k})$ decay exponentially:
$$
\ve \psi(x;\rrl\iz{k})\ve \mly C_k\exp \{ -M\vert x-z\iz{k}\vert
\}\quad (x\in \ttZ^\nu)
$$
for some random $C_k >0$, $M> 0$ and $z\iz{k}\in \ttZ^\nu$ (the
localization center), provided $\e^{-Q}$ is H\" older continuous
and $\xi(0)$ has some finite statistical moments. Moreover, for
small $\rrk $ or the one-dimensional case $\nu=1$, the whole
spectrum $\spect(\cH )$ is dense purely point with a complete set
of eigenfunctions in $l^2(\ttZ^\nu)$ that decay exponentially with
probability~1. This phenomenon is known as \textit {Anderson
localization} for disordered systems; see, e.g., (Fr\" ohlich and
Spencer 1983; Fr\" ohlich et al.~1985; Simon and Wolff 1986;
Carmona et al.~1987; Aizenman and Molchanov 1993; Aizenman et al.~2001) for the proof of the above assertions under various
conditions on potential distributions. Recall that, for the
periodic $\xi(\cdot)$, all the spectrum $\spect (\cH)$ is
absolutely continuous in arbitrary dimension $\nu \dly 1$; and
this is quite a contrast to the Anderson localization in the case
of random potential. See the monographs (Pastur and Figotin 1992; 
Stolz 2011; Kirsch 2008) and references therein for more
discussions on the subject. In the present survey as well as in
(Astrauskas 2007; 2008; 2012; 2013) and (Biskup and K\"onig 2016),
the phenomenon of Anderson localization is illustrated for the top
eigenvalues and eigenfunctions of finite-volume models. We here
emphasize the relationship between asymptotic geometric properties
of $\xi_V$-extremes and localization properties of the leading
eigenfunctions of the operator $\cH_V$ regarding their
localization strength, localization centers, etc., as $V\uparrow
\ttZ^\nu$.

\medskip

Let us discuss  briefly the basic ideas and methods explored in
the study of the Anderson localization phenomenon in the
multidimensional case $\nu \dly 1$. Fix an open bounded interval
$I\subset \ttR$ which is covered almost surely by the spectrum
$\spect(\cH)$. The {\it proof of Anderson localization} in the
spectral intervals $I$ relies heavily on the study of the
resolvents $\cG(\rrl+i\rre):=(\rrl+i\rre-\cH)^{-1}$ and the
corresponding Green functions $\cG(\rrl+i\rre; x,
y):=\cG(\rrl+i\rre)\delta _y(x)$ ($x,y \in \ttZ^{\nu}$) in the
complex domain $\rre >0$, $\rrl \in I$, where $i=\sqrt{-1}$.
Alternatively, their finite-volume versions
$\cG\iz{V}(\rrl+i\rre):=(\rrl+i\rre-\cH\iz{V})^{-1}$ and
$\cG\iz{V}(\rrl+i\rre; x, y):=\cG\iz{V}(\rrl+i\rre)\delta _y(x)$
($x,y \in V$) are also explored. The main task here is to prove
that, for $\rrl \in I$, the Green functions $\cG(\rrl+i\rre; x,
y)$ or $\cG\iz{V}(\rrl+i\rre; x, y)$ decay exponentially in $\ve
x-y\ve$ uniformly in $\rre>0$, provided $I$ is chosen close to the
spectral edge or the diffusion constant $\rrk$ is small. For
simplicity, assume throughout that the potential distribution has
a bounded density $p(\cdot)$ with bounded support.

Fr\" ohlich and Spencer (1983) developed the {\it multiscale
method}. They constructed inductively a sequence of relevant cubes
$V^{\prime}$, $V^{\prime}\uparrow \ttZ^\nu$ with the following
properties: For any fixed $\rrl \in I$ and $V^{\prime}\uparrow
\ttZ^\nu$, with high probability the Green function
$\cG\iz{V^{\prime}}(\rrl+i\rre; x, y)$ decays exponentially in $\ve
x-y\ve$ for all $y \in \partial V^{\prime}$, the \lk boundary\rk of $V^{\prime}$, and all
$x\in V^{\prime}$ far away from $\partial V^{\prime}$, and this estimate holds uniformly
in $\rre \dly 0$. By finite-rank perturbation formulas, this estimate
implies an absence of absolutely continuous spectrum in $I$ with
probability one. Even the stronger form of the construction of the
relevant cubes $V^{\prime}$ (\lk uniformity\rk\ in $\rrl \in I$) is
applied to prove the exponential decay of the eigenfunctions
$\psi(\cdot;\rrl)$ associated with (generalized) spectral values
in $I$ and, consequently, the presence of the pure point spectrum
in $I$ with probability one (Fr\" ohlich et al.~1985). See also
the survey by Kirsch (2008) for a detailed discussion on the
multiscale analysis.

In the above considerations, one should apply the {\it Wegner
estimate}. Loosely speaking, this estimate states that the mean
number of eigenvalues of $\cH_V$ in the interval $I$ does not
exceed $\ve I\ve \ve V\ve C$, provided the density $p(\cdot)$ of
potential distribution satisfies $p(\cdot)\mly C$. In particular,
the latter guarantees the bound of probability to find at least
one eigenvalue $\rrl\iz{l,V}$ in a small spectral interval. The
Wegner estimate and its modifications are the basic probabilistic
tools in the proof of Anderson localization; see, e.g., (Kirsch
2008).

Simon and Wolff (1986) applied the so-called spectral averaging
methods in multiscale analysis to obtain the effective condition
for the Anderson localization in the interval $I$, for the
distributional density as above. Recall the {\it Simon-Wolff
criterion}: If for each $x\in \ttZ^\nu$ and for Lebesgue-almost
every $\rrl\in I$ with probability one \beq \lim_{\rre\downarrow
0}\sum_{y\in \ttZ^\nu}\lf\ve \cG(\rrl+i\rre; x, y)\rg\ve^2
<\infty,
\eeq then the operator $\cH$ has only pure point spectrum in $I$
with probability one. If in (1.15) one claims a square-summability of the Green function with weights $\e^{m \ve x-y\ve}$ for some nonrandom
$m>0$, then the corresponding eigenfunctions decay exponentially
with probability one.

Aizenman and Molchanov (1993) and Aizenman et al.~(2001) developed
the {\it fractio-\ nal moment method} to prove the Anderson
localization in the spectral intervals $I$. The task here is to
obtain (and explore) the exponential decay of the averaged Green
functions, where the average is taken over the random potential.
This method enables to avoid a complicated dependency of the
previous constructions on individual potential configurations in
the almost sure setting. Under certain continuity conditions on
potential distribution, the key statement in (Aizenman and
Molchanov 1993) is the following {\it fractional-moment criterion}:
If for fixed $0<\sigma <1$, there are constants $C_1 >0$, $M_1>0$
such that \beq \bE \lf (\ve \cG(\rrl+i\rre; x, y)\ve^\sigma \rg)
\mly C_1\e ^{ -M_1\vert x-y\vert }\quad \hbox{for\ all}\quad x,\
y\in \ttZ^\nu,
\eeq for all $\rrl \in I$ and uniformly in $\rre >0$, then one has, with probability one, the Anderson localization in the interval $I$ for the operator $\cH$.

This implication can be proved by using the Simon-Wolff criterion
like in (1.15). On the other hand, $\sigma < 1$ is chosen to
depend only on continuity assumptions for potential distribution.
In particular, the H\" older continuity implies that the left-hand
side of (1.16) is finite. 

Eq. (1.16) can be proved by using
suitable finite rank perturbation arguments, i.e., Krein formulas.
These formulas imply that the kernel $\cG(\rrl+i\rre; x, y)$ is equal to a
simple rational function in the variable $\xi(x)$ with
coefficients depending on potential values outside $x$. Now one
can estimate the left-hand side of (1.16) as an integral of
rational function. These estimates are shown to form a certain
iteration procedure, from which one deduces (1.16). When applying
the iteration scheme, a crucial fact is the assumption that the
diffusion constant $\rrk$ is small or the interval $I$ is near the
upper edge of the spectrum.

Applying similar arguments, Aizenman et al.~(2001) deduced a {\it
fractional-moment finite-volume criteria} for the Anderson
localization in the interval $I$. Roughly speaking, these criteria state that, if for some
$0<\sigma <1$ and some $V$, the expectation 
\begin{equation*}
\begin{split}
&\bE \lf (\ve
\cG_V(\rrl+i\rre; 0, y)\ve^\sigma \rg)\ \ \mbox{is}\ \ \mbox{sufficiently}\ \ \mbox{small}\\
&\mbox{for}\ \ \mbox{all}\ \ y \in \partial V,\ \mbox{uniformly}\ \ \mbox{in}\ \ (\rrl;\rre) \in
I\times \ttR_+, \mbox{\hspace{2.5cm}}
\end{split}
\end{equation*}
then the exponential decay (1.16) holds true in
$I$. The finite-volume criteria and their implications are shown to
hold under continuity assumptions on potential distribution mentioned
above.

The criteria of (Aizenman et al.~2001) also imply the
exponential decay of eigenfunctions of the operators $\cH_V$,
associated with the eigenvalues in (compact) spectral intervals $I$ of Anderson localization,
in particular, for $I$ at the upper edge of the spectrum
$\spect(\cH)$. Under the continuity conditions on potential
distribution as above, we are able to formulate this result in a precise
form: There are nonrandom constants $c>0$ and $M>0$ such that with
probability one \beq \vert \psi(x; \rrl \iz{k,V})\vert \mly \ve
V\ve^c \exp \{ -M\vert x-z\iz{\tau (k), V}\vert \}\quad (x\in V)
\eeq for some $z\iz{\tau (k), V}\in V$ (localization center), for all
$\rrl\iz{k,V} \in I$ and all $V$ large enough; cf.~(Klopp
2011). See also the surveys (Hundertmark
2008; Stolz 2011) for detailed discussions on the
fractional-moment methods and their applications.

\medskip

{\bf (II)} {\it Local fluctuations of eigenvalues in the spectral
regions of Anderson localization}. We denote by $I_{{\rm
pp}}:=(a;b)$ an open interval of real axis such that a certain 
fractional-moment finite-volume criterion is fulfilled in 
$I_{{\rm pp}}$. Therefore, the spectrum in
the whole of $I_{{\rm pp}}$ is purely point, so $\spect (\cH)\cap
I_{{\rm pp}} \subset \spect_{{\rm pp}} (\cH)$ with probability one 
(i.e., Anderson localization in $I_{{\rm pp}}$). The spectral
intervals $I_{{\rm pp}}$ are distinguished by the {\it Poissonian
asymptotic behavior} of eigenvalues $\rrl\iz{l,V}$, close to a
fixed $\rrl^0\in I_{{\rm pp}}$, of the finite-volume model
$\cH\iz{V}$ as $V\uparrow \ttZ^\nu$. Here and to the end of this
section the potential distribution is again assumed to have a bounded smooth
density with bounded support. These limit theorems are formulated
in terms of {\it the integrated density of states}, viz.
$N(\rrl)$, and {\it the density of states}, viz.
$n(\rrl):=N^{\prime}(\rrl)$ ($\rrl \in \ttR^\nu$). Recall that
$N(\cdot)$ is the nonrandom distribution function of eigenvalues
defined as the almost sure limit of the empirical distribution
function $N_V(\rrl):=\# \{k \colon \ \rrl \iz{k,V} \mly \rrl \}/
\ve V\ve$ as $\ve V\ve \to \infty$. Moreover, the Wegner estimate
implies that $N(\cdot)$ is an absolutely continuous function with bounded density $n(\cdot)$, provided the potential distribution
has a bounded density. The support of $n(\cdot)$ coincides almost
surely with the spectrum $\spect(\cH)$; therefore,
$N(\rrl)\rightarrow 0$ (resp., $N(\rrl)\rightarrow 1$) as $\rrl$
approaches the bottom (resp., the upper edge) of the spectrum.
See, e.g., (Kirsch 2008) for the definition of functions
$N(\cdot)$, $n(\cdot)$ and their properties.

Now pick a number $\rrl^0$ from the interval $I_{{\rm pp}}$ such that $n(\rrl^0)> 0$. We
consider the normalized eigenvalues \beq \Lambda_{k,V}^0:=\ve V\ve
n(\rrl^0)\lf( \rrl_{k,V}-\rrl^0 \rg) \quad (1 \mly k \mly \ve
V\ve),
\eeq and define the corresponding point process $\cM_V^0$ on
$\ttR$ by
$$
\cM_V^0 :=\sum_{k=1}^{\ve V\ve} \delta_{\Lambda_{k,V}^0}.
$$
With these assumptions and abbreviations, one needs to show that
the point process $\cM_V^0$ converges weakly, as $V\uparrow
\ttZ^\nu$, to the Poisson point process on $\ttR$ with intensity
measure $\dr \rrl$, i.e.~the Lebesgue measure. The first result in
this direction was proved by Molchanov (1981), who considered the
one-dimensional spatially continuous random Schr\" odinger
operators. In the case of the Anderson Hamiltonians in $\ttZ^\nu$
with arbitrary $\nu \dly 1$, this Poisson limit theorem was
established by Minami (1996). Killip and Nakano (2007) proved the
Poisson convergence of both the normalized spectral values (1.18)
and localization centers of the corresponding eigenfunctions,
extending Minami's result. The proof in the multidimensional case
$\nu \dly 1$ relies on applications of the fractional-moment
criteria for the Anderson localization, combined with the Wegner and
Minami estimates that control the probability of finding, respectively, 
at least one and two eigenvalues of $\cH_V$ in a small interval. In particular, these
estimates imply the upper and lower bounds of order $\ve
V\ve^{-1}$ for the gap between successive eigenvalues close to a
fixed $\rrl^0$ as above. From Poisson limit theorems for (1.18),
one can extract some information on the statistical properties of
eigenvalues in the spectral regions of Anderson localization. For
example, one can obtain the limiting distribution for the
normalized spacings of eigenvalues, the limiting joint
distribution for the normalized eigenvalue and its localization
center, and other important limiting distributions for eigenvalue
statistics.

Recently, Germinet and Klopp (2013; 2014) presented a
comprehensive study on limit theorems for eigenvalues of
large-volume Hamiltonians in the spectral regions $I_{{\rm pp}}$
of Anderson localization.  For fixed $\rrl^0 \in I_{{\rm pp}}$ as
above, the authors considered the {\it unfolded eigenvalues} \beq
\Upsilon_{k,V}^0:=\ve V\ve(N(\rrl_{k,V})-N(\rrl^0)) \quad (1 \mly
k \mly \ve V\ve),
\eeq i.e., the eigenvalues under nonlinear renormalization. They
proved the following limit theorems for eigenvalue statistics:

1) Poisson limit theorem for the point process based on the
unfolded eigenvalues (1.19), where the limiting Poisson process
coincides with that for (1.18);

2) Poisson limit theorem for the point process based on both the
unfolded eigenvalues (1.19) and the normalized localization centers
of the corresponding eigenfunctions;

3) Limit theorems for the empirical distribution function of the
normalized spacings of eigenvalues close to $\rrl^0$;

4) Limit theorems for the normalized distance between localization
centers of the corresponding eigenfunctions; \\ 
and other important limit theorems for various
statistics related to the eigenvalues and eigenfunctions in the
intervals $I_{{\rm pp}}$. Moreover,
Germinet and Klopp (2013) considered more general random
Hamiltonians $\cH$ with convolution-type long-range kinetic
operators instead of the Laplacian: For such Hamiltonians, the
Poissonian asymptotic results were shown to hold also for the
eigenvalues at the spectral edges (in particular, for the $K$th
largest eigenvalues with fixed $K\dly 1$). At the spectral edges
of the Schr\" odinger operators, this result was shown to hold for
the one-dimensional case $\nu=1$; cf.~Section 6.2 below.

The proof of the limit theorems in (Germinet and Klopp 2013; 2014)
relies heavily on the techniques of the  general theory of
Anderson localization, including applications of the
fractional-moment criteria in the large-volume setting, as well as
various versions of the Wegner and Minami estimates.

In the proof of the above results, the following crucial
observation is related to localization properties of the
corresponding eigenfunctions of the operator $\cH_V$ (and is also
quite close to the context of our paper treating the top
eigenvalues and eigenfunctions of $\cH_V$): Namely, the
eigenvalues $\rrl\iz{k,V}$ near $\rrl^0$ of the operator $\cH_V$
can be approximated, with a \lk good\rk\ error, by independent local
eigenvalues of the operators restricted to much smaller disjoint
cubes $\wdt V(z)\subset V$ centered at $z$. This approximation is
feasible due to the facts that the localization centers of the
corresponding eigenfunctions of $\cH_V$ are located far away from
each other (so the probability of having at least two centers in cubes $\wdt
V(z)$ is asymptotically negligible), and the eigenvalues of $\cH_V$ \lk live\rk on potential values in small spatial neighborhoods of the
corresponding localization centers because of the exponential
decay of eigenfunctions. One needs to consider the small
spectral interval $I_V\subset I_{{\rm pp}}$, centered at $\rrl^0$,
such that the number of eigenvalues $\rrl\iz{k,V}$ in $I_V$
unboundedly increases as $V\uparrow \ttZ^\nu$. The latter is
fulfilled if, say, $\ve I_V\ve \asymp \ve V\ve^{-c}$ for some
$0<c<1$. The study of the structure of eigenvalues $\rrl\iz{k,V}$
inside $I_V$ is based on the Wegner and Minami estimates. In
particular, using the eigenvalue approximation described above,
one has with high probability that the number of eigenvalues in
$I_V$ is roughly approximated by $N(I_V)\ve V\ve$ (large deviation
principle); here $N(I)$ is the probability measure associated with the
integrated density of states. Also, note that the above
continuity conditions on potential distribution imply another
important bound: \beq N(I_V)\dly \const \ve I_V\ve^{1+\vartheta}
\eeq for some $\vartheta\dly 0$ and for all $V$ large enough. These bounds
are crucial in estimating the probabilities of the occurrence of a
single or several eigenvalues in small spectral intervals for
operators over cubes $V$ and $\wdt V(z)\subset V$ introduced above.

In (Germinet and Klopp 2013), the results on local fluctuations
for eigenvalues $\rrl\iz{k,V}$ are extended up to spectral edges.
Here the proofs rely on the improved versions of Wegner and Minami
estimates which ensure the more explicit control of eigenvalues
$\rrl \iz{k,V}$ in small spectral intervals $I_V$, since the
amount $N(I_V)$ is now allowed to be exponentially small in $\ve
I_V\ve^{-1}$ instead of (1.20). Thus, this case includes
situations, where the measure $N(\cdot)$ is extremely small, for
example, at the edges of spectral bands  where the phenomenon of
\lk Lifshits tails\rk occurs.

We finally notice that the Poisson limit theorems for the lower
eigenvalues (1.18) or (1.19) (in the spectral regions $I_{{\rm
pp}}$ of Anderson localization) agree with the results of the
present paper and our earlier papers on the extreme value theory for
the largest eigenvalues $\rrl\iz{K,V}$ of $\cH\iz{V}$, with $K
\dly 1$ fixed. However, in limit theorems for the largest
eigenvalues, the choice of normalizing constants depends
strongly on the regularity and tail decay conditions of the
potential distribution, in contrast to limit theorems for
eigenvalues near $ \rrl^0 \in I_{{\rm pp}}$, where the normalizing
constants are simply expressed in terms of the (integrated)
density of states. In particular, the spacings of the top
eigenvalues have an asymptotic order, which is much larger than
$\ve V\ve^{-1}=$ the order of the gaps between successive
eigenvalues close to $\rrl^0 \in I_{{\rm pp}}$. It is worth mentioning that the proof of Poisson limit theorems for eigenvalues
in $I_{{\rm pp}}$ relies heavily on the methods and ideas of the
general theory of Anderson localization; and neither the extreme value
theory nor links to the asymptotic geometric properties of random
potential are explored.

 We finally mention the following open
problems regarding the {\it Anderson transition phenomenon} for
the infinite-volume Hamiltonians $\cH=\rrk\Delta +\xi(\cdotp )$:
The first conjecture is that for $\nu \dly 3$, the spectral bands
outside some neighborhood of the spectral edges consist of purely
absolutely continuous spectrum and the corresponding
eigenfunctions are delocalized. The second conjecture is that the
eigenvalues of $\cH\iz{V}$ in the spectral intervals of
delocalization obey non-Poissonian asymptotic behavior as
$V\uparrow \ttZ^\nu$, rather the limit theorems like in the theory
of Wigner random matrices with light-tailed entries; cf.~Section 1.4 below. The first
problem is partially solved for the very special models on the
Bethe lattice as well as for the Schr\" odinger operators with
sparse potential (e.g., Kirsch 2008; Molchanov and Vainberg 1998,
2000).

\medskip

\subsection{\normalsize{\textsf {Relations to random matrices } }}

{\bf (I)} {\it Wigner random matrices}. Another important model of
disordered quantum systems (in particular, heavy nuclei atoms) is
presented by real symmetric random matrices
$$
H_N=\big(h\iz{i,j}\big)\iz{1\mly i,j\mly N}
$$
with i.i.d.~centered entries $h\iz{i,j}$ ($i\mly j$) and $N
\rightarrow \infty$; i.e., large {\it Wigner matrices} (Mehta
2004; Anderson et al.~2010). The extreme eigenvalues (i.e., high
energy levels) $\rrl\iz{K,N}$ and the corresponding
$l^2$-normalized eigenvectors of $H_N$ are here interpreted as the
basic states of quantum systems.

Recently, there has been much progress toward the extreme value theory
for the eigenvalues $\rrl\iz{K,N}$ of Wigner matrices $H_N$ as
$N\to \infty$ and $K\dly 1$ fixed. It has been turned out that
there are two different regimes of asymptotic behavior of the
largest eigenvalues, depending on the tail decay rate of the
entries in absolute value:

(${\rm I}_1$) For polynomially decaying distributions \beq
\bP\big(\big |h\iz{i,j}\big|>t\big)= t^{-\beta}(1+\mro(1))\quad
\hbox {as}\ \ t\to \infty
\eeq with $\beta <4$ (very heavy tails), Auffinger et al.~(2009)
proved that with high probability, the $K$th largest eigenvalue
$\rrl\iz{K,N}$ of Wigner matrices $H_N$ is approximately equal to
the $K$th largest value among $\vert h\iz{i,j}\vert$ ($1\mly i\mly
j\mly N$), $K \dly 1$ fixed. This in turn implies Poisson limit
theorems for the normalized eigenvalues $\rrl\iz{K, N}A_N$, where
the normalizing constants $A_N>0$ are chosen the same as in the
corresponding limit theorems for extremes of $\vert
h\iz{i,j}\vert$ ($1\mly i\mly j\mly N$). Recall that
distributional tails (1.21) are in the domain of attraction of the
max-stable Fr\'{e}chet law $G_\beta$, therefore, the eigenvalue
$\rrl\iz{K, N}$ is of the order $A_N^{-1}= N^{2/\beta}(\const+\mro(1))$;
cf.~Example 6.11 below. Moreover, with high probability, the $K$th
eigenvector is asymptotically concentrated on two coordinates,
i.e., it behaves like a superposition of two delta functions in
limit as $N \to \infty$. See also (Soshnikov 2004) for the case
$\beta <2$. To prove these assertions, one first observes that,
under condition (1.21) with $\beta <4$, the larger entries (in
absolute value) are extremely sparse and strongly pronounced in
comparison to other entries in $H_N$. Thus, the standard
perturbation theory for symmetric matrices is applied to conclude
that the top eigenvalues and eigenvectors of the former matrix
$H_N$ are approximated by the corresponding eigenvalues and
eigenvectors of a very sparse (symmetric $N\times N$) matrix whose
entries are the larger values of $ h\iz{i,j}$'s in absolute value.
This fact in turn enables us to apply the extreme value theory for the random
variables $\ve h_{i,j}\ve$ and, as a consequence, to prove Poisson
limit theorems for the eigenvalues of Wigner matrices $H_N$.

(${\rm I}_2$) Assume that the $\ve h_{i,j}\ve$'s have lighter
tails (including (1.21) with $\beta >4$ and Bernoulli entries), or
the $ h_{i,j}$'s have some finite statistical moments of higher
order and satisfy additional conditions on a distributional
symmetry. Then the largest eigenvalues $\rrl\iz{K,N}$ of Wigner
matrices $H_N$ are distinguished by non-Poissonian asymptotic
behavior, rather the Tracy-Widom limit law; see, e.g., (Soshnikov
1999; Lee and Yin 2014; Bourgade et al.~2014). In particular, the
normalized top eigenvalues $\rrl\iz{K,N}/\sqrt{N}$ tend almost
surely to the nonrandom constant $2 (\ttE h_{1,2}^2 )^{1/2}$,
i.e., the right endpoint of the support of their limiting spectral
distribution density; cf.~(Bai and Yin 1988). Moreover, the
corresponding $l^2$-normalized eigenvectors are completely
delocalized; i.e., with high probability their sup-norm does not
exceed $N^{-1/2}(\log N)^{\const}$. See, e.g., (Tao and Vu 2010; 
Erd\H{o}s et al.~2013a; Vu and Wang 2015; G\"{o}tze et al.~2015), where this delocalization property is
extended to all the eigenvectors of $H_N$ provided the tails of $\ve
h_{i,j}\ve$ are lighter than the exponential, or
the $ h_{i,j}$'s have a large enough number of moments. It is worth
mentioning that limiting distributions for eigenvalues or
eigenvectors (in particular, the Tracy-Widom limit law for the
largest eigenvalues) can be explicitly computed for Wigner
matrices with Gaussian entries; see, e.g., (Anderson et al.~2010).
Thus, the usual comparison methods (four moments theorem, Green
function comparison method, etc.) can be used to extend the asymptotic
results for the Gaussian case to general Wigner matrices; e.g.,
(Tao and Vu 2014).

The transition from Poisson limit theorems to Tracy-Widom
asymptotics for the top eigenvalues of Wigner random matrices was
discussed in detail by Biroli et al.~(2007). The value $\beta=4$
in (1.21) or, roughly speaking, the fourth statistical moment
indicates here the threshold separating these two different regimes
of asymptotic behavior. Note that the Wigner matrix model with
light-tailed entries reflects the global or mean-field interaction;
thus, the asymptotic geometric properties of entries do not play any
role in limit behavior of eigenvalues and eigenvectors. The latter
is in a sharp contrast to the one-dimensional Anderson Hamiltonian
in $V \subset \ttZ$, which is a random band (tridiagonal) matrix
reflecting local interaction: the diagonal elements are i.i.d.~random 
variables and the deterministic off-diagonal elements are
given by the Laplacian. In view of discussions on the Anderson
model (Section 1.2 above), it turns out that the diagonal operator
necessitates concentration properties of eigenvectors; meanwhile,
the Laplacian forces these properties to be less expressed (thus,
the geometric features of the model play here a crucial role).

{\bf (II)} {\it Random band matrices}. Recently, there has been a
considerable attention drawn to symmetric {\it random band matrices}
$H_N^{(W)}$ of size $N \to \infty$, where the matrix entries $
h\iz{i,j}$ vanish if $|i-j|$ exceeds $W$, and other entries
(above the diagonal) are i.i.d.~centered random variables; here
$0\mly W \mly N$ is a band width. Random band matrices are natural
interpolations between Anderson Hamiltonians ($\nu =1$) and Wigner
matrices. For band models, the asymptotic behavior of the top
eigenvalues and eigenvectors depends strongly on the growth rate
of the band width $W=W_N \to \infty$ as well. Benaych-Georges and
P\'ech\'e (2014) considered the random band matrices $H_N^{(W)}$,
whose entries have polynomially decaying distributions (1.21) with
arbitrary $\beta
>0$ and band width $W=N^\mu$ with $0<\mu \mly1$. They established
that the band model exhibits a phase transition depending on $\mu$
and $\beta$, with $\beta=2(1+\mu^{-1})$ as the threshold
separating two different regimes of asymptotic behavior of the
largest eigenvalues $\rrl\iz{K,N}$ and the corresponding
eigenvectors ($K \dly 1$ fixed):

(${\rm II}_1$) Assume (1.21) and $W=N^\mu$ such that
$\beta<2(1+\mu^{-1})$, i.e., either the distributional tails are
sufficiently heavy or the band of matrix is sufficiently narrow.
This case includes heavy-tailed Wigner matrices, i.e., $\mu=1$ and
$\beta <4$, considered in (${\rm I}_1$) above. For $0<\mu<1$ and
$\beta<2(1+\mu^{-1})$, the asymptotic results are similar to that
in (${\rm I}_1$). I.e., with probability $1+\mro(1)$, the $K$th
largest eigenvalue of $H_N^{(W)}$ is approximately equal to the
$K$th largest value of the sample $\vert h\iz{i,j}\vert$ ($1\mly
i\mly N$, $0\mly j-i \mly W$). Therefore, Poisson limit theorems
for the normalized eigenvalues $\rrl\iz{K, N}A_N$ hold true, where
$A_N= N^{-(1+\mu)/\beta}(\const^\prime +\mro(1))$ (cf.~Example 6.11
below), and the $K$th eigenvector is asymptotically localized on
two coordinates. The proof of these assertions is again heavily
based on techniques of the extreme value theory, in particular,
describing asymptotic geometric properties of entries of band
matrices. The latter is combined with the perturbation theory for
matrices to derive simple asymptotic formulas for the largest
eigenvalues and eigenfunctions of $H_N^{(W)}$.

(${\rm II}_2$) Assume the $h\iz{i,j}$'s are symmetrically
distributed with tails (1.21) and let $W=N^\mu$ with
$\beta>2(1+\mu^{-1})$, i.e., either the distributional tails are
sufficiently light or the band of matrix is sufficiently wide. In
this case, the band models $H_N^{(W)}$ posses the mean-field
features, like in the light-tailed Wigner models considered in
(${\rm I}_2$) above. Thus, each eigenvector associated with the
upper spectral edge is asymptotically delocalized in the sense
that its $l^2$ mass is more or less uniformly spread over $N$
coordinates. Moreover, the extreme eigenvalues do not longer obey
Poissonian asymptotic behavior; in particular, they tend to a
nonrandom positive constant when divided by $N^{\mu/2}$.

Earlier, Sodin (2010) considered the band matrices $H_N^{(W)}$
whose entries are symmetrically distributed with tails lighter
than the Gaussian tails, including Bernoulli entries. He studied
the transition from localization to delocalization for
eigenvectors at the upper spectral edge (as well as the
corresponding limit theorems for eigenvalues) with varying degrees
of strength and generality. In particular, the results of this
paper suggest that the eigenvectors associated with the largest
eigenvalues are delocalized, provided $W=N^\mu$ with $\mu
> \mu_{{\rm cr}} =5/6$. See also (Erd\H{o}s et al.~2013b) for
similar delocalization results for eigenvectors associated with
the inner part of the spectrum $\spect (H_N^{(W)})$.

In view of these observations on the localization properties at the
upper spectral edge, we distinguish two paradigmatic models in the
theory of random matrices. First, the general theory of Anderson
localization suggests that, with probability 1, the
one-dimensional Schr\" odinger operators $\cH=\rrk\Delta +\xi$
(i.e., the particular model of the band matrices with $W=1$) have
exponentially localized eigenvectors at the spectral edges,
provided the potential distribution is arbitrary satisfying very mild
continuity conditions (Carmona et al.~1987). On the other hand,
the large Wigner matrices $H_N$ (i.e., the band matrices with
$W=N$) have localized eigenvectors associated with the largest
eigenvalues if only entries are very heavy-tailed like in (1.21)
with $\beta < 4$; meanwhile, for the light tails, eigenvectors of
$H_N$ are typically delocalized. See also (Spencer 2011) for a
discussion on Anderson-type models ($M=\rO(1)$), band matrices ($M
=\mro(N)$) and Wigner matrices ($M=N$).

Finally, it is worth emphasizing that the above asymptotic results
for eigenvalues remain valid for the corresponding complex
Hermitian random matrices with i.i.d.~entries, instead of the real
symmetric matrices.

\medskip

\subsection{\normalsize{\textsf {The earlier literature on extremes of i.i.d.~random fields } }}

As already mentioned, most statements of the present paper on the $\xi_V$-extremes and
the corresponding RV classes were announced
in (Astrauskas 2007; 2008; 2012; 2013). We now provide a brief
overview of the earlier literature on the related asymptotic
results for extreme order statistics of i.i.d.~random sequences
and fields.

High-level exceedances consisting of single rare $\xi_V$-peaks
were studied in (Astrauskas 2001). Related asymptotic results (in
particular, the so-called longest head runs in coin tossing) for
Bernoulli distributed i.i.d.~random variables $\xi (x)$, $x \in
\ttZ$, were discussed, e.g., in (Binswanger and Embrechts 1994).

In the case of exponentially distributed $\eta(0)$, strong limit
theorems for the spacings $\eta\iz{K,V}-\eta\iz{K+1,V}$ ($K$
fixed) were proved by Astrauskas~(2006). Devroye (1982) derived
strong and weak limit theorems for $\min_{1\mly k \mly \ve
V\ve}(\zeta\iz{k,V}-\zeta\iz{k+1,V})$ where $\zeta(0)$ is uniformly
distributed.

In the case of $\xi(0)$ with arbitrary distribution, 
strong asymptotic bounds for $\xi\iz{K,V}$ are given in (Shorack
and Wellner 1986) where $K$ is fixed, and in (Deheuvels~1986) where
$K=K_V \to \infty$. Wellner (1978) derived strong asymptotic
bounds for the uniform $k$th order statistics $\zeta\iz{k,V}$
(thus, for $\eta_{k,V}$) uniformly in $k \dly 1$.

For the Gaussian random fields $\{\xi(x) \colon x \in
\ttZ^\nu\}$ with correlated values, Astrauskas (2003) studied some
asymptotic geometric properties of $\xi_V$-extremes almost surely,
in particular, high level exceedances and minimum of spacings. The
geometry of high level excursion sets of smooth Gaussian random
fields in $\ttR^\nu$ was investigated in the monograph by Adler
and Taylor (2007). See also (G\" artner et al.~2000) for some
geometric aspects of high peaks of smooth Gaussian random
potentials related to the long-time asymptotics for the spatially
continuous parabolic Anderson models.

For the extreme value theory for random variables, in particular,
characterization of the domains of attraction of max-stable
distributions, we refer to the monographs by Resnick (1987), de
Haan and Ferreira (2006), Leadbetter et al.~(1983), Embrechts et
al. (1997). See also the monograph by Shorack and Wellner (1986)
for a detailed account of strong and weak limit theorems for order
statistics and their functions related to mathematical statistics.
Finally, the monograph by Bingham et al.~(1987) provides a
detailed account of the theory of regularly varying functions.

In the proof of a number of our statements on $\xi_V$-extremes, we
explore the representation $\xi\iz{k,V}=f(\eta\iz{k,V})$, where
$f:=Q^{\leftarrow }$ is the generalized inverse function of $Q$
and, as above, $\eta\iz{k,V}$ stands for the $k$th extreme value
among independent exponentially distributed random variables
$\eta(x)$ ($x\in V$) with mean 1. Due to the nice properties of
$\eta\iz{k,V}$ (for instance, $\eta\iz{k,V}$ is a sum of
independent exponentially distributed random variables), we first
obtain the asymptotic results for $\eta\iz{k,V}$, which are then
transferred to $\xi\iz{k,V}$ under appropriate conditions on $f$.
These conditions are formulated in terms of regular variation (RV)
of $f(s)$ as $s\to \infty$. We further give a characterization of
RV classes, in particular, their links to continuity and tail
decay of the distribution $1-\e^{-Q}$ at the right endpoint; see
Appendix A. They are also compared with the well-known RV classes
including the domains of attraction of max-stable laws,
$\rO$-regular variation, asymptotically balanced, etc; see
Appendix A.

An interesting further problem is an extension of the present
asymptotic results for functionals (EX), (SP), (MIN), (NEI) to other
classes of random fields $\xi(\cdotp)$ including: 1) independent
non-identically distributed random variables; 2) random fields
with correlated values, in particular, Gaussian fields (Astrauskas
2003) and moving average fields defined as a linear combination of
i.i.d.~random variables with nonrandom real coefficients. See,
e.g., the review papers by Elgart et al.~(2012), Tautenhahn and
Veseli\'c (2015) for a detailed background of the random alloy
type models $\rrk \Delta +\xi(\cdot)$ with the moving average
potential $\xi(\cdot)$.

\medskip

\subsection{\normalsize{\textsf {Notation. Representation of i.i.d.~random fields } }}

Let us introduce the further notation and remarks we use
throughout the paper. We denote by $\ttR_+$ the positive half-axis
and by $\ttN$ positive integers. Let $\log _j$ stand for the
$j$times iterated natural logarithm. For real $a$, $b$, we write
$a\vee b:=\max(a, b)$ and $a\land b:=\min(a, b)$, and $[a]$ for
the integer part of $a$. Given a subset $U\subset \ttZ^\nu $, we
write $\vert U\vert $ for the number of its elements. Let $\dist
(U, U\pr )$ stand for the lattice $l^1$-distance between subsets
$U, U\pr \subset \ttZ^\nu$. The summation over $x\in V\colon\
a\mly \vert x\vert \mly b$ is abbreviated to $\sum _{a\mly \vert x
\vert \mly b}$. By $t_0$, $\vert V_0\vert $, etc. we denote
various large numbers, values of which may change from one
appearance to the next. Similarly, $\const $, $\const\pr $ etc.
stand for various positive constants. We write $1/0+=\infty$,
$\log(0+)=-\infty$ and $1/ \infty =0$. Let $g\circ h=g(h(\cdot))$
stand for a composition of real functions $g$ and $h$. Also, for
$g>0$ and $h>0$, we write $g(t)\asymp h(t)$ as $t \to \infty $, if
the ratio $g(t)/h(t)$ is bounded away from zero and from above for
all large $t$.

By $\cG\iz{V}(\rrl;\zeta\iz{V}; x, y)$ ($x\in V$, $y\in V$) we
denote the Green function of the Hamiltonian $\rrk\Delta\iz{V} +
\zeta\iz{V}$ in $l^2(V)$, viz.
$$
\cG\iz{V}(\rrl;\zeta\iz{V}; x,
y):=\cG\iz{V}(\rrl;\zeta\iz{V})\delta
_y(x):=(\rrl-\rrk\Delta\iz{V} - \zeta\iz{V})^{-1}\delta _y(x).
$$
Here $\delta _y(\cdotp )$ is the Kronecker delta function, i.e.,
$\delta _y(x):=1$ if $x=y$, and $\delta _y(x):=0$ if $x\neq y$.

Throughout the paper we suppose that all random variables are
defined on a common probability space $(\rrO ,\cF ,\bP )$. Let
$\bE $ stand for the expectation with respect to $\bP $. Recall
$Q(t):=-\log \bP(\xi(0)>t)$ is the cumulative hazard function of
an i.i.d.~random field $\xi(x)=\xi^{(\omega)}(x)$ ($\omega \in
\rrO$; $x\in \ttZ^\nu$), and let $t_Q$ denote its right endpoint
$t_Q:=\sup \lf \{t\colon\ Q(t)<\infty \rg \}$. Without loss of
generality, we shall assume throughout that $0<t_Q\mly \infty$.
Clearly $Q:(-\infty;t_Q)\rightarrow \ttR_+\cup\{0\}$ is a right-continuous
nondecreasing function such that $Q(-\infty)=0$ and
$Q(t_Q)=\infty$. Most of the conditions of our results are
formulated in terms of the inverse of the cumulative hazard
function defined by  \beq f(s):=Q^{\leftarrow }(s):=\inf\lf \{
t\colon\ Q(t)\dly s\rg \} \quad (s \in \ttR_+ )
\eeq (thus $f:\ttR_+\rightarrow (-\infty;t_Q)$ is a
left-continuous nondecreasing function such that $f(s)$ tends to
$t_Q$ as $s\to \infty$). The reason for this is the following
useful {\it representation of order statistics} $\xi\iz{k, V}$:
\beq \xi\iz{1, V}:=f(\eta\iz{1, V})\dly \xi\iz{2, V}:=f(\eta\iz{2,
V})\dly \ldots \dly \xi\iz{\vert V\vert , V}:=f(\eta\iz{\vert
V\vert , V}),
\eeq where \beq \eta\iz{1, V}:=\eta (z\iz{1, V})>\eta \iz{2,
V}:=\eta (z\iz{2, V})>\ldots >\eta \iz{\vert V\vert ,V}:=\eta
(z\iz{\vert V\vert ,V})
\eeq is the variational series based on the sample $\eta_V :=\{
\eta(x)\colon\ x\in V\}$ of exponential i.i.d.~random variables
with mean 1.

\medskip

\subsection{\normalsize{\textsf {Outline } }}

In Section 2, we collect conditions on deterministic functions
$\xi_V$ in terms of functionals (EX), (SP), (MIN), and (NEI), which
yield expansion formulas for the largest eigenvalues
$\rrl\iz{K,V}$ of the discrete Schr\" odinger operator
$\cH_V=\rrk\Delta_V +\xi_V$ on $l^2(V)$ as $V\uparrow \ttZ^\nu$.
Section 2.1 provides rough bounds for $\rrl\iz{K,V}$. We then
study $\rrl\iz{K,V}$ in the cases of $\xi_V$ with extremely sharp
peaks (Section 2.2), dominating single peaks (Section~2.3),
dominating large islands of high $\xi_V$-values the diameter of
which unboundedly increases (Section~2.4) and, finally, dominating
islands of high $\xi_V$-values the diameter of which is bounded
(Section~2.5). The results of Sections 2.2--2.5 follow simply from
the more general statements of (Astrauskas 2008; 2012) and Section
2.4 in (G\"artner and Molchanov 1998).

Sections 3--5 contain the main results of the paper dealing with
asymptotic behavior of extremes of the i.i.d.~random field
$\xi(\cdotp)$ with the distribution function satisfying certain RV and
continuity conditions at the right endpoint. Functionals (EX),
(SP), (MIN) and (NEI) are studied in Sections 3, 4.1-4.2, 4.3 and 5,
respectively.

Section~6 provides an overview of current results on extreme value
theory for the spectrum of the Anderson Hamiltonian $\cH_V=\rrk
\Delta_V +\xi_V$, $V\uparrow \ttZ^\nu$, with an i.i.d.~potential
$\xi(\cdot)$. The issues under discussion include the asymptotic
expansion formulas and Poisson limit theorems for the largest
eigenvalues and their localization centers. We consider separately
three cases of the distribution tails $\e^{-Q}$ of $\xi(0)$: the
tails are heavier than the double exponential function
(Section~6.1); the tails are lighter than the double exponential
function (Section~6.2), and the double exponential tails
(Section~6.3). As already mentioned, we give proof sketches of most theorems of this
section demonstrating their connections to the results of Sections
3--5 on $\xi_V$-extremes. In Section 6.4, we comment and compare
the proofs of Poisson limit theorems stated in Sections~6.1
and~6.3 and proved in the earlier papers by Astrauskas and
Molchanov (1992), Astrauskas~(2007; 2008; 2012; 2013) and Biskup
and K\"onig (2016).

In Section 7, we discuss the long-time intermittent behavior of
the solutions to the parabolic problems associated with the
Anderson Hamiltonian $\cH=\rrk\Delta +\xi(\cdot)$. We focus on the
representation of the solutions  in the spectral terms of the
operators $\cH_V=\rrk\Delta_V +\xi_V$. In view of this
representation, we discuss some techniques of the extreme value theory
for eigenvalues of $\cH_V$, that can be applied to study the
intermittency properties of time-dependent Anderson models.

Finally, in Appendix A, we characterize and compare the RV classes
of distributions introduced in Sections 3--6.

\medskip

\setcounter{equation}{0} \addtocounter{sct}{1}


\section{\normalsize{\textsf {ASYMPTOTIC EXPANSION FORMULAS FOR THE LARGEST EI\-GEN\-VA\-LUES OF DETERMINISTIC HAMILTONIANS} }}

\setcounter{prop}{1}

Let $V=[-n;n]^\nu \cap \ttZ^\nu$ ($n \in \ttN$) be a sequence of
cubes. By introducing the periodic norm $\ve x\ve:=\ve x\ve_n:=
\min_{y\in (2n+1)\ttZ^\nu }\ve x-y\ve$, $V$ may be considered as a
sequence of tori tending to $\ttZ^\nu$. We are interested in the
finite-volume Schr\" odinger operators $\cH_V=\rrk\Delta_V +\xi_V$
on $l^2(V)$ with periodic boundary conditions. Recall that  $\rrk
>0$ is a diffusion constant, $\Delta_V $ denotes the lattice
Laplacian on $l^2(V)$ (i.e., a restriction of the operator $\Delta
\psi(x):=\sum _{\vert y-x\vert =1}\psi(y)$ to torus $V$) and
$\xi_V :=\{ \xi_V(x)\colon\ x\in V\}\in [-\infty;\infty)^{\ve
V\ve}$ are deterministic functions, i.e., potential. The values
$-\infty$ of $\xi_V$ (i.e., \lk hard obstacles\rk\ ) are allowed
to include the cases which are interesting from a physical point
of view; see, e.g., (Biskup and K\"onig~2001; K\"onig~2016). Write
$V_b:= \{x \in V \colon\ \xi_V(x) > -\infty \}$. Then $\cH_V$ is
interpreted as an operator on $l^2(V)$ with zero boundary
conditions outside $V_b$. The spectral problem \beq
\cH_V\psi=\rrl\psi \quad \quad (\rrl \in \ttR; \ \psi \in l^2(V))
\eeq has $\ve V_b\ve$ solutions $\rrl\iz{1,V}\dly \rrl\iz{2,
V}\dly \ldots \dly \rrl\iz{\vert V_b\vert ,V}$, i.e., the ordered
eigenvalues of the operator $\cH_V$.

In this section, we provide asymptotic expansion formulas for the
first $K$ largest eigenvalues $\rrl\iz{k,V}$ under conditions on
the first terms of the variational series $\xi\iz{1,V}\dly
\xi\iz{2, V}\dly \ldots \dly \xi\iz{\vert V\vert ,V}$ of the
sample  $\xi_V$ and their coordinates $z\iz{k,V} \in V$ defined by
$\xi\iz{k,V}=\xi\iz{V}(z\iz{k,V})$ ($1\mly k \mly \ve V\ve$); here
$V= \{z_{k, V} \colon\ 1\mly k \mly \ve V\ve \}$. The results of
this section follow simply from more general results of
(Astrauskas 2008; 2012) and Section 2.4 of (G\" artner and
Molchanov 1998), where one finds more discussions on the
relationship between $\xi_V$-extremes and the top eigenvalues of
$\cH_V$.

In this section, the proof of the statements relies on deterministic 
spectral arguments. It is worth mentioning that the
(probabilistic) Feynman-Kac representations of the Green function
and the principal eigenfunction of Schr\" odinger operators $\cH_V$ as
well as the related path decomposition techniques present powerful
probabilistic tools for deriving the explicit upper bounds for the
principal eigenvalue and eigenfunction of $\cH_V$ (G\"artner and
Molchanov 1998; G\" artner et al.~2007). However, this method is
not explained in the present section; see Section 7 below for some aspects 
of these techniques related to the parabolic Anderson models.

\subsection{\normalsize{\textsf {Preliminaries: rough bounds } }}

We start with the following simple bounds for eigenvalues $\rrl
\iz{k,V}$, provided $\vert V_b\vert \dly 2$.

\begin{thm}
{\rm (i)} For any $V$ and any $\xi_V$, \beq \xi\iz{1, V}\mly
\rrl\iz{1, V}\mly \xi\iz{1, V}+2\nu \rrk \quad \hbox{and}\quad \ve
\rrl\iz{l, V}-\xi\iz{l, V}\ve \mly 2\nu \rrk \quad (2\mly l \mly
\ve V_b \ve).
\eeq

{\rm (ii)} For any $V$, any $K \mly \ve V_b\ve$ and $\xi_V$ such
that $\min_{1\mly k<l \mly K} \ve z\iz{k,V}-z\iz{l,V}\ve \dly 2$,
we have that \beq \xi\iz{l, V}\mly \rrl\iz{l, V}\mly \xi\iz{l,
V}+2\nu \rrk \quad \hbox{for all} \quad 1 \mly l \mly K.
\eeq
\end{thm}

\medskip

\begin{proof}
We repeatedly use the fact that the $K$th  eigenvalue
$\rrl_{K,V}=\rrl_{K,V}(\xi_V)$ of the operator $\rrk
\Delta_V+\xi_V$ is a nondecreasing function in each variable
$\xi_V(x)$ tending to infinity ($K \dly 1$, $x \in V$); i.e., the
monotonicity property of eigenvalues (Lankaster~1969, Theorem
3.6.3).

(i) To estimate $\rrl\iz{1, V}$, we abbreviate $\xi'(x):=\xi\iz{1,
V}$ if $x=z\iz{1, V}$, and $\xi'(x):=-\infty $, otherwise. Note
that $\xi'(\cdotp )\mly \xi\iz{V}(\cdotp) \mly \xi\iz{1, V}$ in
$V$. Therefore, $\rrl\iz{1, V}$ is bounded from below by
$\xi\iz{1, V}$, i.e., the principal eigenvalue of the operator
$\rrk\Delta_V +\xi'_V$. Moreover, $\rrl\iz{1, V}$ is bounded from
above by the principal eigenvalue of the operator $\rrk\Delta_V
+\xi\iz{1, V}$ on $l^2(V)$, which in turn does not exceed $2\nu
\rrk+\xi\iz{1, V}$, since the norm of the Laplacian $\rrk\Delta_V
$ is less than $2\nu \rrk $. Similarly, since each eigenvalue
$\rrl\iz{l, V}$ is bounded from above (resp., from below) by the
$l$th eigenvalue of the diagonal operator $\xi_V +2\nu \rrk $ on
$l^2(V)$ (resp., $\xi_V-2\nu \rrk $), we obtain (2.2) for $l\dly
2$.

(ii) We need to show the lower bound in (2.3). Without loss of
generality, we assume that $\xi \iz{K,V}>0$ (this may be achieved
by shift transform of $\xi_V$ and $\rrl$ in the spectral problem
(2.1)). Write $\cE^K_V:=\{z\iz{1, V},\ldots, z\iz{K, V}\}$. We
introduce the following functions: $\zeta(x):=\xi_V(x)$ if $x\in
\cE^K_V$, and is zero, otherwise; and further on, $\wdt
\zeta(x):=0$ if $x\in \cE^K_V$, and $\wdt \zeta(x):=-\infty$,
otherwise. Then $\xi_V(\cdotp) \dly \zeta(\cdotp)+\wdt
\zeta(\cdotp)$ in $V$, therefore, each eigenvalue $\rrl \iz{l, V}$
is bounded from below by the corresponding eigenvalue
$\underline{\rrl}_{l,V}$ of the operator $\rrk\Delta_V +\zeta_V
+\wdt \zeta_V$; here $1 \mly l\mly K$. To estimate
$\underline{\rrl}_{l,V}$, we rewrite the corresponding spectral
problem in the form: \beq (\rrl -\rrk\Delta_V -\wdt \zeta_V)\psi
=\zeta_V \psi \ \ \quad (\rrl>0, \psi \in l^2(V))
\eeq and apply the resolvent operator $\cG_V(\rrl; \wdt
\zeta_V):=(\rrl-\rrk\Delta\iz{V}-\wdt \zeta_V)^{-1}$ to both sides
of~(2.4). Since $\cG_V(\rrl; \wdt
\zeta_V)\delta_z=\rrl^{-1}\delta_z$ for $z \in \cE^K_V$, equation
(2.4) is transferred to
$$
\psi =\sum_{z\in \cE^K_V}\xi_V(z)\psi(z)\rrl^{-1}\delta_z \ \
\quad (\rrl >0);
$$
here $\delta _y(\cdotp )$ is the Kronecker delta function.
Clearly, for each $1 \mly l \mly K$, the pair
$\underline{\rrl}_{l,V}=\xi _{l, V}$ and $\psi(\cdotp;
\underline{\rrl}_{l,V})=\delta_{z \iz{l,V}}(\cdotp)$ solves this
equation. Summarizing, we have that $\rrl _{l,V} \dly
\underline{\rrl}_{l,V}=\xi _{l, V}$ ($1 \mly l \mly K$), as
claimed. Theorem 2.1 is proved.\hfill$\square$
\end{proof}

In Sections 2.2--2.5 below, we consider three classes of functions
$\xi_V$:

\textbf{(J)} Sparse distinct $\xi_V$-peaks dominate in the
landscape of $\xi_V$ as $V\uparrow \ttZ^\nu$, i.e., $\xi_V$
possess properties like (1.9)--(1.11). Then the $K$th largest
eigenvalue $\rrl\iz{K,V}$ is associated with an isolated peak
$\xi_{\tau(K),V}$, so that $\rrl_{K,V} \leftrightarrow
z_{\tau(K),V}$ for some $\tau(K)=\tau_ V (K) \dly 1$ (Section
2.3). In particular, if the functions $\xi_V$ possess extremely
sharp peaks like (1.2)--(1.4), then the eigenvalue $\rrl\iz{K,V}$
is associated with the $K$th largest value of $\xi_V$, viz.,
$\rrl_{K,V} \leftrightarrow z_{K,V}$ (Section 2.2). In both cases,
the lower bounds in (2.3) are achieved as $V\uparrow \ttZ^\nu$.

\textbf{(JJ)} The landscape of $\xi_V$ is dominated by flat
islands of large values with an unboundedly increasing diameter.
Then the largest eigenvalues are associated with such relevant
islands. In this case, the upper bounds in (2.2) are achieved as
$V\uparrow \ttZ^\nu$ (Section 2.4).

\textbf{(JJJ)} Similarly as in (JJ), bounded islands of large
values prevail in the landscape of $\xi_V$. Then the asymptotic
expansion terms of the principal eigenvalue $\rrl\iz{1,V}$ fill
the gap between its lower and upper bounds in (2.2) (Section 2.5).

In the case of (J) we obtain the explicit expansion formulas for
eigenvalues in terms of $\xi_V$-values. Meanwhile, for (JJ) and
(JJJ) we restrict ourselves to a derivation  of the second order
expansion formulas for eigenvalues.

\medskip

\subsection{\normalsize{\textsf {Potentials with extremely sharp single peaks} }}

For $N \dly 2$, let us write \beq \cE_V^N :=\{z \iz{1,V}, z
\iz{2,V},\ldots,z \iz{N,V}\} \subset V
\eeq for the subset of coordinates of the first $N$ largest values
of $\xi_V$, and \beq r \iz{N,V} =\min_{1\mly l <k \mly N}\ve
z\iz{l,V}-z\iz{k,V} \ve=\min_{x,y \in \cE_V^N \atop{x\neq y}}\ve
x- y\ve
\eeq for the minimum distance between sites in $\cE_V^N$. For
natural $1\mly K=K_V <N=N_V <\ve V\ve$ and $p\dly 0$, we introduce
the following conditions on functions $\xi_V$: \beq &&\lim _V \min
_{1\mly l\mly K}\xi_{(l+1)\wedge K,V}^p(\xi \iz{l,V}-
\xi\iz{l+1,V})=\infty\ \ \hbox{where} \ \ \lim_V \xi\iz{K,V}=\infty, \\[8pt]
&&C:=\lsup _V \frac {\xi\iz{N,V}} {\xi\iz{K,V}} <1, \\[8pt]
&&\lim _V \frac {r\iz{N,V}} {\log N}=\infty , \\[8pt]
&&{\ds M :=\lsup _V \max_{1\mly l\mly K}\max
_{\br{x-z_{l,V}}=1}\ve \xi_V(x)\ve < \infty}
\eeq and, finally, \beq
 \ds \lim _V \min_{K\!+\!1 \mly l\mly N}\xi_{K,V}^2 \bigg
(\xi\iz{K,V}\!+\!\frac
{2\nu\rrk^2}{\xi\iz{K,V}}\!-\!\xi\iz{l,V}\!-\!\rrk^2\sum
_{\br{x\!-\!z_{l,V}}=1}\frac{1}{\xi\iz{K,V}\!-\!\xi_V(x)}\bigg
)\!=\!\infty.
\eeq We write $\xi\iz{0,V}:=\infty$, and
$s_V(l):=(\xi\iz{l-1,V}-\xi\iz{l,V})\wedge
(\xi\iz{l,V}-\xi\iz{l+1,V})$ for $1 \mly l \mly K$.

\begin{thm}
{\rm (i)} Under {\rm (2.7)} with $p=0$, we have that
$$
\lsup_V \max_{1 \mly l\mly K}\big|\rrl\iz{l,V}-\xi\iz{l,V}\big|
s\iz{V}(l) \mly \const_1(\rrk,\nu).
$$
and
$$
\lsup_V \max_{1\mly l \mly K}\max_{x \neq z_{l,V}}\frac
{\log\big|\psi(x;\rrl\iz{l,V})\big|}
{|x-z\iz{l,V}| \log s\iz{V}(l)} \mly -1.
$$

{\rm (ii)} Under {\rm (2.7)--(2.9)} with $p=1$, we have that
$$
\lsup_V \max_{1 \mly l\mly K}\big|\rrl\iz{l,V}-\xi\iz{l,V}\big|
\xi\iz{l,V} \mly  \frac{\const_2(\rrk, \nu)}{1-C}.
$$
and \beq
\lsup_V \max_{1\mly l \mly K}\max_{x \neq z_{l,V}}\frac
{\log\big|\psi(x;\rrl\iz{l,V})\big|}
{|x-z\iz{l,V}| \log \xi\iz{l,V}} \mly -1.
\eeq

{\rm (iii)} If $\xi_V$ satisfies {\rm (2.7)--(2.11)} with $p=2$,
then
$$
\lsup_V \max_{1 \mly l\mly K} \bigg| \rrl\iz{l,V}-
\xi\iz{l,V}-\frac {2\nu\rrk^2}{\xi\iz{l,V}}\bigg|\xi_{l,V}^2 \mly
\frac{M \cdotp \const_3(\rrk, \nu)}{(1-C)^2}+ \const_4(\rrk, \nu).
$$
and {\rm (2.12)} holds true.
\end{thm}

\medskip

\begin{proof} We first note that condition (2.10) implies
$\xi \iz{N,V} >-\infty$ for any $V \supset V_0$. On the other
hand, if $\xi \iz{N,V}=-\infty$ and $r \iz{N,V} >1$, then
$\rrl\iz{l,V}=\xi\iz{l,V}$ for all $1\mly l\mly N$. The latter is
shown by the same arguments as in the proof of the lower bound in
(2.3).

Now, assuming $\xi \iz{N,V} >-\infty$ and letting $V\uparrow
\ttZ^\nu$, the assertions of Theorem~2.2(i), (ii) and (iii) are
derived from Theorem A.1(i), (ii) and (iii), respectively, in
(Astrauskas~2012, Appendix A) with the abbreviations $\Pi
:=\cE_V^N$, $L :=\xi \iz{N,V}$ and $r:= r \iz{N,V}$.\qed
\end{proof}

\medskip

 In the case of conditions (2.7)--(2.9) with $p=2$ (part
(iii) of the theorem), we have imposed additional restrictions
(2.10) and (2.11) to control the influence of the lower
$\xi_V$-values on the correspondence $\rrl_{K,V}\leftrightarrow
z_{K,V}$. If $\xi_{K,V}^2(\xi_{K,V}-\xi_{K+1,V})=\rm O(1)$, then
the lower $\xi_V$-values may essentially contribute to the
expansion of the eigenvalues $\rrl_{K,V}$ and, therefore, the
correspondence $\rrl_{K,V}\leftrightarrow z_{K,V}$ fails. See
Section 6.1 of the present paper where we consider the case of
i.i.d.~samples $\xi (\cdotp)$ in $V \uparrow \ttZ^\nu$ with \lk
weakly~\rk\ pronounced asymptotic peaks.

\medskip

Let $\xi(\cdotp)$ be an i.i.d.~random field with the distribution
function $1-\e^{-Q}$. We will show that, if $Q$ satisfies the
condition $Q(t)=\mro(t^{p+1})$ for $p=0,1$ and 2 (\lk heavy
tails\rk \ $\e^{-Q}$) and additional RV conditions as $t
\to\infty$, then with high probability $\xi_V$ satisfies the
assumptions of Theorem 2.2(i), (ii) and (iii), respectively, where
$K \in \ttN$ is fixed and $N=[\ve V\ve^\theta]$ for some $0<
\theta <1/2$; see Theorems 4.3(i), 4.5, 3.1 ($R=0$), 5.3 and 5.4
with $0 <\rre < \theta$. Therefore, Poisson limit theorems for the
largest eigenvalues $\rrl \iz{K,V}$ are reduced to those for
extreme values of i.i.d.~random fields $\xi(\cdotp)$ or
$\xi(\cdotp)+2\nu \rrk^2/(\xi(\cdotp)\vee 1)$ (Theorem 6.9).

\medskip

\subsection{\normalsize \textsf {Potentials with dominating single peaks: the general case} }

To simplify the proceedings, we need some notation and remarks.
For $ N \dly 2$ and $\cE^N_V$ as in (2.5), we introduce the
following function: $\wdt \xi_V(x):=0$ if $x \in \cE^N_V$, and
$\wdt \xi_V(x):=\xi_V(x)$ if $x \in V \setminus \cE^N_V$. Then
$$
\xi_V=\sum_{z \in \cE^N_V}\xi_V(z)\delta_z +\wt \xi_V,
$$
i.e., $\xi_V$ is a superposition of $\xi_V$-peaks and the noise
component $\wdt \xi_V$. To exclude the trivialities, we assume
that $\xi \iz{N,V} >-\infty$ for each $V$ (for the case when $\xi
\iz{N,V} =-\infty$ and $r\iz{N,V}>1$, see the proof of Theorem 2.2
above). For each $z \in V$, let $\wdt \rrl _V (z)$ be the
principal eigenvalue of the \lk single peak \rk \ Hamiltonian
$\rrk \Delta_V+\xi_V(z)\delta_z+\wdt \xi_V (1-\delta_z)$ on
$l^2(V)$. We associate the sites $z\iz{\tau (l),V} \in V$ with the
variational series \beq \wdt \rrl\iz{1, V}\!:=\!\wdt
\rrl(z\iz{\tau (1), V})\!\dly \!\wdt \rrl\iz{2, V}\!:=\!\wdt
\rrl(z\iz{\tau (2), V})\!\dly \!\ldots \!\dly \!\wdt \rrl\iz{\vert
V\vert , V}\!:=\! \wdt \rrl(z\iz{\tau (\vert V\vert) , V})
\eeq based on the sample $\wdt \rrl_V$; here $V=\big \{z\iz{\tau
(l),V} \colon\ 1 \mly l \mly \ve V\ve \big \}$.

\begin{thm}
Assume that there are natural numbers $1 \mly K=K_V < N=N_V < \ve
V\ve$ such that the functions $\xi_V$ satisfy condition {\rm
(2.9)} and the following conditions: \beq \lim_V
\big(\xi_{K,V}-\xi_{N,V}\big )=\infty
\eeq and \beq
 \ds \linf_V \min_{1 \mly l \mly K}\frac{\log \big (\wt \rrl_{l,V}-
 \wt\rrl_{l+1,V}\big )}{r_{N,V}\log \big(\xi_{(l+1)\wedge K,V}-\xi_{N,V}\big )} \dly 0.
\eeq Then \beq
 \ds \lsup_V \max_{1 \mly l \mly K}\frac{\log \big \ve \rrl_{l,V}-
 \wt\rrl_{l,V}\big \ve}{r_{N,V}\log \big(\xi_{l,V}-\xi_{N,V}\big )} \mly -2.
\eeq
and
$$
\lsup_V \max_{1\mly l \mly K}\max_{x \neq z_{\tau (l),V}}\frac
{\log\big|\psi(x;\rrl\iz{l,V})\big|}
{|x-z\iz{\tau(l),V}| \log (\xi\iz{l,V}-\xi\iz{N,V})} \mly -1.
$$

\end{thm}

\medskip

\begin{proof} We write
\beq \wt \cE\iz{h,V}\!:=\!\{z \in \cE_V^N: \wt\rrl_V(z) \!\dly
\!\xi\iz{N,V} \!+\!2\nu \rrk \!+\!h\} \ \ \hbox{where} \ \
h\!:=\!\frac{\xi\iz{K,V}\!-\!\xi\iz{N,V}}{2}.
\eeq By the first bound in Theorem 2.1(i), $\wdt \rrl_V(z) \dly
\xi_V(z)$ for all $z \in \cE_V^N$. This combined with (2.14) gives
\beq \big \ve\wt \cE\iz{h,V}\big \ve \dly K
\eeq for any $V\supset V_0$. Finally, according to (Astrauskas
2008, Section 2.2 and Appendix~B.2), \beq \xi_V(z) \mly \wdt
\rrl_V(z) \mly \xi_V(z)+ 2 \nu \rrk^2/h \ \ \hbox{for any} \ \ z
\in \wt \cE\iz{h,V}.
\eeq In view of (2.18) and (2.19), the assertion of Theorem 2.3
is derived similarly as in the proof of Theorem~B.3 in (Astrauskas 2008, Appendix B) with the
abbreviations $L:=\xi\iz{N,V}$, $\Pi:=\cE_V^N$~(2.5), $\wt \Pi:
=\wt \cE\iz{h,V}$ (2.17) and $r:=r\iz{N,V}$ (2.6).\qed
\end{proof}

\medskip

Note that conditions (2.14) and (2.15) of Theorem 2.3 are
substantially weaker than~(2.8) and (2.7), respectively, in
Theorem 2.2. According to (2.16) and (2.19) with $h \to \infty$ as
in (2.17), we obtain that $\rrl_{l,V}=\xi_{l,V}+ \mro(1)$
uniformly in $1 \mly l \mly K$, so that the eigenvalues
$\rrl_{l,V}$ achieve their lower bounds in (2.3) as $V \uparrow
\ttZ^\nu$.

On the other hand, from (Astrauskas~2008, Appendices A and B) we
know that, for each $z \in \wt \cE \iz{h,V}$, the eigenvalue $\wdt
\rrl_V(z)$ is the maximal solution to the equation \beq \cG_{V}
(\rrl ;\wdt \xi_V; z, z)=\frac {1}{\xi_V(z)};
\eeq here $\cG_{V} (\rrl ;\wdt \xi_V; \cdot, \cdot)$ is the
Green function of the Hamiltonian $\rrk\Delta_{V} +\wdt \xi_V$
on $l^2(V)$, so that $\cG_{V} (\rrl ;\wdt \xi_V; z, z)$ is
expanded over paths: \beq \cG_{V} (\rrl ;\wdt \xi_V; z,
z)=\sum_{\Gamma}\rrk^{\vert \Gamma\vert }\prod_{v \in V}\big
(\rrl-\widetilde \xi(v)\big )^{-n_v(\Gamma)},
\eeq where the sum $\sum_{\Gamma }$ is taken over all paths
$\Gamma : v_0:=z \to v_1 \to \cdots \to v_m :=z$ in $V$ such that
$\vert v_i-v_{i-1}\vert =1$ for each $1 \mly i \mly m$ and each
$m\in \ttN$, $n_v(\Gamma)$ denotes the number of times the path
$\Gamma $ visits the site $v \in V$, $|\Gamma |:=\sum_{v \in
V}n_v(\Gamma )-1 \dly 0$. Substituting (2.21) to the left-hand
side of (2.20) and iterating this with respect to the eigenvalue
$\rrl=\wdt \rrl_V(z)$, we obtain the explicit expansion formulas
for $\wdt \rrl_V(z)$ ($z \in \wt \cE \iz{h,V}$) presented as a
power series in the variables $\xi_V(z)$ and $\wt \xi_V(x)$ ($\ve
x-z \ve\! \dly\! 1$), in particular,
\begin{eqnarray}
\hspace{-.8cm}&&\wt \rrl_V(z)=\xi_V(z) +\rrk^{2}\sum_{\vert
x-z\vert
=1} \frac {1}{\xi_V(z)-\wt \xi_V(x)}+ \nn \\
\hspace{-.6cm}&&{}+\!\rO\bigg(\!\sum _{{\vert x-z\vert =1\atop
\vert y-z\vert =1}\atop {\vert u-z\vert \mly 2}}\!\!\frac
{1}{(\xi_V(z)\!-\!\wdt \xi_V(x))(\xi_V(z)\!-\!\wdt
\xi_V(y))(\xi_V(z)-\wdt \xi_V(u))}\!\bigg)\quad
\end{eqnarray}
as $V \uparrow \ttZ^\nu$.

\medskip

\begin{rem}
With notation at the beginning of Section 2.2, assume that there are natural numbers $1< N=N_V < \ve
V\ve$ such that the functions $\xi_V$ satisfy the  conditions: $\lim _V r\iz{N,V}=\infty$ and $\lim_V
(\xi\iz{1,V}-\xi\iz{N,V} )=\infty$. Then there are constants $\const\!_i=\const\!_i(\rrk,\nu)>0$ such that
\beq  \ve \rrl_{1,V}-
 \xi_{1,V}\ve \mly \frac{\const\!_1}{\xi_{1,V}-\xi_{N,V}} +\frac{\const\!_2}{r_{N,V}}=\mro(1)
\eeq
as $V\uparrow \ttZ^\nu$.
\end{rem}

\medskip

As shown above, limit (2.23) follows from the assumptions of
Theorem 2.3 or Theorem 2.2(ii) with $K=1$. {\it To prove (2.23)}
under the (weaker) conditions of Remark~2.4, we first observe from
Theorem 2.1(i) that $\rrl\iz{1,V} \dly \xi\iz{1,V}$ for all $V$.
Second, we apply Lemma 2.8 below with $R:=\frac{1}{3}r\iz{N,V}\to
\infty$ to see that the eigenvalue $\rrl\iz{1,V}$ is bounded from
above (with accuracy $\rO(R^{-1})$) by the maximum of local
principal eigenvalues of the single-peak Hamiltonians over all
balls in $V$ of radius $R$. Using formulas (2.20)--(2.22) for such
local principal eigenvalues, we finally obtain the bound
$\rrl_{1,V} \mly \xi_{1,V}+ \const\!_1 /(\xi_{1,V}-\xi_{N,V})+
\const/R$, as claimed. \qed

\medskip

\begin{rem}
Theorems 2.2, 2.3 and Remark 2.4 include
the condition on asymptotic sparseness of $\xi_V$-peaks: for instance,
$r\iz{N, V}\rightarrow \infty$ as $V \uparrow \ttZ^\nu$.
Molchanov and Vainberg (1998; 2000) studied the existence and location of
spectral components (pure point, absolutely continuous, etc) of the
Schr\" odinger operators $\rrk \Delta+\xi(\cdot)$ in $l^2(\ttZ^{\nu})$ with sparse
deterministic potential
$\xi(\cdotp)=\sum_{k \dly 1}a_k\delta_{z_k}(\cdotp)$, where amplitudes
$a_k$ are bounded and $\{ z_k\}$ is a rare subset of $\ttZ^\nu$,
for example, $\wt r_n :=\min\iz{l\neq n}\ve z_l-z_n\ve\rightarrow \infty$ as
$n\rightarrow \infty$.
\end{rem}
\qed

Let $\xi(\cdotp)$ be an i.i.d.~random field with the distribution
function $1-\e^{-Q}$. We will show that, if the tails $\e^{-Q}$
are heavier than the double exponential function (i.e., $\log
Q(t)=\mro(t)$) and satisfy additional regularity and continuity
conditions at infinity, then with probability one $\xi_V$
satisfies the assumptions of Theorem 2.3, where $K \in \ttN$ is
fixed and $N=[\ve V\ve^\theta]$ for some $0< \theta <1/2$; see
Theorems 3.1 ($R=0$), 4.6~($\rho=\infty$) and~4.8. Therefore,
Poisson limit theorems for the largest eigenvalues $\rrl \iz{K,V}$
are reduced to those for extremes of nonlinear functions (2.22) on
$\xi_V$ (Theorem~6.2).

\medskip

\subsection{\normalsize \textsf {Potentials with
dominating flat increasing islands of high values} }

Let $\ttB_R(z):=\{y \in V \colon\ \ve y-z \ve \mly R\}$ denote the
ball in $V$ with center $z \in V$ and radius $R \dly 0$. The
following theorem gives a simple condition on $\xi_V$ which
ensures that the largest eigenvalues $\rrl\iz{K,V}$ achieve their
upper bounds in (2.2) as $V \uparrow \ttZ^\nu$.

\begin{thm}
If \beq \lim_{R \to \infty}\lsup_V \min_{z \in V} \big
(\xi\iz{1,V}-\min_{x \in \ttB_R(z)} \xi_V(x)\big)=0,
\eeq then, for arbitrarily fixed $K \in \ttN$,
$$
\lim_V \big(\rrl\iz{K, V}-\xi\iz{K, V}\big )= 2\nu \rrk.
$$
\end{thm}

\medskip

\begin{proof}
Because of Theorem 2.1(i), we only need to show the lower limit
bound \beq \linf_V \big(\rrl\iz{K, V}-\xi\iz{1, V}\big )\dly 2\nu
\rrk.
\eeq From (2.24) we see that there exist a sequence $0 < \rre_R
\to 0$ and sites $z_V \in V$ such that $\xi_V(\cdot) \dly \xi
\iz{1,V}-\rre_R$ in $\ttB_R(z_V)$ for any $R \dly R_0$ and any $V
\supset V_0(R)$. Abbreviate $\xi_V^{(R)}(x):=\xi_V(x)$ if $x\in
\ttB_R (z_V)$, and $\xi_V^{(R)}(x):=-\infty $, otherwise. Since
$\xi_V(\cdot) \dly \xi_V^{(R)}(\cdot)$ in $V$ and
$\xi_V^{(R)}(\cdot) \dly \xi \iz{1,V}-\rre_R$ in $\ttB_R (z_V)$
for $R$ and $V$ as above, the monotonicity property of eigenvalues
implies that $\rrl_{K, V}\dly \rrl _{K,V}^{(R)}+\xi_{1,V}-\rre_R$,
where $\rrl _{K,V}^{(R)}$ is the $K$th eigenvalue of the operator
$\rrk\Delta $ on $l^2(\ttB_R (z_V))$ with zero boundary
conditions. Since $\rrl _{K,V}^{(R)}$ tends to $2\nu \rrk $
letting first $V \uparrow \ttZ^\nu$ and then $R \to \infty$
(Kirsch~2008, Section 3.1), this estimate implies (2.25), as
claimed. \qed
\end{proof}

Clearly condition (2.24) is fulfilled if and only if there are a
sequence $R_V \to \infty$ and sites $z_V \in V$ such that
$$
\lim_V \big (\xi\iz{1,V}-\min_{\ve x\! -\!z_V\ve \!\mly \!R_V}
\xi_V(x)\big)=0.
$$
From the proof of the theorem we know that the eigenvalues
$\rrl\iz{K,V}$ of the operator $\cH_V=\rrk \Delta_V+\xi_V$  in $V$
are approximated by the corresponding local eigenvalues in the
regions $\ttB_{R_V; \rm{opt}}:=\ttB_{R_V}(z_V) \subset V$ where
$\xi_V(\cdot)$ is close to $\xi\iz{1,V}$, i.e., relevant regions.

Let $\xi(\cdotp)$ be an i.i.d.~random field with the distribution
function $1-\e^{-Q}$. We will show that, if the tails $\e^{-Q}$
are lighter than the double exponential function (i.e., $t^{-1}
\log Q(t) \to \infty$) and satisfy additional RV conditions at
infinity, then with probability one $\xi_V$ satisfies the
assumption of Theorem 2.6; see Theorem 4.6 with $\rho=0$ and
Theorem 3.1(i) for any large $R$ and
$\theta(\cdot)\equiv\theta=\const$. In this case, we obtain the
second order expansion formulas for the largest eigenvalues
$\rrl\iz{K,V}$ (Theorem~6.14).

\medskip

\subsection{\normalsize \textsf {Potentials with
dominating bounded islands of high values} }

In this section, we describe a class of deterministic functions
(potential) $\xi_V :V \rightarrow [-\infty; \infty)$ for which the
asymptotic terms for the principal eigenvalue $\rrl\iz{1,V}$ ($V
\uparrow \ttZ^\nu$) fill the gap between its lower and upper
bounds in (2.2). We use the variational arguments developed by G\"
artner and Molchanov (1998). To formulate the results, we need
some abbreviations and remarks related to the variational
problems. To emphasize the dependence of $\rrl_{1,V}$ on the
sample $\xi_V$, we denote by $\rrl(\xi_V):=\rrl_{1,V}$ the
principal eigenvalue of the operator $\cH_V=\rrk \Delta_V+ \xi_V$
on $l^2(V)$. As in Section 2.4, let $\ttB_R(z) \subset V$ be the
closed ball of radius $R \dly 0$ centered at $z \in V$, and let
$\ttB_R: =\ttB_R(0)$.

Given a ball $\ttB \subset V$, let $\xi^\ttB_V(x) :=\xi_V(x)$ if
$x \in \ttB$, and $\xi^\ttB_V(x):=-\infty$, otherwise. As before,
$\cH^\ttB_V:= \rrk \Delta_V +\xi^\ttB_V$ is interpreted as an
operator with zero boundary conditions outside $\ttB$. We write
$$
\xi _V(x)= \xi\iz{1, V} +h_V(x) \quad (x \in V),
$$
where the function $h_V \mly 0$ admits the interpretation as the
shape of $\xi_V$-values close to the maximum $\xi\iz{1,V}$. Note
that \beq \rrl(\xi_V)=\xi_{1,V}+ \rrl(h_V) \quad \mbox {and} \quad
\rrl(\xi_V^\ttB)=\xi_{1,V}+ \rrl(h_V^\ttB). \eeq
For a fixed constant $0 < \rho <\infty$, we are interested in the
following supremum of $\rrl(h^\ttB)$ over $h: \ttB
\rightarrow[-\infty; 0]$:
$$
\sup \Big \{ \rrl(h^\ttB)\colon\ \sum_{x \in \ttB} \e^{h(x)/\rho}
< 1 \Big \}.
$$
This variational problem is equivalent to the corresponding
variational problem in terms of the functionals
$$
S^\ttB(p) \!:=\!\sum_{x \in \ttB} \sqrt{p(x)}\Delta \sqrt{p}(x) \
\ \mbox{where}\ \ p(y)\!=\!0\ \ \mbox{for}\ \ y \in \ttZ^\nu
\backslash \ttB,
$$
and
$$
 I^\ttB(p) :=-\sum _{x \in \ttB}p(x) \log p(x)
$$
for $p(\cdot) \in \cP (\ttB)$, the set of probability measures on
$\ttB$. More precisely, for a sequence of balls $\ttB_R \subset
\ttZ^\nu$, the following formulas hold true according to  Rayleigh--Ritz theorem and (G\"
artner and Molchanov~1998, Lemmas 2.17 and 1.10):

\beq &&\sup \Big \{ \rrl(h^{\ttB_R})\colon\ \ \sum_{x \in \ttB_R}
\e^{h(x)/\rho}
< 1 \Big \} \n\\
&&\quad {}=\sup \Big \{ \rrl(h^{\ttB_R})\colon\ \ \sum_{x \in
\ttB_R} \e^{h(x)/\rho} = 1 \Big \} \\
&&\quad {}=\sup_{p \in \cP(\ttB_R)} \Big (\rrk S^{\ttB_R}(p)- \rho
I^{\ttB_R}(p) \Big ),\n
\eeq where the right-hand side of (2.27) converges (as $R \to
\infty$) to \beq \sup_{p \in \cP(\ttZ^\nu)} \big (\rrk S(p)- \rho
I(p) \big )=:2\nu \rrk q(\rho /\rrk).
\eeq Here $S(p)$ and $I(p)$ are the corresponding functionals on
$\cP(\ttZ^\nu)$= the set of probability measures on lattice
$\ttZ^\nu$. It is easy to check that $q :\ttR_+ \rightarrow (0;1)$
is convex, strictly decreasing and surjective function; $q(0)=1$
and $q(\infty)=\lim_{\rho \to \infty}q(\rho)=0$. Moreover,
$q(\rho)=(2\rho \log \rho)^{-1}(1+\mro(1))$ as $\rho \to \infty$
(Astrauskas 2008, Proposition 2.1 and Corollary 4.5). The supremum on both
sides of (2.27) and (2.28) is attained (G\" artner and Molchanov
1998, Sections 1 and 2.4). Denote by $h^{\ttB_R}_{\rm{opt}}$ the
maximizer for the variational problem on the left-hand side of
(2.27). Then $p^{\ttB_R}_{\rm{opt}}$ is the maximizer for the
right-hand side of (2.27) if and only if $h^{\ttB_R}_{\rm{opt}}=
\rho \log p^{\ttB_R}_{\rm{opt}}$. For this and further properties
of the maximizes in (2.27) and (2.28) in the limit case
$R=\infty$, see (G\" artner and den Hollander 1999, Sections 0.3
and~0.4) and (G\" artner et al.~2007, Sections 1.3 and~3).

The following theorem tells us that, under reasonable conditions
on $\xi_V$, the principal eigenvalue $\rrl\iz{1,V}$ of the
operator $\cH_V= \rrk \Delta_V +\xi_{1,V}+ h_V$ in $V$ is
approximated (letting first $V \uparrow \ttZ^\nu$ and then $R
\to\infty$) by the local principal eigenvalue of the operator restricted to the regions
$\ttB_{R;\rm{opt}}:= \ttB_R(z_V) \subset V$ where $h_V$ is close
to $h^{\ttB_R}_{\rm{opt}}$, i.e., relevant regions with optimal
potential shape.

\begin{thm}
Given a constant $0 < \rho < \infty$ and a sequence $R \to
\infty$, assume that functions $\xi_V$ satisfy the following
conditions: \beq \lim_{R \to \infty} \lsup_V \max_{z \in V}
\sum_{y \in \ttB_R(z)} \exp \Big \{ \frac
{\xi_V(y)-\xi\iz{1,V}}{\rho} \Big \} \mly 1
\eeq and  \beq
 \ds \linf_{R \to \infty}\linf_V \max_{z \in V} \min_{y \in \ttB_R(z)}
 \big (\xi_V(y)
 -\xi\iz{1,V}-h^{\ttB_R}_{\rm{opt}}(y-z) \big ) \dly 0.
\eeq Then \beq \lim_V \big(\rrl\iz{1, V}-\xi\iz{1, V}\big )= 2\nu
\rrk q(\rho/ \rrk).
\eeq
\end{thm}

\medskip

\begin{proof} Limit (2.31) follows from the results of
(G\" artner and Molchanov 1998, Section
 2.4) under the stronger conditions on $\xi_V$
 including sparseness of clusters of $\xi_V$-extremes. To prove
 (2.31) under conditions (2.29) and (2.30), we apply the same
 arguments as in (G\" artner and Molchanov 1998, the proof of Theorem
 2.16) combined with the following lemma by Biskup and K\"onig
 (2001), which is slightly modified for the operator $\cH_V$ with periodic
 boundary conditions:

\begin{lem}
{\rm (Biskup and K\"onig 2001, Lemma 4.6).} For each $R \in \ttN$,
$V \supset V_0(R)$ and each $\xi_V$,
$$
\rrl(\xi_V) \mly \xi \iz{1,V} +\max_{z \in V} \rrl \big(
h^{\ttB_R(z)}_V \big )+ \const R^{-1}
$$
for some (universal) $\const > 0$.
\end{lem}

\medskip

We first obtain the upper bound for $\rrl(\xi_V)$. Condition
(2.29) implies that there is a sequence $0 < \rre_R \to 0$ such
that
$$
\max_{z \in V} \sum_{y \in \ttB_R(z)} \exp \big \{h_V(y)/ \rho\big
\} < \exp \big \{\rre_R /\rho \big \} \quad \mbox{for\ each} \quad
V \supset V_0(R).
$$
In view of (2.26), this estimate and Lemma 2.8 yield that
\begin{align*}
& \rrl(\xi_V) -\xi \iz{1,V}\\
& \mly \sup \Big \{ \rrl(h^{\ttB_R})\colon\ h(\cdot) \mly 0,\
\sum_{y\in \ttB_R} \e^{h(y)/\rho}
< \e^{\rre_R /\rho} \Big \} + \frac{\const}{R}\\
& \mly \sup \Big \{ \rrl(h^{\ttB_R})\colon\ h(\cdot) \mly 0,\
\sum_{y\in \ttB_R} \e^{h(y)/\rho} < 1 \Big \} +\rre_R
+\frac{\const}{R}
\end{align*}
for $V$ as above. Taking the limit as first $V \uparrow \ttZ^\nu$
and then $R \to \infty$, and using (2.27)--(2.28), we arrive at
\beq \lsup_V \big(\rrl\iz{1, V}-\xi\iz{1, V}\big ) \mly 2\nu \rrk
q(\rho/ \rrk).
\eeq By combining condition (2.30), the monotonicity property of
eigenvalues and assertions (2.27)--(2.28), similarly as in the
proof of (2.25) we obtain the lower bound
$$
\rrl(\xi_V) -\xi\iz{1, V} \dly \rrl(h^{\ttB_R}_{\rm{opt}})+\mro(1)
\to 2\nu \rrk q(\rho/ \rrk)
$$
letting first $V \uparrow \ttZ^\nu$ and then $R \to \infty$. This
and (2.32) conclude the proof of Theorem~2.7.\qed
\end{proof}

\medskip

If $\rho =\rho_V \to \infty$ in Theorem 2.7, then the \lk
relevant\rk\ regions $\ttB_{R;\rm{opt}}$ shrink to single sites
and, therefore, we are in the situation of Theorem 2.3. Meanwhile,
if $\rho =\rho_V \to 0$, then we stick to the result of Theorem
2.6.

Let $\xi(\cdotp)$ be an i.i.d.~random field with the distribution
function $1-\e^{-Q}$. It will be shown that, if the tails
$\e^{-Q}$ are the double exponential (i.e., $t^{-1} \log Q(t) \to
1/ \rho$) and satisfy additional RV conditions at $\infty$, then
with probability one $\xi_V$ satisfies the assumptions of Theorem
2.7; see Theorem 4.7 and Theorem 3.1(i) for arbitrarily large $R$
and $\theta_R(y)\approx 1- \exp \{h^{\ttB_R}_{\rm{opt}}(y)/\rho \}$
($y \in \ttB_R$). (For continuous $Q$, see Corollaries~2.7, 2.12
and 2.15 in (G\" artner and Molchanov 1998).) In this case, the
second order expansion formula for $\rrl\iz{1,V}$ holds true
(Theorem 6.19).

\medskip

\setcounter{equation}{0} \addtocounter{sct}{1}

\section{\normalsize{\textsf {CLUSTERING OF HIGH-LEVEL EXCEEDANCES OF I.I.D.~RANDOM\\
[5pt]  FIELDS} }}

\setcounter{prop}{1}

Let $\xi(x)$, $ x\in \ttZ^{\nu}$, be an i.i.d.~random field with
the cumulative hazard function~$Q$. The main task of the present
section is to investigate the almost sure asymptotic structure of
clusters (\lk islands\rk \ ) of bounded size formed by exceedances
of the sample $\xi_V$ as $V \uparrow \ttZ^\nu$. With the
abbreviations in Section 1, we also need additional notation. For
$\theta < 1$, put $L\iz{V, \theta }:=f((1-\theta )\log \vert
V\vert ) $ where $f:=Q^{\la}$. (Without loss of generality, we
write $L \iz{V,1}:=\sup \{t \colon\ Q(t-)=0\}$, so that almost
surely $\xi(x)\dly L \iz{V,1}$ for each $x$.) Let $\ttB_R(z):= \{
x \in \ttZ^\nu \colon\ \ve x -z\ve \mly R \}$, and $\ttB_R:=
\ttB_R(0)$. For fixed $R\in \ttN \cup \{ 0\}$ and a function
$\theta_R(\cdot) \colon\ \ttB_R\rightarrow (-\infty; 1]$, we
denote by $\ttV_R$ the set of balls $\ttB_R(z) \subset V$, and
$$
\cE_{V, \theta}^R:=\{ \ttB_R(z) \subset V \colon\ \xi(y)\dly
L_{V, \theta_R(y-z) }\ \ \mbox{for\ all}\ \ y \in \ttB_R(z) \},
$$
the subset of clusters of $\xi_V$-exceedances in $\ttV_R$ over the
level function $L_{V,\theta_R(\cdot)}$. We abbreviate
$$
r (\cE_{V, \theta}^R):= \min  \{ \dist (\ttB, \ttB\pr)\colon\
\ttB\!\in\! \cE_{V, \theta}^R,\ \ttB\pr\!\in\! \cE_{V, \theta}^
R,\ \ttB\!\neq\!\ttB\pr  \}\ \ \mbox{if}\ \ \vert \cE_{V, \theta}^
R\vert \dly 2,
$$
and $r(\cE_{V, \theta}^R):=\br{V}^{1/\nu }$ if $\br{\cE_{V,
\theta}^R}\mly 1$, by convention; here $\dist (\ttB, \ttB\pr )$
stands for the lattice $l^1$-distance between balls $\ttB, \ttB\pr
\subset V$. If $R=0$ and $\theta:=\theta(0)$, then $\cE_{V, \theta
}:=\cE_{V, \theta}^0$ shrinks to the subset of single
$\xi_V$-exceedances, so that
$$
r(\cE \iz{V, \theta })= \min  \{ \vert x-y\vert \colon\ x\in
\cE\iz{V, \theta },\ y\in \cE\iz{V, \theta }, \ x\neq y \}.
$$
To formulate the main result of this section, we also need the
following abbreviations
$$
\mu_R :=\sum_{y \in \ttB_R} (1-\theta_R(y)) > 0\ \ \mbox{and}\
\ \theta_{\max,R}:=\max_{y \in \ttB_R} \big \{\theta_R(y)\colon\
\theta_R(y)<1 \big \}.
$$

\begin{thm} {\rm (cf.~Theorems 2.2--2.7).} For arbitrarily fixed
$R\in \ttN \cup \{ 0\}$, the following almost sure limits hold
true.

{\rm (i)} If $\mu_R < 1$, then
$$
\linf_V \frac {\log\big\vert \cE_{V, \theta}^R\big\vert } {
\log\br{V}} \dly 1-\mu_R.
$$

{\rm (ii)} If $\mu_R < 1$ and, in addition, $Q$ satisfies the
condition \beq \lim_{t \uparrow t_Q}\frac{Q(t-)}{Q(t)}=1,
\eeq then
$$
\lim_V\frac{\log\big\vert \cE_{V, \theta}^R\big\vert
}{\log\br{V}}=1-\mu_R.
$$

{\rm (iii)} If $\mu_R> 1$ and $Q$ satisfies {\rm (3.1)}, then
$$
\lim_V \big\vert \cE_{V, \theta}^R\big\vert =0.
$$

{\rm (iv)} If $\theta_{\max,R} <\mu_R < 1$ and $Q$ satisfies {\rm
(3.1)}, then
$$
\lim_V\frac{\log r \big(\cE_{V, \theta}^
R\big)}{\log\br{V}}=\frac{2\mu_R-1}{\nu }.
$$

\end{thm}

\medskip

\begin{rem}
(a) Clearly, for arbitrary $Q$ and $\mu_R < 1$,
\begin{eqnarray*}
&&\bE\big\vert \cE_{V, \theta}^R\big\vert \quad (=\mbox{the\ mean\
number\ of\ clusters\ of\
exceedances}) \\
&&=\vert \ttV_R\vert \prod_{y \in \ttB_R} \bP \big( \xi(y) \dly L_{V, \theta_R(y) }\big ) \\
&&\dly \const \vert V\vert \exp\big \{-\sum_{y \in \ttB_R} Q(L
_{V, \theta_R(y)}-)\big \} \dly \const\vert V\vert ^{1-\mu_R}
\to \infty,
\end{eqnarray*}
according to Lemma A.11(iii) in Appendix.

(b) On the other hand, by Lemma A.11(iii), condition (3.1) implies
the asymptotic formula \beq Q(f(s))=s +\mro(s) \quad \mbox{as}
\quad s \to \infty,
\eeq which in turn yields, for $\mu_R \dly 0$, the upper bound
$$
\log \bE\vert \cE_{V, \theta}^R\vert \mly (
1-\mu_R+\mro(1))\log\vert V\vert
$$
as $\ve V \ve \to \infty$.
\end{rem}

\medskip

\begin{rem} {\rm (see part (iv)).} If $\theta_{\max,R} <\mu_R < 1$, then $1/2 <\mu_R < 1$.
\end{rem}

\medskip

\begin{rem}
In the case $R=0$ and $0 < \theta <1/2$, i.e., single rare
$\xi_V$-peaks, Theorem 3.1 was proved by Astrauskas (2001). For
the Gaussian random field $\xi(\cdot)$ with correlated values, the
case $R=0$ was studied by Astrauskas (2003). Here the results
depend slightly on the correlation function of $\xi( \cdot)$.
Finally, assertion (i) generalizes Corollary 2.15(b) in (G\"artner
and Molchanov 1998) where the continuity of $Q$ is assumed.
\end{rem}

\medskip

\medskip

{\it Proof of Theorem 3.1.} To simplify the proof, we assume
throughout that $\theta_R(\cdot)\equiv \theta(\cdot)<1$ in
$\ttB_R$. The general case is treated similarly.

(i) We denote by $\wdt \ttV_R\subset \ttV_R$ the maximal subset of
nonintersecting balls $\ttB_R(z)$ in $V$, so that $\br{\wdt
\ttV_R}\asymp \br{V}$. The claimed bound is proved by estimating
$\big \ve \cE_{V,\theta}^R\cap \wdt \ttV_R\big \ve $ similarly as
in the proof of Theorem 1 in Astrauskas (2001), where the
exceedances $\{\xi(x)\dly L_{V,\theta} \}$ ($x\in V$) are
replaced by mutually independent (multiple) exceedances $\{
\xi(\cdot)\dly L_{V, \theta(\cdot-z) }$ in $\ttB_R(z) \}$
($\ttB_R(z)\in \wdt \ttV_R$). In particular, if $Q(t_Q-)=\infty$,
we obtain that, for any $-1 <\delta <0$ and $V \uparrow \ttZ^\nu$,
\begin{eqnarray*}
&&\bP\lf (\big\ve \cE_{V,\theta}^R\cap \wdt \ttV_R\big\ve \mly
(1+\delta) \bE \big\ve \cE_{V,\theta}^R\cap
\wdt \ttV_R\big\ve \rg )\\[6pt]
&&\mly \exp \big \{ -\const (\delta)\cdot\bE \big\ve \cE_{V,\theta}^
R\cap \wdt \ttV_R\big\ve (1+\mro(1))\big \}
\end{eqnarray*}
for some $\const(\delta) >0$. Since the right-hand side is
summable over $V$ according to the assertion of Remark 3.2(a), we
conclude the proof of (i) by using the Borel-Cantelli lemma.

(ii) We only need to estimate $\ve \cE_{V,\theta}^R\ve $ from
above. Fix a function $\theta\pr(\cdot) \colon\ \ttB_R\rightarrow
(-\infty; 1)$ such that $\theta \pr (\cdot)> \theta(\cdot)$ in
$\ttB_R$, and pick a constant $\delta
>1-\mu \pr$ where $\mu \pr :=\sum_{y \in \ttB_R}(1-\theta \pr(y) )$.
We then apply Chebyshev's inequality and the assertion of Remark
3.2(b) to find that, for any $V\supset V_0$, \beq \bP \lf (\big\ve
\cE_{V,\theta \pr}^R\big\vert
>\br{V}^\delta \rg )\mly \bE \big \vert \cE_{V,\theta \pr}^R\big \vert
\br{V}^{-\delta }\mly \br{V}^{-\mconst }
\eeq where $\const =\const(\theta \pr(\cdot) ,\delta )>0$. Choose
a subsequence $\{ V(l)\colon \ l \in \ttN \} \subset \{V\}$ such
that \beq V(l)\ \ \mbox{monotonously\ increases\ and}\ \ \vert
V(l)\vert =2^l(1+\mro (1)) \ \ \mbox{as}\ \ l\to \infty .
\eeq Since the right-hand side of (3.3) is summable over the
subsequence $\{ V(l)\}$, the Borel--Cantelli lemma implies that
almost surely $\big\ve \cE_{V(l),\theta \pr}^R\big\vert \mly \ve
V(l)\ve ^\delta$ for all $l \dly l_0( \rro)$. Because of the
monotonicity of $\cE_{V,\theta}^R$ in $L _{V, \theta(\cdot)}$,
we obtain that with probability 1, for any $V$ such that $V(l-1)
\subset V \subseteq V(l)$ and any $l \dly l_0(\rro; \theta(\cdot),
\theta \pr(\cdot))$, the set $\cE_{V,\theta}^R$ is contained in
$\cE_{V(l),\theta \pr}^R$, therefore,
$$
\big \vert \cE_{V,\theta}^R\big \vert \mly \big \vert
\cE_{V(l),\theta \pr}^R\big \vert \mly \ve V(l) \ve^\delta \mly
\const \ve V \ve ^\delta .
$$
Since $\theta \pr (\cdot)> \theta(\cdot)$ and $\delta
>1-\mu \pr$ are chosen arbitrarily, this estimate yields the upper
limit bound for $\log \big\ve \cE_{V,\theta}^R\big\vert$, as
claimed.

As in part (ii), it suffices to prove the assertions of
(iii)--(iv) for the subsequence $\{ V(l) \}$ (3.4) instead of $\{
V \}$.

(iii) We note that $ \cE_{V,\theta}^R \neq \tus$ if and only if
there exists $\ttB_R(z) \subset V$ such that $\xi(\cdot) \dly
L_{V, \theta(\cdot-z)}$ in $\ttB_R(z)$. According to the
assertion of Remark 3.2(b), the probability of the last event does
not exceed $\bE \big \vert \cE_{V,\theta}^R\big \vert \mly
\br{V}^{-\rho }$ for some $0 < \rho <-1+\mu$. Therefore, $\bP \lf
(\cE_{V,\theta}^R \neq \tus \rg ) \mly \ve V\ve^{-\rho}$. Since
the latter is summable over $\{ V(l)\}$ (3.4), the Borel-Cantelli
lemma yields that almost surely $\cE_{V(l),\theta}^R = \tus$ for
all $l \dly l_0(\rro)$, as claimed.

(iv) With $\wdt \ttV_R\subset \ttV_R$ defined in part (i), the
almost sure upper bound for $r\big(\cE_{V,\theta}^R\cap \wdt
\ttV_R\big)\dly r\big(\cE_{V,\theta}^R\big)$ is derived similarly
as in the proof of Theorem 2 of Astrauskas (2001) where the
exceedances $\{\xi(x)\dly L_{V,\theta} \}$ ($x\in V$) are
replaced by mutually independent (multiple) exceedances $\{
\xi(\cdot)\dly L_{V, \theta(\cdot-z) }$ in $\ttB_R(z) \}$
($\ttB_R(z)\in \wdt \ttV_R$).

To obtain the lower bound for $r\big(\cE_{V,\theta}^R\big)$, we
first note that the event $\big \{r\big (\cE_{V,\theta}^R\big)=0,
\\ \big\ve \cE_{V,\theta}^R\big\vert \dly 2 \big \}$ implies that
there exists $\ttB_R(z) \subset V$ such that $\xi(\cdot) \dly
L_{V, \theta(\cdot-z)}$ in $\ttB_R(z)$ and $\xi(y) \dly L_{V,
\theta_{\max}}$ for some $y \in \big (\ttB_{3R}(z)
\setminus\ttB_R(z)\big )\cap V$. Therefore, as in the proof
of~(iii) we obtain that, for fixed $y \in \ttZ^{\nu}\setminus
\ttB_R$ and for any $V \supset V_0$,
\begin{eqnarray*}
&&\bP \lf (r\big (\cE_{V,\theta}^R\big)=0,
\ \big\ve \cE_{V,\theta}^R\big\vert \dly 2 \rg ) \\
&&\mly  \const \ve V\ve \bP \lf (\xi(\cdot) \dly L_{V,
\theta(\cdot)}\ \rm{in}\ \ttB_{\it {R}},\ \xi(\it {y}) \dly \it {L}_{V,
\theta_{\max}} \rg) \mly \br{V}^{-\rho}
\end{eqnarray*}
for some $0 < \rho < \mu -\theta_{\max}$. Second, similarly as in
the proof of part~(i), we find that $\bP \lf (\big \ve
\cE_{V,\theta}^{R}\big \ve <2 \rg ) \mly \br{V}^{-{\rm
const}}$ for any $V \supset V_0$ and some $\const > 0$.
Summarizing these bounds and picking $0 < \rre < (2\mu -1)/\nu$
arbitrarily, we get that, for any $V \supset V_0$,
\begin{eqnarray}
&&\hspace{-.7cm}\bP \lf
(r\big(\cE_{V,\theta}^{R}\big)< \ve V \ve^\rre \rg ) \n \\
&&\hspace{-.7cm}\mly \bP \lf ( 1 \mly
r\big(\cE_{V,\theta}^{R}\big) < \ve V \ve^\rre,
\big\ve\cE_{V,\theta}^R\big \ve \dly 2 \rg )\! +\!\ve V
\ve^{-\mconst_1}\! \mly \! \ve V \ve^{-\mconst_2}
\end{eqnarray}
for some $\const_i >0$, where the last probability is estimated
similarly as in the proof of Theorem 2 of (Astrauskas 2001) with
mutually independent (multiple) exceedances $\{ \xi(\cdot)\dly
L_{V, \theta(\cdot-z) }$ in $\ttB_R(z) \}$ instead of
$\{\xi(x)\dly L_{V,\theta} \}$. Since the right-hand side
of~(3.5) is again summable over $\{V(l)\}$ (3.4), we conclude from
the Borel--Cantelli lemma that almost surely
$r\big(\cE_{V(l),\theta}^ R\big) \dly \ve V(l) \ve ^\rre$ for any
$l \dly l_0(\rro; \rre)$, as claimed. This completes the proof of
Theorem 3.1. \qed

\medskip

\begin{rem}
By the same arguments as in the proof above, the assertions of
Theorem 3.1 are extended to the following class of high-level
exceedances:

For fixed $R\in \ttN \cup \{ 0\}$, we denote by
$\ttS_R$ the set of all subsets $U \subset \ttZ^\nu$, the diameter
of which does not exceed $R$. Let $\ttV_R:=\lf \{U \in \ttS_R
\colon\ U \subset V \rg \}$. For a fixed set of functions
$\Theta_R:=\lf \{ \theta_{U,R}(\cdot) \in (-\infty;1)^{\ve U\ve}
\colon\ U \in \ttS_R \rg\}$, let $\cE_{V, \Theta}^R \subset
\ttV_R$ be the subset of elements $U \in \ttV_R$ such that
$\xi(\cdot)\dly L_{V, \theta_{U,R}(\cdot)}$ in $U$. I.e.,
$\cE_{V, \Theta}^R$ consists of clusters of exceedances in
$\ttV_R$ over level functions $L\iz{V,\Theta}$. Denote by
$r\lf(\cE_{V, \Theta}^R\rg)$ the minimum distance among elements
$U,U\pr \in \cE_{V, \Theta}^R$, $U \neq U\pr$. Finally, let $\mu_R
:=\sum_{y \in U} (1-\theta_{U,R}(y))$ be a positive constant
independent of $U \in \ttS_R$, and write
$\theta_{\max,R}:=\sup_{U\in \ttS_R}\max_{y \in U}
\theta_{U,R}(y)$. With these notation for $\cE_{V, \Theta}^R$ and
$r \lf (\cE_{V, \Theta}^R\rg)$, the almost sure assertions
(i)--(iv) of Theorem~3.1 hold true.
\end{rem}

\medskip

\setcounter{equation}{0} \addtocounter{sct}{1}

\section{\normalsize{\textsf {SPACINGS OF ORDER STATISTICS OF I.I.D.~RANDOM FIELDS} }}

\setcounter{prop}{1}

\subsection{\normalsize{\textsf {Spacings of consecutive order statistics} }}

We first formulate the results for the exponential order
statistics $\eta \iz{K, V}$ and their spacings, which are then
transferred to $\xi \iz{K, V}=f(\eta \iz{K, V})$ under appropriate
conditions for $f$.

Note that the random variables
$$
\eta \iz{1, V}-\eta \iz{2, V},\ldots ,(\br{V}-1) (\eta
\iz{\br{V}-1, V}-\eta \iz{\br{V}, V} ), \br{V}\eta \iz{\br{V}, V}
$$
are mutually independent exponentially distributed with mean 1;
see, e.g., (Shorack and Wellner 1986, pp. 336). This property
immediately implies the first assertion of the following lemma.

\begin{lem}
{\rm (i)} For fixed $K\in \ttN$,
\begin{eqnarray*}
&&\lim _V \bP\! (\eta \iz{1, V}-\eta \iz{2, V}\!>\!t_1,\ldots
,\eta \iz{K-1,
V}-\eta \iz{K, V}\!>\!t\iz{K-1}, \eta \iz{K, V}-\log \br{V}\!>\!t ) \\
&&\quad {}=\bigg (\prod _{l=1}^{K-1}\e^{-lt_l}\bigg
)\frac{1}{(K-1)!}\int_t^\infty \exp\lf \{ -Ks-\e ^{-s}\rg \} \dd s
\end{eqnarray*}
for all $t_l\dly 0$ ($1\mly l\mly K-1$) and all $t\in \ttR $.

{\rm (ii)} For an arbitrary sequence $\{K\iz{V}\}$ such that
$1\mly K\iz{V}\mly \br{V}$,
$$
\lsup _V\sqrt{K\iz{V}}\max _{K_V\mly k\mly \br{V}}\lf \vert \eta
\iz{k, V}-\log \frac {\br{V}}{k}\rg\vert <\infty \quad \hbox{in\
probability}.
$$
\end{lem}

\begin{proof}
Let us show (ii). Write $\eta \iz{\br{V}+1, V}:=0$. By
Kolmogorov's inequality in (Shorack and Wellner 1986, pp. 843), we
have that
\begin{eqnarray*}
&&\bP\Bigg (\max _{K_V\mly k\mly \br{V}}\bigg \vert \sum
_{l=k}^{\br{V}}\lf (\eta \iz{l, V}-\eta \iz{l+1, V}-\frac
{1}{l}\rg )\bigg \ve >\lf (\frac
{C}{K\iz{V}}\rg )^{1/2}\Bigg ) \\
&&\quad {}\mly \frac {K\iz{V}}{C}\sum _{l=K_V}^{\br {V}} \bE\lf
(\eta \iz{l, V}-\eta \iz{l+1, V}-\frac {1}{l}\rg )^2\mly \frac
{2}{C}
\end{eqnarray*}
for any $C>C_0$ and any $V\supset V_0(C)$. Combining this bound
with the following simple estimate
$$
\max _{K_V\mly k\mly \br{V}}\bigg (\sum _{l=k}^{\br {V}}\frac
{1}{l}-\log \frac {\br {V}}{k}\bigg )\mly \frac {1}{K\iz{V}}\quad
(V\supset V_0),
$$
we obtain the claimed assertion of (ii).\qed
\end{proof}

\smallskip

The almost sure asymptotic behavior of the random variables $\eta
\iz{K, V}$ and $\eta \iz{K, V}-\eta \iz{K+1, V}$  ($\br{V}\to
\infty $) is more intricate.

\begin{lem} For any fixed constants $K\in \ttN $ and $m\in \ttN \setminus \{ 1\}$,
the following almost sure limits hold true.

{\rm (i)} \vskip-2.5\baselineskip
\begin{eqnarray*}
&&\linf_V\frac{\log (\eta \iz{K, V}-\eta \iz{K+1, V} )+\sum
_{i=2}^{m-1}\log_i\br{V}}{\log _m\br{V}}=-1, \\
\noalign{{\rm (ii)} \vskip-\baselineskip} &&\lsup_V\frac{\eta
\iz{K, V}-\eta \iz{K+1, V} -K^{-1} \sum
_{i=2}^{m-1}\log_i\br{V}}{\log
_m\br{V}}=\frac{1}{K} \\
\noalign{\noindent \mbox{and}} \noalign{{\rm
(iii)}\vskip-\baselineskip} &&\lsup _V \max_{1 \mly l \mly \ve V
\ve} \lf\vert \eta \iz{l, V}-\log \frac {\br{V}}{l} \rg\vert \frac
{1}{\log \log \br{V}} = 1;
\end{eqnarray*}
here $\sum _2^1\ldots :=0$.
\end{lem}

\begin{proof} Assertions (i) and (ii) follow from more general
results for exponential spacings in (Astrauskas 2006, Corollary
12). Assertion (iii) follows from the corresponding strong limits
for the uniform order statistics $\zeta\iz{ \ve V \ve -k +1, V}$
(Shorack and Wellner 1986, pp. 408 and pp. 420--424) via
transformation $\eta \iz{k,V} = -\log \zeta\iz{ \ve V \ve -k +1,
V}$ ($1 \mly k \mly \ve V \ve$).

\qed
\end{proof}

\smallskip

We now turn to the case $\xi \iz{k,V}=f(\eta \iz{k,V})$. For $p
\dly 0$, we denote by $A\Pi^{p}_{\infty}$ the class of functions
$f:=Q^{\la}$ that satisfy \beq\lim_{s \to \infty} f(s)^p
\big(f(s+c)-f(s) \big)
 =\infty\ \ \mbox{for\ any}\ \ c > 0,
\eeq and by $A\Pi^{p}_{0}$ the class of functions $f$ that satisfy
\beq\lim_{s \to \infty} f(s)^p \big(f(s+c)-f(s) \big)
 =0\ \ \mbox{for\ any}\ \ c > 0,
\eeq and, finally, $\rO A\Pi^p$ stands for the class of $f$
satisfying \beq
 f(s)^p
\big(f(s+c)-f(s) \big)
 \asymp 1\ \ \mbox{as}\ \ s \to \infty,\ \ \mbox{for\ any}\ \ c >
 0.
\eeq

 We see that, if $f$ is in $A\Pi^{p}_{\infty}$ or $\rO A
\Pi^p$, then the right endpoint $t_Q$ is infinity or,
equivalently, $f(s) \to \infty$ as $s \to \infty $. Of course,
$A\Pi^{p}_{0}$ includes the trivial case of finite $t_Q > 0$. The
characterization of $A\Pi^{p}_{\infty}$, $A\Pi^{p}_{0}$ and $\rO
A\Pi^p$ is given in Lemmas A.6, A.7 and A.8 of Appendix A
respectively. In particular, the functions $f:=Q^\la \in \rO
A\Pi^p$ are associated with Weibull type distributions
$1-\e^{-Q}$, where $Q(t) \asymp t^{p+1}$ as $t \to \infty$.

\begin{thm} {\rm (Cf.~(2.7) and (Astrauskas 2012; 2013)).}
For fixed natural $K>l\dly 1$ and real $p\dly 0$, we have the
following limits in probability.

{\rm (i)} If $f \in A\Pi^{p}_{\infty}$, then \beq \lim_V\xi_{K,
V}^p(\xi\iz{l, V}-\xi\iz{K, V})=\infty.
\eeq

{\rm (ii)} If $f \in A\Pi^{p}_{0}$, then $ \lim_V\xi_{K,
V}^p(\xi\iz{l, V}-\xi\iz{K, V})=0$.

{\rm (iii)} If $f \in \rO A\Pi^p$, then $ \xi_{K, V}^p(\xi\iz{l,
V}-\xi\iz{K, V})\asymp 1$ as $\br{V}\to \infty$.
\end{thm}

\medskip

\begin{proof}
Using notation (1.22)--(1.24), rewrite the left-hand side of
(4.4) in the form
$$
f(\eta\iz{K, V})^p\lf (f (\eta \iz{K, V}+(\eta \iz{l, V}-\eta _{K,
V}) )-f(\eta \iz{K, V})\rg ).
$$
The claimed assertions follow by applying Lemma 4.1(i).\qed
\end{proof}

\medskip

To obtain these limits with probability 1, we need the stronger
conditions for $f$. Let us abbreviate
$$
d_{m,\gamma}(s):=s\Big (\prod_{i=1}^{m-1} \log_{i}s \Big )(\log_m
s)^{1+\rrg}\quad (s \dly s_0).
$$
\begin{thm} {\rm (Cf.~(2.7)).}
For fixed constants $K \in \ttN$ and $p\dly 0$, the following
almost sure limits hold true.

{\rm (i)} If $\lim_{s \to \infty}f(s)^p
\big(f(s+1/d_{m,\gamma}(s))-f(s) \big) = \infty$ for some $m \in
\ttN$ and $\rrg
> 0$, then
$$
\lim_V \xi_{K+1, V}^p(\xi\iz{K, V}-\xi\iz{K+1, V})= \infty.
$$

{\rm (ii)} If $\lim_{s \to \infty}f(s)^p \big(f(s+K^{-1} \log
d_{m,\gamma}(s))-f(s) \big) = 0$  for some $m \in \ttN$ and $\rrg
> 0$, then
$$
\lim_V \xi_{K+1, V}^p(\xi\iz{K, V}-\xi\iz{K+1, V})= 0.
$$

{\rm (iii)} If $\lim_{s \to \infty}\big (f(s+\log s)-f(s)\big)=
0$, then
$$
\lim_V \big( \xi\iz{K,V}-f(\log \br{V})\big)=0.
$$

\end{thm}

\medskip

\begin{proof}
Assertions (i)--(ii) follow by the same arguments as in the proof
of Theorem~4.3, where one applies Lemma 4.2 instead of Lemma
4.1(i). Assertion (iii) follows from Lemma 4.2 (iii).\qed
\end{proof}

\medskip

\subsection{\normalsize{\textsf {Spacings of intermediate order statistics} }}

We denote by $PI_{<2}$  the class of functions $f:=Q^{\leftarrow}$
satisfying the condition \beq \lsup_{s \to
\infty}\frac{f((1-\rre)s)}{f(s)}< 1 \quad \mbox {for\ some} \quad
0 < \rre <1 /2.
\eeq Class  (4.5) is characterized in Lemma A.13 (Appendix A).

\begin{thm} {\rm (Cf.~(2.8) ).}
Assume that $f \in PI_{<2}$ {\rm (4.5)}. Then for fixed $K \in
\ttN$ and $\theta > \rre$, almost surely \beq \lsup_V \xi\iz{[\ve
V\ve^\theta],V} \big / \xi \iz{K, V} < \const <1.
\eeq
\end{thm}

\medskip

\begin{proof} By Lemma 4.2(iii), with probability one the random variable $\xi \iz{K, V}=
f(\eta \iz{K, V})$ is bounded from below by $f( \log \ve V\ve
-2\log \log \ve V\ve)$ and $\xi\iz{[\ve V\ve^\theta],V}$ is
bounded from above by  $f( (1-\theta) \log \ve V\ve +2\log \log
\ve V\ve)$ for each $V\supset V_0(\rro)$. Substituting these
bounds into the left-hand side of (4.6) and using (4.5), we obtain
the claimed assertion.

\qed
\end{proof}

\medskip

We denote by $RV_\rrr$ the class of nondecreasing functions
$g:\ttR_+ \rightarrow \ttR_+$ such that, for any $c > 1$, $\lim_{s
\to \infty} g(cs)/g(s) = c^\rho$. I.e., $g$ is regularly varying
at infinity with index $0 \mly\rho \mly \infty$. The case $\rho =
\infty$ (resp., $\rho = 0$) indicates a rapid variation (resp.,
slow variation) of the function $g$. See Lemma A.3 in Appendix A
for a summary of the well-known properties of the class $RV_\rho$.

\begin{thm}{\rm (Cf.~(2.14) and (2.24) ).} For some $0\mly\rho\mly\infty$,
assume that $\e^f \in RV_\rho$. Then, for all constants $0 \mly
\rre < \theta <1$,  almost surely
$$
\lim_V \big(\xi\iz{[\ve V\ve^\rre],V} -\xi\iz{[\ve V\ve^\theta],V}
\big )=\rho \log \frac{1-\rre}{1-\theta}.
$$
\end{thm}

\medskip

\begin{proof} To prove this assertion, use Lemma 4.2(iii) and Lemma A.3(ii)
similarly as in the proof of Theorem 4.5.\qed
\end{proof}

\medskip

The following statement is closely related to the result of
Theorem 4.6 with $0 <\rho <\infty$. As in Section 2.4, let
$\ttB_R(z)$ denote the closed ball in $V$ with the center $z \in
V$ and the radius $R \dly 0$.

\medskip

\begin{thm} {\rm (Cf.~(2.29) and (G\" artner and Molchanov 1998)).}
For some $0 <\rho <\infty$, assume that $\e^f \in RV_\rho$. Then,
for any fixed $R \in \ttN$, almost surely
$$
\lsup_V \max_{z \in V}\sum _{y \in \ttB_R(z)} \exp \big
\{\big(\xi(y) -\xi\iz{1,V} \big )/\rho \big \} \mly 1.
$$
\end{thm}

\medskip

\begin{proof} For continuous $Q$, this assertion is a straightforward
consequence of Corollary 2.12 in (G\" artner and Molchanov 1998)
and Theorem 4.4(iii) above. (In view of Lemma A.3(ii), the
condition of Theorem 4.4(iii) follows from the assumption of
Theorem 4.7). If the continuity condition on $Q$ is dropped, one
applies slightly mo\~di\~fied arguments based on the technique of
function inversion, e.g., Lemma A.11 in Appendix~A.\qed
\end{proof}

Recall that the conditions of Theorems 4.6 and 4.7 are discussed
in Lemma A.3. In particular, the assumption of Theorem 4.7 implies
that $\log Q(t)=t/\rho +\mro(t)$ as $t \to \infty$, i.e., the
double exponential tails $\e^{-Q}$.

\medskip

\subsection{\normalsize{\textsf {Minimum of spacings} }}

We first recall some notation from Section 2.3. For fixed $0 <
\theta < 1/2$, we write $\wt \xi(x):=\xi(x)$ if $\xi(x) <
f((1-\theta)\log \ve V \ve)$, and $\wt \xi(x):=0$, otherwise. For
any $z \in V$, let $\wt \rrl(z)$ denote the principal eigenvalue
of the \lk single peak\rk \ Hamiltonian $\rrk\Delta_{V}
+\xi(z)\delta \iz{z} +\wt\xi_V (1-\delta \iz{z})$ in $l^2(V)$. As
in (2.13), let $\wt \rrl_{K,V}$ denote the $K$th extreme order
statistics of the random field $\wt \rrl_V$. For $\rrk=0$ and $0 <
\rre <\theta$, we know from Theorem 3.1(i)($R=0$) that with
probability one $\wt \rrl_{k,V}\equiv \xi_{k,V}$ for all $1 \mly k
\mly \ve V \ve^{\rre}$ and all large $V$. For any $\rrk \dly 0$,
we are interested in the asymptotic behavior of the minimum of the
gaps $\wt \rrl _{k, V}-\wt \rrl _{k+1, V}$ ($1 \mly k \mly \ve V
\ve^\rre$) defined by
$$
S\iz{V, \rre}:= \min \big\{ \wt \rrl _{k, V}-\wt \rrl_{k+1, V}
\colon\ 1 \mly k \mly \ve V \ve^\rre \big\}.
$$

Given a constant $\mu>0$, we say that the function
$F:\ttR\rightarrow\ttR$ is {\it log-H\" older con\-ti\-nu\-ous of
order $\mu >0 $ at infinity}, if $F$ satisfies the following
condition:
\begin{equation}
\begin{split}
&\lf \ve F(t+s)-F(t-s)\rg \ve \vert \log s\vert ^\mu =\rO(1)\\
&\mbox{as}\ \ t\to \infty\ \mbox{and}\ \ s\downarrow 0\ \
\mbox{simultaneously.}\mbox{\hspace{2.5cm}}
\end{split}
\end{equation}

\begin{thm} {\rm (Cf.~(2.15)).} Let $t_Q=\infty$, $\rrk \dly 0$
and $0 <\rre < \theta <1/2$, and assume that the distribution
tails $\e^{-Q}$ are log-H\"older continuous of order $\mu>0$ at
infinity. For $\rrk > 0$, assume additionally that $\e^f \in
RV_\infty$. Then almost surely
$$
 \lsup _V \frac{\log\{
-\log (S\iz{V, \rre }\land 1 )\} } {\log\br{V}}\mly \frac{1+\rre
}{\mu}.
$$
\end{thm}

\begin{proof} The assertion follows from Lemmas
3.5 and 4.3  in (Astrauskas 2008) and Theorem 3.1(ii) above, where
$R=0$ and $0< \theta <1/2$.\qed
\end{proof}

\medskip

In (Astrauskas 2003), the results of Theorem 4.8 are extended to the
Gaussian random fields with correlated values.

We end this section with some generalization of Theorem 4.3(iii)
for the functions $f \in \rO A \Pi^p$ (4.3) associated with
Weibull type distributions $1\!-\!\e^{-Q}$, where $Q(t) \asymp
t^{p+1}$ as $t \to \infty$.

\begin{thm} {\rm (Cf.~Lemma 4.2 in (Astrauskas 2013)).}
For some $p \dly 0$, assume that $f \in \rO A \Pi^p$ {\rm (4.3)}.
Then, for arbitrarily fixed constants $K \in \ttN$, $0<\rre <1$
and any sequence $\{ n_V\} \subset \ttN$ such that
$n\iz{V}=\rO(\br{V}^\rre )$, we have the following limits in
probability:

{\rm (i)} \vskip-2.5\baselineskip
\begin{eqnarray*}
&&\xi\iz{n_V, V}\asymp \lf (\log\br{V}\rg )^{1/(p+1)}\ \
\mbox{as}\  \ \br{V}\to \infty , \\
\noalign{\noindent \mbox{and}} \noalign{{\rm
(ii)}\vskip-\baselineskip} &&0<\linf _V\min_{K+1\mly l\mly
\br{V}^\rre }\xi_{l, V}^p(\xi
\iz{K, V}-\xi\iz{l, V})\frac{1}{\log l} \\
&&\quad {}\mly \lsup _V\max_{K+1\mly l\mly \br{V}^\rre }\xi_{l,
V}^p(\xi \iz{K, V}-\xi\iz{l, V})\frac{1}{\log l} <\infty .
\end{eqnarray*}
\end{thm}

\medskip

\begin{proof} Assertion (i) follows from a combination of the formula
$\xi \iz{k,V}=f(\eta \iz{k,V})$, Lemma 4.1(ii) and the limit $f(s)
\asymp s^{1/(p+1)}$ as $s \to \infty$ (the latter follows from
Lemma A.8(iii) with $a(\cdot) \equiv \const$). Assertion (ii) is
shown by combining the formula $\xi \iz{k,V}=f(\eta \iz{k,V})$ and
Lemmas 4.1(ii), A.8(iii) similarly as in the proof of Theorems~4.3
and~4.5 above.\qed
\end{proof}

\smallskip

\setcounter{equation}{0} \addtocounter{sct}{1}

\section{\normalsize{\textsf {NEIGHBORING EFFECTS FOR EXTREMES OF I.I.D.~RANDOM\\
[5pt] FIELDS} }}

\setcounter{prop}{1}

We finally study the asymptotic properties of $\eta \iz{V}$-values
neighboring to $\eta \iz{V}$-peaks. It is then straightforward to
extend the results for $\eta (\cdotp )$ to $\xi(\cdotp )=f(\eta
(\cdotp ))$.

The following lemma tells us that, for fixed $y\neq 0$ and for
small $\rre
>0$, asymptotic properties of the random variables $\eta (z\iz{k,
V}+y)$ ($1\mly k \mly\vert V\vert^{\rre}$) and their extremes are
the same as in the case of exponential i.i.d.~random variables.

\begin{lem}
For fixed $y\in \ttZ^\nu \backslash \{0\}$, $0 < \rre < 1/2$ and a
sequence of integers $K:= K_V =\rO (\ve V \ve^\rre)$, the
following assertions
hold true.\\

{\rm (i)} $\lim_V \bP \big (\eta(z\iz{K, V}+y) > t \big )
=\e^{-t}$ for
all $ t \dly 0$.\\

{\rm (ii)} If, in addition, $K:=K\iz{V}\to \infty $, then
$$
\lim_V \bP \big ( \max_{1\mly l\mly K}\eta (z\iz{l, V}+y)-\log K
\mly t \big)=\exp \{-\e^{-t} \} \quad \mbox {for\ all} \quad t \in
\ttR.
$$

{\rm (iii)} \*
$$
\lim_{M \to \infty} \lsup_V \Big \ve \max_{M \mly l \mly \vert
V\vert^\rre} \frac {\eta(z\iz{l,V}+y)}{\log l} - 1 \Big \ve=0\ \
\mbox{in\ probability}.
$$
\end{lem}

\begin{proof}
Here and in the sequel, we need the following key statement (which
is frequently used in (Astrauskas 2013) as well).

\begin{lem}
Fix a finite subset $U\subset \ttZ^\nu\backslash \{0\}$, $U\neq
\tus $, and a sequence of nonrandom real functions
$\{D_l(t\iz{U})\colon\ t\iz{U}\in \ttR^{\br{U}}\}$ ($l \in \ttN$).
Abbreviate $\eta (z; l)\!:=\!D_l(\{\eta (z+x)\colon x\!\in \!U\})$
for $z\in \ttZ^\nu $ and $l \in \ttN$. Finally, pick a sequence of
integers $K:=K\iz{V}=\rO(\br{V}^\rre )$ for some $0<\rre
<\frac{1}{2}$.  Then, for any $V$ and any $t\in \ttR$,
$$
\bigg \ve \bP\lf(\max_{1\mly l\mly K}\eta (z\iz{l, V}; l)\mly t\rg
)-\prod _{l=1}^{K}\bP(\eta (0; l)\mly t)\bigg \ve \mly 3
\br{V}^{-{\rm const}},
$$
where $\const>0$ does not depend on $V$ and $t$.
\end{lem}

Now, part (i) of Lemma 5.1 follows from Lemma 5.2 with $U:=\{y\}$,
where $D_K(t\iz{U})\equiv t_y$ and $D_l(t\iz{U})\equiv 0$ for
$l\neq K $. Part (ii) follows from Lemma 5.2, where
$D_l(t\iz{U})\equiv t_y$ ($l\in \ttN$), combined with Lemma
4.1(i). Finally, by Lemma 5.2 with $D_l(t\iz{U})\equiv t_y / \log
l$ ($l\dly 2$), we derive that, for any small $\delta >0$,
$$
\lsup_V \bP \Big ( \max_{M \mly l \mly \ve V\ve^\rre}\frac
{\eta(z\iz{l,V}+y)}{\log l} > 1+\delta \Big ) \mly
\sum_{l=M}^{\infty} \e^{-(1+\delta) \log l} \to 0
$$
and
$$
\lsup_V \bP \Big ( \max_{M \mly l \mly \ve V\ve^\rre}\frac
{\eta(z\iz{l,V}+y)}{\log l} < 1-\delta \Big ) \mly
\prod_{l=M}^{\infty} \big(1-\e^{-(1-\delta) \log l}\big) = 0
$$
as $M \to \infty$, i.e., assertion (iii) of Lemma 5.1 is
proved.\qed
\end{proof}

\medskip

{\it Proof of Lemma 5.2}. Fix a constant $\theta  \in (\rre,\frac
12)$,
 so that $K:=K\iz{V}\mly \frac 12|V|^\theta $ for each $V \supset
V_0$. Write $L\iz{V}:=(1-\theta )\log|V|$. Denote by
$\cE\iz{V}\subset V$ the subset consisting of sites at which
$\eta(\cdot)$ exceeds the level $L\iz{V}$, and let $r(\cE_V)$ be
the minimum distance among sites in $\cE_V$; cf.~the notation at
the beginning of Section 3. We abbreviate by $I$ the intervals
$(-\infty,t]$ or $(t,\infty)$, where $t \in \ttR$. Further, pick
$\delta$ to satisfy $0 < \delta <(1-2 \theta )/\nu$. Now
\begin{eqnarray*}
&&\bP\Big( \max_{1 \mly l \mly K}\eta(z\iz{l,V};l)  \in I \Big)\\
&&\mly  \bP\Big( \max_{1 \mly l \mly K}\eta (z\iz{l,V};l)\in I,\
2^{-1}|V|^\theta \mly |\cE\iz{V}|\mly 2|V|^\theta ,\ r(\cE\iz{V})>
|V|^\delta  \Big)\\
&&\quad +\bP\big( |\cE\iz{V}| < 2^{-1} |V|^\theta  \big) +\bP\big(
|\cE\iz{V}|
> 2|V|^\theta  \big)  +\bP\big( r(\cE\iz{V})\mly |V|^\delta  \big)\\
&&=\colon p(I)+p^{(1)}+p^{(2)} +p^{(3)}.
\end{eqnarray*}
Using the continuity of exponential distribution, similarly as in
the proof of Theorem~3.1(ii)-(iv) with $R=0$ and $0< \theta <1/2$,
we obtain that $p^{(i)} \mly |V|^{-{\rm const}}$ for some $\const
> 0$. Thus, to show the assertion of Lemma 5.2, we need to check
that \beq p(I)\mly \bP\Big( \max_{1 \mly l \mly K}\eta (\wt
x\iz{l};l) \in I \Big) +|V|^{-{\rm const}_1 }
\eeq for a fixed (nonrandom) subset $\wt V:=\{\wt x\iz{l}\colon 1
\mly l \mly K\} \subset \ttZ^\nu$ such that $r(\wt  V)>
|V|^\delta$.

Let $\sum\iz{V'}$ be the sum over all subsets $V' \subset V$ with
the properties: $\frac 12 |V|^\theta  \mly |V'| \mly 2|V|^\theta$
 and $r(V')> |V|^\delta$.
We denote by $\sum\iz{\{x_l\}}$ the sum over all permutations
$x_1,\ldots,x\iz{|V'|}$ of the subset $V'$. Write $p\iz{V}\colon
=|V|^{\theta-1}$. Then, for $V \supset V\iz{0}$, \beq
&&\hspace{-1cm}p(I)\mly \sum\limits_{V'}\sum\limits_{\{x_l\}}
\bP\Big(\max\limits_{1 \mly l \mly K} \eta
(x\iz{l};l) \in I, \nn  \\
&&\qquad\quad  \eta(x_1) \dly \eta (x_2) \dly \ldots \dly
\eta(x\iz{|V'|}) \dly L\iz{V},\max\limits_{x \in V\backslash V'}
\eta(x) <
L\iz{V}\Big)\nn \\
&&\hspace{-1cm}\ \ \mly \sum\limits_{V'}\sum\limits_{\{x_l\}}
\bP\Big( \max\limits_{1 \mly l \mly K}\eta (x\iz{l};l)\in I,
\eta(x_1) \dly \eta(x_2) \dly \ldots \dly \eta(x\iz{|V'|})
\dly  L\iz{V},\quad  \nn \\
&&\qquad\quad \max \{\eta(x)\colon x \in V\big\backslash (
(V'+U)\cup V')\} < L\iz{V}\Big),
\eeq where $V'+U$ denotes the algebraic sum of the subsets $V'$
and $U$. Since all the random variables are mutually independent,
the double sum on the right-hand side of (5.2) is equal to
\begin{eqnarray*}
&&\bP\Big( \max_{1 \mly l \mly K}\eta(\wt  x_l;l)\in I\Big)
\sum_{V'}p_{V}^{|V'|}(1-p\iz{V})^{|V| -(|U| +1)| V'|} \\
&&\quad {}\mly\bP\Big( \max_{1 \mly l \mly K}\eta(\wt  x_l;l) \in
I\Big)+\const |V|^{2\theta -1}\quad(V \supset V_0),
\end{eqnarray*}
since $|V'|p\iz{V} \asymp |V|^{2\theta -1}$ via the notation. This
completes the proof of (5.1). Lemma~5.2 is proved.\qed

\medskip

We now turn to the case $\xi(\cdotp)=f(\eta (\cdotp))$.

\begin{thm}{\rm (Cf.~(2.10)).} Fix $y\in \ttZ^\nu \backslash
\{0\}$ and $K \in \ttN$. Then
$$
\lsup_V \big|\xi (z\iz{K, V}+y)\big| < \infty \quad \mbox{in \
probability }.
$$
\end{thm}

\medskip

\begin{proof}
The assertion follows from Lemma 5.1(i).\qed
\end {proof}

To the end of this section, let us fix constants
$0<\rre<\theta<1/2$. With $L\iz{V,\theta}$ as in Section 3, we
write $\wt \xi(x):=\xi(x)$ if $\xi(x) < L\iz{V,\theta }$, and $\wt
\xi(x):=0$, otherwise. For natural $K \dly 1$ and $l >K$, we put
\setcounter{equation}{2}

\begin{equation}
\begin{split}
\chi\iz{K,V}(l):=&\xi_{K,V}^2\bigg(
\xi_{K,V}+2\nu\rrk^2 \xi_{K,V}^{-1}-\xi_{l,V}\\
&-\rrk^2\mathop{\sum}\limits_{|x|=1}\big(
\xi_{K,V}-\wt\xi(z_{l,V}+x)
 \big)^{-1} \bigg)\vien_ {\big \{
\xi\iz{K,V} > L\iz{V,\theta }\big \} };
\end{split}
\end{equation}
here $\vien_{ \rrO^{\prime} }:= \vien_{ \rrO^{\prime} }(\rro)$
denotes the indicator of $\rrO^{\prime} \subset \rrO$. To study
the asymptotic behavior of variables (5.3), we introduce the class
$S A \Pi_\infty ^2$ of functions $f:=Q^\la$ such that
\begin{equation}
\begin{split}
&\lim_{s \to \infty} \inf_{a \in (c,\theta s)}\bigg(
f(s)^2\big(f(s+a)-f(s)\big)-\frac{f(2a)}{c}
\bigg) =\infty \\
&\mbox{for\ any}\ 0 < c < 1\quad  \mbox{and\ some}\ 0 <\theta
<1/2.
\end{split}
\end{equation}
The class $SA\Pi_\infty ^2$ is a strict subset of $A\Pi_\infty ^2$
(4.1). The following theorem provides some generalization of limit
(4.4) for $p=2$.

\begin{thm}{\rm (Cf.~(2.11)).} Fix $K \in \ttN$.
If $f$ belongs to the classes  $SA\Pi^{2}_{\infty}$ {\rm (5.4)}
and $PI_{<2}$ {\rm (4.5)} with $\rre$ and $\theta$ as above, then
$$ \lim_V \min_{K+1 \mly l \mly 2\vert
V\vert^\theta} \chi\iz{K,V}(l)=\infty\ \ \mbox{in\ probability}.
$$
\end{thm}

\medskip

\begin{proof}
We begin with estimating $\chi(l):=\chi\iz{K,V}(l)$ for $K+1 \mly
l \mly M$; $M \dly M_0$. Let $\rrO_{V,M}^{(1)} \in \cF$ denote the
subset of configurations $\xi^{(\rro)}_V$ satisfying the following
three inequalities: \beq \max_{|x|=1}\max_{K \mly l \mly
M}\xi(z\iz{l,V}+x) \mly f(2 \log M),
\eeq \beq \xi\iz{K+1,V}>0 \quad \mbox{and} \quad
\frac{\wt\xi(\cdotp)}{\xi\iz{K,V}} \mly
\frac{L\iz{V,\theta}}{\xi\iz{K,V}} \mly \const\pr<1
\eeq for $\const\pr>\const(\rre)$ specified in Theorem 4.5.
According to Lemma 5.1(ii) and Theorem 4.5, we obtain that
$\lsup_V\bP(\rrO\backslash \rrO_{V,M}^{(1)}) \to 0$ as $M \to
\infty$. On the other hand, expanding the sum $\sum_{|x|=1}$ in
(5.3) over powers of $\wt\xi(z\iz{l,V}+x)/\xi\iz{K,V}$ with $K+1
\mly l \mly M$, we get that, for any $M \dly M_0$ and any $V
\supset V_0(M)$, the inequalities (5.5) and~(5.6) imply the
following estimate
$$
\min_{K+1\mly l \mly M}\chi(l) \dly \big(\xi_{K,V}\big)^2\big(
\xi_{K,V}-\xi_{K+1,V} \big)-\const f(2 \log M),
$$
where $\const >0$ does not depend on $V$ and $M$. From this
implication and Theorem~4.3(i) with $p=2$, we obtain that, for any
$C
> 0$, \beq \lsup_{V}\bP\big( \min_{K+1 \mly l \mly M}\chi(l)\mly C
\big) \mly \lsup_V \bP\big( \rrO\backslash \rrO_{V,M}^{(1)}
\big)\to 0
\eeq as $M \to \infty$.

It only remains to estimate $\chi(l):=\chi\iz{K,V}(l)$ for $M \mly
l \mly 2\vert V\vert^\theta$. Using formulas~(1.23) and (1.24), we
represent $\xi\iz{l,V}$ and $\xi(z\iz{l,V}+x)$ in the form:
$$
\xi\iz{l,V}=f\bigg( \log\frac{|V|}{l}+\rho\iz{V}(l) \bigg)\quad
\hbox{and}\quad \xi(z\iz{l,V}+x)=f(\eta(z\iz{l,V}+x)),
$$
where $\rho\iz{V}(l):=\eta\iz{l,V}-\log (|V| /l)$. Denote by
$\rrO_{V,M}^{(2)} \in \cF$ the subset of configurations
$\xi^{(\rro)}(\cdot)=f(\eta^{(\rro)}(\cdot))$ satisfying (5.6) and
the following three inequalities: \beq |\rho\iz{V}(K)| \mly \frac
13\log M,\quad \max_{M \mly l \mly 2\vert
V\vert^\theta}|\rho\iz{V}(l)| < 1
\eeq and \beq \max_{|x|=1}\max_{M \mly l \mly 2\vert
V\vert^\theta} \frac{\eta(z\iz{l,V}+x)} {\log l} < \frac 32.
\eeq We then apply Lemma 4.1(ii) to  $\rho\iz{V}(l)$ and
Lemma~5.1(iii) to the left-hand side of~(5.9) to obtain that
$\lsup_V\bP(\rrO\backslash \rrO_{V,M}^{(2)}) \to 0$ as $M \to
\infty$. We now write $\xi_+:=\xi\vee 0$ and note that, for any $M
\dly M_0$ and any $V \supset V_0(M)$, inequalities (5.6), (5.8)
and (5.9) imply the following estimate:
\beq\!\!\!&&\!\!\!\min\limits_{M \mly l \mly 2\vert V\vert^\theta}
\chi(l) \nn \\
&& \dly \min\limits_{M \mly l \mly 2\vert V\vert^\theta}\bigg[
\xi_{K,V}^2\big( \xi_{K,V}-\xi_{l,V} \big)\!
-\!\const\!\sum\limits_{|x|=1}
\xi_+(z\iz{l,V}\!+\!x) \bigg]\nn \\
&&{}\dly \min_{M\mly l\mly 2\vert V\vert^\theta}\bigg[ f\bigg(
\log\ds\frac{|V|}{\sqrt{M}}
 \bigg)^2\bigg( f\bigg( \log\ds\frac{|V|}{\sqrt{M}} \bigg)-f\bigg(
\log\ds\frac{|V|}{l}+1 \bigg) \bigg) \nn \\[12pt]
&&\qquad\qquad\qquad {}-{\rm const'}f \bigg(\ds\frac 32 \log l
\bigg) \bigg].
\eeq Using this implication combined with the fact that, by
condition (5.4), the right-hand side of (5.10) tends to infinity
as $|V| \to \infty$,
 we obtain that, for any $C >0$,
$$
\lsup_{V}\bP \lf( \min_{M \mly l \mly 2\vert
V\vert^\theta}\chi(l)\mly C \rg)\mly \lsup_{V} \bP\big(
\rrO\backslash \rrO_{V,M}^{(2)} \big) \to 0
$$
as $M \to \infty$. This limit and (5.7) yield the assertion of
Theorem 5.4.\qed
\end{proof}

\smallskip

\setcounter{equation}{0} \addtocounter{sct}{1}

\section{\normalsize{\textsf {POISSON LIMIT THEOREMS FOR THE LARGEST EIGENVALUES} }}

\setcounter{prop}{1}

This section is to provide an overview of current results on
the extreme value theory for the spectrum of the Anderson Hamiltonian
$\cH_V=\rrk \Delta_V +\xi_V$, $V\uparrow \ttZ^\nu$, with an i.i.d.~
potential $\xi(\cdot)$. The results under consideration are taken
from (Astrauskas and Molchanov~1992), (Astrauskas 2007; 2008;
2012; 2013), (Bishop and Wehr~2012), (Germinet and Klopp 2013),
(G\" artner and Molchanov 1998) and (Biskup and K\"onig~2016).

In Section~6.1, we give Poisson limit theorems for the largest
eigenvalues and the corresponding localization centers, provided
the distribution tails $\e^{-Q}$ of $\xi(0)$ are heavier than the
double exponential function (Theorems 6.2 and 6.9). These limit
theorems are then complemented and illustrated by the
distributions with polynomially decaying tails, Weibull
distributions and those with fractional double exponential tails
(resp., Examples 6.11, 6.12 and 6.13).

In Section 6.2, we first give the second order expansion formulas
for the largest eigenvalues, provided the tails $\e^{-Q}$ are
lighter than the double exponential function (Theorem~6.14). For
bounded $\xi(0)$, further extensions of this result are discussed.

Section 6.3 provides the second order expansion formulas for the
principal eigenvalue in the case of double exponential tails
(Theorem~6.19), which are further extended up to Poisson limit
theorems for eigenvalues.

In Section 6.4, we comment and compare the proofs of Poisson limit
theorems stated in Sections~6.1 and~6.3. We mention, en passant,
that Theorems 6.2, 6.9, 6.14 and 6.19 simply follow from the
corresponding results of Sections 2--5.

\subsection{\normalsize{\textsf {Distribution tails heavier than the double exponential function } }}

The extreme value theory for i.i.d.~random variables $\xi(x)$ deals
with the asymptotic behavior of the $K$th largest values $\xi
\iz{K,V}$ of the sample $\xi_V$ as $V\uparrow\ttZ^\nu$. It is well
known that for suitable normalizing constants $a_V>0$ and $b_V$,
the non-trivial limiting (max-stable) distributions $G(\cdotp)$
for $\bP \big((\xi \iz{1,V}-b\iz{V})a\iz{V} \mly \cdotp \big)$ are
either Weibull law $D_\beta(t):=\exp \lf \{ -(-t)^{\beta} \rg \}$
($t<0$) or Fr\'{e}chet law $G_\beta (t)= \exp\{-t^{-\beta}\}$
($t\dly 0$) for some $\beta
>0$, or Gumbel law $G_{\exp}(t)=\exp\{-\e^{-t}\}$ ($-\infty<t<\infty$);
see, e.g., (Resnick 1987). Note that the weak convergence of
maxima to Gumbel law is equivalent to the limit \beq
\lim\limits_V|V| \bP \big( \xi(0)
>b_V+t/a_V \big)=\e^{-t} \quad \mbox{for\ all} \quad t \in \ttR.
\eeq On the other hand, limit (6.1) implies that the point process
$\cN_V^\xi$ on $[-1/2;1/2]^\nu \times \ttR$, defined by \beq
\cN_V^\xi :=\sum_{z \in V} \delta_{\Xi \iz{V}(z)} \quad \hbox
{where} \quad \Xi\iz{V}(z):=\big(z|V|^{-1/\nu},(
\xi(z)-b\iz{V})a\iz{V} \big),
\eeq converges weakly (as $V\uparrow \ttZ^\nu$) to the Poisson
process on $[-1/2;1/2]^\nu \times \ttR$ with the intensity measure
$\dr x\times \e^{-t} \dr t$, i.e., the product of Lebesgue 
measure on $[-1/2;1/2]^\nu$ and that defined by the increasing 
function $\log G_{\exp}(\cdot)$ on $\ttR$; see (Leadbetter et al.~1983).

The necessary and sufficient conditions for (6.1) to hold are
generally formulated in terms of $\Gamma$-variation of the
function $\e^Q$ at the right endpoint $t_Q$  or, equivalently, in
terms of $\Pi$-variation of its inverse $f\!\circ\!\log$ (Resnick
1987). We say that $f:=Q^\la$ is in the class $A\Pi$, if there
exists a function $a:(-\infty; t_Q) \to \ttR_+$ such that \beq
\lim_{s\to \infty}\frac{f(s+c)-f(s)}{a(f(s))}= c  \quad \hbox{for
any}\quad c\in \ttR_+.
\eeq Here $a(\cdot)$ is called an auxiliary function. The class
$A\Pi$ (6.3) is an argument-additive version of the original class
of $\Pi$-varying functions $f\circ \log$ considered, e.g., in
(Resnick 1987, Section 0.4.3). In Lemma A.1 of Appendix A, we
recall the well-known characterization of the class $A\Pi$ in
terms of $Q$.

\begin{lem}{\rm (Resnick 1987, Sections 0.4.3 and
1.1).}  Limit {\rm (6.1)} holds true if and only if $f \in A\Pi$
for some an auxiliary function $a(\cdot)$. In this case, the
normalizing constants can be chosen $b\iz{V}=f(\log\ve V\ve)$ and
$a\iz{V}=1/a(b\iz{V})$.
\end{lem}

\medskip

We now formulate Poisson limit theorems for the largest
eigenvalues $\rrl\iz{K,V}$ of the random Schr\"odinger operator
$\cH_V=\rrk \Delta_V +\xi_V$ introduced in Section 2.1. Throughout
this subsection, we assume that $\e^f \in RV_{\infty}$, so that
$\log Q(t)=\mro(t)$ as $t \rightarrow \infty$ by Lemma A.3 with
$\rho=\infty$. This class of distributions includes Weibull
distributions (1.5) for arbitrary $\rra>0$ and those with
fractional double exponential tails (1.8) with $\rrg<1$.
Using the notation from Section 4.3, for fixed
small $0 < \theta < 1/2$, we write $\wt \xi(x):=\xi(x)$ if $\xi(x)
< L\iz{V, \theta}:=f((1-\theta) \log \ve V \ve)$, and $\wt
\xi(x):=0$, otherwise. For any $z \in V$, denote by $\wt \rrl(z)$
the principal eigenvalue of the \lk single peak\rk \ Hamiltonian
$\rrk\Delta_{V} +\xi(z)\delta \iz{z} +\wt\xi_V (1-\delta \iz{z})$.
Let $\wt \rrl_{K,V}$ be the $K$th order statistics of the
stationary random field $\wt \rrl(\cdot)$ in $V$, and let
$z_{\tau(K),V} \in V$ stand for its location defined by $\wt
\rrl(z_{\tau(K),V}):= \wt \rrl_{K,V}$. (Recall that the sites
$z\iz{l,V} \in V$ $(1 \mly l \mly |V|)$ are associated with the
variational series~(1.1) based on $\xi_V$.) Note that, for
$Z:=z_{\tau(K), V}$ and $K \in \ttN$ fixed, the eigenvalues $\wdt
\rrl(Z)$ are expanded into a certain power series in the variables
$\xi(Z)$ and $\wt \xi(x)$ ($x \in V$); cf.~(2.20)--(2.22).

\begin{thm} {\rm (see Theorem 4 in (Astrauskas 2007) and
Theorem 5.2 in (Astrauskas 2008)).} Let $t_Q =\infty$ and $\e^f \in RV_{\infty}$,
and assume
that $\e^{-Q}$ is log-H\" older con\-ti\-nu\-ous of order
$\mu\!>\!(1\!+\!\theta )\nu /(1\!-\!2\theta )$ at infinity for
some small $\theta >0$ as above, i.e., {\rm (4.7)} holds true. Then the following assertions {\rm (I)--(II)} hold true:

{\rm (I) (Poisson limit theorem)} Assume, additionally, that there
exist the normalizing constants $A\iz{V}
> 0$ and $B\iz{V}$ such that \beq
\lim_V |V| \bP(\wt \rrl(0) > B\iz{V} +A_V^{-1}t)=\e^{-t} \quad
\mbox{for\ any} \quad t \in \ttR,
\eeq and define the point process $\cN_V^\rrl$ on $[-1/2;1/2]^\nu
\times \ttR$ by \beq \cN_V^\rrl :=\sum_{k=1}^{\ve V\ve}
\delta_{\Lambda \iz{V}(k)} \quad \hbox {where} \quad \Lambda
\iz{V}(k):=\bigg(\frac {z\iz{\tau (k),V}}{|V|^{1/\nu}},\big(
\rrl\iz{k,V}-B\iz{V}\big)A\iz{V} \bigg);
\eeq then $\cN_V^\rrl$ converges weakly to the Poisson process
$\cN$ on $[-1/2;1/2]^\nu \times \ttR$ with the intensity measure
$\dr x\times \e^{-t} \dr t$.

{\rm (II) (Exponential localization)} Fix a constant $\rre$ such that $0<\rre<\theta$ and write $M_V:=\log( L\iz{V, \rre}-L\iz{V,
\theta})$ (so that $M_V\to \infty$); then with probability one
\beq \lsup_V\max_{1 \mly K \mly \ve V\ve^\rre}\max_{x \neq
z_{\tau(K),V}}\frac {\log\big|\psi(x;\rrl\iz{K,V})\big|}
{M\iz{V}|x-z\iz{\tau(K),V}|} \mly -1.
\eeq
\end{thm}

\medskip

{\it Sketch of the proof of Theorem~6.2(I)}. Using  Theorem 4.6 with
$\rho=\infty$, Theorems~3.1 and 4.8 with $R=0$ and $0 <\rre
<\theta <1/2$, we obtain that almost surely $\xi_V$ satisfies the
assumptions of Theorem 2.3, where $K \in \ttN$ is fixed and
$N:=[\ve V\ve^\theta]$. Theorem~2.3 implies that almost surely
$\rrl_{K,V}=\wt \rrl_{K,V}+ \rO\big (\exp \{-\ve
V\ve^{(1+\theta)/\mu}\} \big)$ as $\ve V\ve \to \infty$, for fixed
$K \in \ttN$. This asymptotic formula in turn yields that the
point process $\cN_V^\rrl$ is approximated by the corresponding
point process $\wt \cN_V^{\rrl}$ where $\rrl\iz{k,V}$ are replaced
by $\wt \rrl _{k,V}$ ($1 \mly k \mly \ve V \ve$); see the proof of
Theorem 4 in (Astrauskas 2007). The weak convergence of $\wt
\cN_V^{\rrl}$ to $\cN$ is shown by checking Leadbetter's mixing
conditions for the random field $\wt \rrl(\cdot)$ (Astrauskas
2007, Lemma 6). This concludes the proof of the theorem.\qed

\medskip

In Corollaries 6.3--6.7 below, we give
the alternative conditions on $Q$ (where $\log Q(t)=\mro(t)$) for
the Poisson convergence of the largest eigenvalues to hold.

\begin{cor}
({\it Specification of the normalizing constants $A\iz{V}> 0$
and $B\iz{V}$ in~{\rm (6.4)} for some examples of potential
distributions}). Let $Q(t)=t^{\rra}$  for
$t\dly 0$ where $\rra >0$ (Weibull distribution), or
$Q(t)=\e^{t^{\rrg}}$ for $t\dly t_0$; $ 0<\rrg <1$ (fractional
double exponential distribution). Consequently, $Q$ satisfies the
regularity and continuity conditions of Theorem 6.2. Moreover, the
equations for the normalizing constants $A_V
> 0$ and $B_V$ in (6.4) are derived by applying a certain
iteration scheme for $\wt \rrl_{1,V}$ as in (2.20)--(2.22)
combined with Laplace's method for the corresponding integrals
(Astrauskas 2008, Section 6), (Astrauskas~2016); see also
(Astrauskas 2013, Section~3). From these equations one derives the
explicit expansion formulas for $B\iz{V}$ and hence those for the
eigenvalues $\rrl\iz{K,V}$ up to the random max-stable
fluctuations of order $\rO(A^{-1}_V)$; cf.~Examples 6.12 and 6.13
below.
\end{cor}

\medskip

\begin{cor}
({\it Specification of the normalizing constants $A\iz{V}> 0$
and $B\iz{V}$ in~{\rm (6.4)} under general RV conditions on potential
distributions}). Let $t_Q= \infty$ and, for some large $t_0$, assume that
$Q:[t_0; \infty) \rightarrow \ttR_+$ is (locally) absolutely
continuous with the positive density $Q^{\prime}:[t_0; \infty)
\rightarrow \ttR_+$ obeying the following conditions: \beq
\lim_{t\to \infty }\frac {Q^{\prime}(t+C)} {Q^{\prime}(t)}= 1\quad
\mbox{for\ any} \quad C>0,
\eeq and \beq \linf_{t\to \infty }Q^{\prime}(t)>0.
\eeq Consequently, the
function $Q$ satisfies the regularity and continuity conditions of
Theorem 6.2. Then the
centralizing constants $B_V$ in (6.4) are defined by the
equation: \beq \bP(\wt \rrl(0) > B\iz{V})=\ve V\ve^{-1} \quad \
(V\supset V_0),
\eeq and the normalizing constants $A_V >0$ in (6.4) and $a_V>0,
b_V$ in (6.1) are specified as follows: \beq
A_V=a_V:=Q^{\prime}(b_V) \quad \mbox{where} \quad b_V:= f(\log \ve
V\ve) \quad \ (V\supset V_0).
\eeq Therefore, for any $t \in \ttR$ and $\ve V\ve \to \infty$,
\beq |V| \bP \Big( \xi(0)
>b_V+\frac{t}{a_V}\Big) \to \e^{-t} \quad \mbox{and} \quad  |V|
\bP \Big(\wt \rrl(0) > B\iz{V} +\frac{t}{a_V} \Big) \to \e^{-t},
\eeq where \beq B_V=b_V+\mro(1).
\eeq
\end{cor}

\medskip

Corollary 6.4 ensures that both the distribution functions $\bP(\xi(0) \mly
t)=1-\e^{-Q(t)}$ and $\bP(\wt \rrl(0) \mly t)$ ($t\in \ttR$) are
in the domain of attraction of the max-stable Gumbel law
$G_{\exp}(\cdot)$; cf.~Lemma 6.1 and the assertions before this
lemma. Note also that the additional condition (6.8) is to exclude the
heavy-tailed (\lk subexponential\rk) distributions $1-\e^{-Q}$
which are considered in Theorem 6.9(C0) and Example 6.11 below.

{\it The proof of Corollary {\rm 6.4}}. We observe from Lemmas
A.4(I) and A.3 with $\rrr=\infty$ that conditions (6.7) and (6.8)
imply the regularity and continuity conditions of Theorem 6.2.
Further, from Lemma A.4(III) and Lemma A.1 we derive that $f \in
A\Pi$~(6.3) with the auxiliary function
$a(\cdot)\equiv1/Q^{\prime}(\cdot)$ in $[t_0; \infty)$; therefore,
by Lemma 6.1 we obtain the first limit in (6.11).

To prove the second limit in (6.11), we first notice that, for
each $V\supset V_0$, there is a solution $B_V$ of equation (6.9)
because of the continuity of the distribution function of $\wt
\rrl(0)$. Let us show (6.12). Since $\wt \rrl(0)\dly \xi(0)$, we
get from (6.9) that
$$
\ve V\ve^{-1}=\bP(\wt \rrl(0) > B\iz{V})\dly \bP(\xi(0)
> B\iz{V})=\e^{-Q(B_V)},
$$
therefore, $B_V \dly b_V=f(\log \ve V\ve)$ for $V\supset V_0$. If
$\xi(0) \dly L\iz{V,\rre}:=f((1-\rre)\log\ve V\ve)$ for some $0<
\rre <\theta$, then we get from (2.20)--(2.22) that almost surely
$\wt \rrl(0)\mly \xi(0)+\beta_V$ for some (nonrandom) $0
<\beta_V\downarrow 0$ as $\ve V\ve \to \infty$. Thus, for
$V\supset V_0$,
\begin{eqnarray*}
&&\ve V\ve^{-1}=\bP(\wt \rrl(0) > B\iz{V})=\bP(\wt \rrl(0) >
B\iz{V}, \xi(0) \dly L\iz{V, \rre}) \\[6pt]
&&\hspace{.7cm} \mly \bP(\xi(0)
> B\iz{V}-\beta_V)=\e^{-Q(B_V-\beta_V)},
\end{eqnarray*}
therefore, $B_V\mly b_V+\beta_V=b_V +\mro(1)$. These estimates
imply (6.12), as claimed. To prove the second limit in (6.11), we
also need the following observations. First,  since $\linf_V a_V
>0$, it follows that $\lsup_V |V|
\bP (\xi(0) \dly b\iz{V} +M) \to 0$ as $M \to \infty$. Second, for
any $M \dly M_0$ and any $V \supset V_0(M)$, if $\xi(0) \mly
b\iz{V} -M$, then $\wt \rrl(0)< b\iz{V} -M/2$. These two
assertions imply that, for any $t \in \ttR$,
\begin{eqnarray}
&&|V| \bP \big(\wt \rrl(0) > B\iz{V} +ta_V^{-1} \big) \n \\[6pt]
&&\hspace{.7cm} = |V| \bP \big(\wt \rrl(0) > B\iz{V} +ta_V^{-1} ,
| \xi(0)-b_V | <M \big) +\mro\iz{V,M}(1),
\end{eqnarray}
where $\mro\iz{V,M}(1) \to 0$ letting first $V\uparrow \ttZ^\nu$
and then $M \to \infty$. Thus, it suffices to check that, for any
$t \in \ttR$, any $M\dly M_0(t)$ and $V\uparrow \ttZ^\nu$,
\begin{eqnarray}
&&\bP \big(\wt \rrl(0) > B\iz{V} +ta_V^{-1},
| \xi(0)-b_V | <M  \big) \n \\[6pt]
&&\hspace{.7cm} = \e^{-t}\bP \big(\wt \rrl(0) > B\iz{V}, |
\xi(0)-b_V | <M \big)(1+\mro(1)).
\end{eqnarray}
To prove (6.14), we follow the arguments of the paper (Biskup and
K\"onig 2016, Section~7.1) which are now simplified and adapted to
our case $\log Q(t)=\mro(t)$. (Recall that this paper considers
the case of double exponential tails of potential, i.e., $\log
Q(t)\approx t/\rrr$.) The main idea here is the observation that
the shift of the eigenvalue $\wdt \rrl(0)$ by $ta_V^{-1}$ is
achieved by the corresponding shift of the single $\xi_V$-peak
$\xi(0)$ on the left-hand side of (6.14). Indeed, write
$\xi^{(t)}:= \xi(0)-ta_V^{-1}$, and denote by $\rrl^{(t)}$ the
principal eigenvalue of the Hamiltonian $\rrk\Delta_{V}
+\xi^{(t)}\delta \iz{0} +\wt\xi_V (1-\delta \iz{0})$ in $l^2(V)$.
Note that $\xi^{(0)}= \xi(0)$ and $\rrl^{(0)}=\wdt \rrl(0)$. Fix
$t \dly 0$. Comparing expansion formulas (2.20)--(2.22) for
$\rrl^{(t)}$ with those for $\wdt \rrl(0)$, we find that, for any
(small) $\rre >0$, any $M \dly M_0(t,\rre)$ and any $V\supset
V_0(M,t,\rre)$,
$$
\mbox{if} \quad |\xi(0)-b_V| < M, \quad \mbox{then} \quad
\rrl^{(t+\rre)} \mly \wdt \rrl(0)-ta_V^{-1} \mly \rrl^{(t)};
$$
therefore, we obtain the following bounds for the left-hand side
of (6.14):
\begin{eqnarray}
&&\bP \big(\rrl^{(t+\rre)} > B\iz{V}, | \xi^{(t+\rre)}-b_V |
<M/2 \big) \n \\[6pt]
&&\hspace{.7cm}\mly \bP \big(\wt \rrl(0) > B\iz{V} +ta_V^{-1},
| \xi(0)-b_V | <M  \big) \\[6pt]
&&\hspace{.7cm} \mly \bP \big(\rrl^{(t)} > B\iz{V}, |
\xi^{(t)}-b_V | <2M \big) \n.
\end{eqnarray}
Since $\rre > 0$ is arbitrarily small, it suffices to prove limit
(6.14) for the upper and lower bounds in (6.15). We write
$\rrl^{(t)}(\xi_V)=\rrl^{(t)}$ to emphasize the dependence of
$\rrl^{(t)}$ on the sample $\xi_V=\{\xi(x)\}_{x \in V}$. Let us
consider the functions $\rrl^{(t)}(s_V)$ of $s_V=\{s(x)\}_{x \in
V} \in \ttR^{\ve V\ve}$ and the corresponding integrals on the
right hand-side of (6.15) with respect to the probability measure $\prod_{x
\in V}\dd F(s(x))$; here $F:=1-\e^{-Q}$ stands for the
distribution function of $\xi(0)$ with the density
$p(\cdot):=F^{\prime}(\cdot)$ in $[t_0; \infty)$. By the change of
variables $u(0):=s(0)-ta_V^{-1}$ and $u(x):=s(x)$ for all $x \in
V\setminus \{0\}$,  we get that
$\rrl^{(t)}(s_V)=\rrl^{(0)}(u_V)$; therefore, for fixed $M>0$ and
$V\uparrow \ttZ^{\nu}$,
\begin{eqnarray}
&&\bP \big(\rrl^{(t)} > B\iz{V}, |
\xi^{(t)}-b_V | <M \big) \n \\[6pt]
&&\hspace{.3cm}= \int_{\ttR^{\ve V\ve}} \vien_{\big
\{\rrl^{(t)}(s_V)
> B_V, | s(0)-ta_V^{-1}-b_V | <M\big\}}\prod_{x
\in V}\dd F(s(x)) \n \\[6pt]
&&\hspace{.3cm} = \int_{\ttR^{\ve V\ve}} \vien_{\big
\{\rrl^{(0)}(u_V)
> B_V, | u(0)-b_V | <M\big\}}\frac {p( u(0)+ta_V^{-1})}
{p( u(0))}\prod_{x
\in V}\dd F(u(x)) \n \\[6pt]
&&\hspace{.3cm} = \ttE \bigg( \vien_{\big \{\rrl^{(0)}(\xi_V)
> B_V, | \xi(0)-b_V | <M\big\}}\frac {p( \xi(0)+ta_V^{-1})}
{p( \xi(0))}\bigg ) \\[6pt]
&&\hspace{.3cm} = \bP \big( \rrl^{(0)}(\xi_V)
> B_V, | \xi(0)-b_V | <M\big)\big(\e^{-t}+\mro(1)\big) \n
\end{eqnarray}
by applying Lemma A.4(IV) to the ratio of the densities in the
last expectation~(6.16) where $a_V=Q^{\prime}(b_V)$. Formula (6.16) combined
with~(6.15) implies~(6.14) for $t \dly 0$, as claimed. Since the
case $t \mly 0$ is treated similarly, this concludes the proof of
assertions of Corollary~6.4. \qed

\medskip

\begin{cor}
({\it Suppression of the $\log$-H\"older continuity of $\e^{-Q}$
in Theorem {\rm 6.2}}). Let $t_Q= \infty$. Assume that
 $\e^f \in RV_{\infty}$, and let
assumption (6.4) be fulfilled with \beq A_V=\rO\big (
(L\iz{V,\rre} -L\iz{V,\rre^{\prime}})^{\ve
V\ve^{(1-2\delta)/\nu}}\big ) \quad \mbox{for\ some\ constants}
\quad 0< \rre < \rre^{\prime}< \delta < \frac 12.
\eeq Then the point process $\cN_V^\rrl$ (6.5)
converges weakly to the Poisson process
$\cN$ as in Theorem 6.2(I) and, moreover, the $K$th eigenfunction obeys
exponential localization~(6.6) in probability.
\end{cor}

\medskip

To {\it prove Corollary 6.5}, we again apply
Theorem 2.3. So we need to show that the samples $\xi_V$ satisfy
limits (2.9), (2.14) and (2.15) in probability with the same
abbreviation as in the proof of Theorem 6.2. First, since the
condition $\e^f \in RV_{\infty}$ implies (3.1) (see Lemma A.3 with $\rrr=\infty$), we may apply Theorem 3.1(iv) with $R=0$ to obtain limit
(2.9) with $N=[\ve V\ve^\theta]$ where $\rre^{\prime} < \theta <
\delta$. Second, by Theorem~4.6 with $\rrr = \infty$, the
condition $\e^f \in RV_{\infty}$ yields (2.14) with fixed $K \in
\ttN$ and $N$ as above. It remains to prove limit (2.15) in
probability with those $K$ and $N$. As mentioned in the proof of
Theorem~6.2, assumption~(6.4) implies that the point process $\wt
\cN_V^{\rrl}$, based on the sample $\wdt \rrl_V$, converges weakly
to the corresponding Poisson process. This convergence in turn yields that
with probability $1+\mro(1)$ the normalized spacings
$A\iz{V}\big(\wdt \rrl\iz{k, V}-\wdt \rrl\iz{k+1, V}\big)$ are
bounded away from zero as $V\uparrow \ttZ^\nu$, for any fixed $k
\in \ttN$ (Astrauskas 2007, Corollary 1(jj)). Combining this with
the upper bound (6.17) for $A_V$ and observing from Lemma~4.2(iii)
and Theorem~3.1(iv) ($R=0)$ that almost surely $\xi \iz{K,V}-\xi
\iz{N,V} \dly L\iz{V,\rre}-L\iz{V,\rre^{\prime}}$ and $r\iz{N,V}
\dly \ve V\ve^{(1-2\delta)/\nu}$ for any $V\supset V_0$, we arrive
at limit (2.15) in probability with $N=[\ve V\ve^\theta]$, as
claimed. The assertions of Corollary~6.5 are proved. \qed

\medskip

\begin{cor}
({\it The second order expansion formula for the top eigenvalues}). Let again $t_Q= \infty$, and  $\e^f \in RV_{\infty}$.  Assume, in addition, that $f(s+\log s)-f(s) \to 0$ as $s \to \infty$. Then, for fixed $K\dly 1$, with probability one
$$
\lim_V\big(\rrl\iz{K,V}-f(\log \ve V\ve)\big)=\lim_V\big(\xi\iz{K,V}-f(\log \ve V\ve)\big)=0.
$$
\end{cor}

\medskip

Notice that the additional condition of Corollary 6.6 is to exclude
the heavy-tailed (\lk subexponential\rk) distributions
$1-\e^{-Q}$, or in other words, the class of i.i.d.~potentials
whose extremes possess sharp random fluctuations; see Lemma
4.4(iii).

{\it The proof of Corollary 6.6}. We first obtain from Theorem 3.1(iv) with $R=0$
and Theorem 4.6 with $\rrr=\infty$ that the samples $\xi_V$
($V\uparrow \ttZ^\nu$) satisfy almost surely the conditions of Theorem 2.1(ii) and
Remark 2.4 with $N:=[\ve V\ve^\theta]$ for some $0<\theta <1/2$. Thus, using the lower bound
for $\rrl\iz{K,V}$ (Theorem 2.1(ii)) and the almost sure limit (2.23)
for $\rrl\iz{1,V}$ (Remark 2.4) combined with Theorem 4.4 (iii), we arrive at the assertion of Corollary 6.6. \qed

\medskip

\begin{cor}
({\it Localization centers}). Assume that $Q$ satisfies the
conditions of Theorem 6.2 on the regular increase (i.e. $\e^f \in RV_{\infty}$) and the
$\log$-H\"older continuity at infinity. Then almost surely the
eigenfunction $\psi(\cdotp; \rrl \iz{K,V})$ is asymptotically
delta-function at the site $z_{\tau(K),V} \in V$ for each $K =
\mro (\ve V\ve^{\rre})$; see Theorem 6.2(II). Consequently, any
site $z_{\tau(k),V}$ in (6.5) can alternatively be defined as a
localization center of the eigenfunction $\psi(\cdotp; \rrl
\iz{k,V})$, viz. \beq \psi(z\iz {\tau (k),V}; \rrl \iz{k,V}):=\max
_{1\mly l\mly \ve V \ve} \psi(z\iz {l,V}; \rrl \iz{k,V})\quad
\mbox{for\ some}\ \ \tau(k)=\tau_V(k),
\eeq for all $1 \mly k \mly \ve V \ve$.
\end{cor}
\qed

The latter definition of the sites $z_{\tau(k),V}$ in (6.5)
is more natural in the context of the localization theory for the
Anderson Hamiltonians.

\medskip

The asymptotic behavior of the localization indices
$\tau(K)=\tau_V(K)$ is studied by Astrauskas (2013).

\begin{lem}{\rm (Astrauskas 2013, Theorem 2.1).} Assume that
the condition of Theorem {\rm 6.2} on the log-H\" older
con\-ti\-nui\-ty of $\e^{-Q}$ at infinity holds true. Fix $K \in
\ttN$.

{\rm (i)} If $f \in OA\Pi^2$ { \rm (4.3)}, then
$$
\lsup_V \tau\iz{V}(K) < \infty \quad \mbox{ in\ probability}.
$$

{\rm (ii)} If $f \in A\Pi_0^2$ {\rm (4.2)} and $\e^f \in
RV_{\infty}$, then
$$
\lim_V \tau\iz{V}(K)=\infty \quad \mbox{and}\ \lim_V \frac{\log
\tau\iz{V}(K)} {\log \ve V\ve}=0 \quad \mbox{ in\ probability}.
$$
\end{lem}

\medskip

Recall that, in Theorem 6.2 and Lemma 6.8, the condition $\e^f \in
RV_{\infty}$ implies $\log Q(t)=\mro(t)$ (see Lemma A.3 with
$\rho=\infty$ ); the condition $f \in OA\Pi^2$ yields $Q(t)\asymp
t^3$ (see Lemma A.8(ii) with $p=2$); finally, the condition $f \in
A\Pi_0^2$ implies $t^{-3} Q(t) \to \infty$ as $t \to \infty$ (see
Lemma A.7 with $p=2$).

In the case $Q(t)=\mro(t^3)$ as $t \to \infty$ (for example,
Weibull distribution (1.5) with $\rra<3$), the eigenvalue
$\rrl_{K,V}$ approaches the $K$th extreme value of $\xi_V$ as $V
\uparrow \ttZ^\nu$, for fixed $K\dly 1$. For such distributions, we
obtain a simplified version of Poisson limit theorems for the
largest eigenvalues:

\begin{thm} {\rm (see Theorem 5 in (Astrauskas 2007) and
Theorem 2.5 in (Astrauskas 2012)).} Let $t_Q=\infty$, and $f \in A\Pi$ {\rm (6.3)}
for some auxiliary function $a(\cdot)$, and assume that either of
the following conditions {\rm (C0)--(C2)} holds true:

{\rm (C0)} $\lim_{s\to\infty}a(s)=\infty$,

{\rm (C1)} $\lim_{s\to\infty}s a(s)=\infty$ and $f \in PI_{<2}$
{\rm (4.5)}\\* or

{\rm (C2)} $f \in S A\Pi_\infty ^2$ {\rm (5.4)} and $f \in
PI_{<2}$. {\rm (4.5)}\\* Write now $b\iz{V}:=f(\log|V|)$,
$A\iz{V}=a\iz{V}:=1/a(b\iz{V})$ and \beqz B\iz{V}:=
\begin{cases}
b_V, & \text{under\ conditions\ {\rm (C0)}\ or\ {\rm (C1)},} \\
b_{V}+2\nu\rrk^2 b_{V}^{-1}, & \text{under\ condition\ {\rm
(C2)}.}
\end{cases}
\eeqz Define the point process $\cN_V^\rrl$ on $[-1/2;1/2]^\nu
\times \ttR$ by \beq \cN_V^\rrl :=\sum_{k=1}^{\ve V\ve}
\delta_{\Lambda \iz{V}(k)} \quad \hbox {where} \quad \Lambda
\iz{V}(k):=\bigg(\frac {z\iz{k,V}}{|V|^{1/\nu}},\big(
\rrl\iz{k,V}-B\iz{V}\big)A\iz{V} \bigg).
\eeq Then $\cN_V^\rrl$ converges weakly to the Poisson process
$\cN$ on $[-1/2;1/2]^\nu \times \ttR$ with the intensity measure
$\dr x\times \e^{-t} \dr t$.
\end{thm}

{\it Sketch of the proof}. Conditions (6.3) and (C0) imply  that
$f \in A\Pi^{0}_{\infty}$ (4.1). Therefore, by Theorem 4.3(i) with
$p=0$, the samples $\xi_V$ satisfy the condition of Theorem~
2.2(i), consequently, $\rrl_{K,V}=\xi_{K,V}+ \mro(1)$ in
probability as $\ve V\ve \to \infty$, for fixed $K \in \ttN$.
Similarly, (6.3) and (C1) imply  that $f$ is in the classes
$A\Pi^{1}_{\infty}$ (4.1) and $PI_{<2}$~(4.5). Therefore, by
combining Theorem 4.3(i) for $p=1$, Theorems 4.5 and~3.1(iv) with
$R=0$ and $0 <\rre <\theta <1/2$, we obtain that the samples
$\xi_V$ satisfy the conditions of Theorem 2.2(ii) with $N=[\ve
V\ve^\theta]$. Consequently, $\rrl_{K,V}=\xi_{K,V}+
\rO(\xi^{-1}_{K,V})$ in probability, for fixed $K \in \ttN$. Using
these asymptotic expansion formulas for $\rrl_{K,V}$, we obtain
that in the cases (C0) and (C1) the point process $\cN_V^\rrl$
(6.19) is approximated by the corresponding point process
$\cN_V^{\xi}$ with $\xi$ instead of $\rrl$; see the proof of
Theorem~2.5(ii) in (Astrauskas 2012). Since $\cN_V^{\xi}$
converges weakly to $\cN$ (Leadbetter et al.~1983), this concludes
the proof of the theorem for (C0) and (C1).

In the case of (C2), we combine Theorem 4.3(i) for $p=2$ and
Theorems 4.5, 3.1(iv),~5.3 and 5.4 for $R=0$ and $0 <\rre <\theta
<1/2$, to find that the samples $\xi_V$ satisfy the conditions of
Theorem 2.2(iii) in probability, with $N$ and $K$ as above.
Consequently, $\rrl_{K,V}=\xi^0_{K,V}+ \rO(\xi^{-2}_{K,V})$ in
probability, where $\xi^0_{K,V}$ is the $K$th extreme value of the
i.i.d.~field $\xi^0(\cdot):= \xi(\cdot)+2\nu
\rrk^2/(\xi(\cdot)\vee 1)$ in $V$. Using this limit and applying
the same arguments as above with $\xi$ replaced by $\xi^0$, we
obtain the assertion of the theorem in the case (C2).\qed

\medskip

From Lemmas A.6 and A.13 we know that condition (C0) (resp., (C1)
or (C2)) implies that $Q(t)=\mro(t)$ (resp., $Q(t)=\mro(t^2)$ or
$Q(t)=\mro(t^3)$) as $t$ tends to infinity.

The following corollary provides some limiting distributions
for the top eigenvalues and the corresponding localization centers,
which immediately follow from the Poisson convergence results
of Theorems 6.2 and 6.9; see (Astrauskas 2007). We also
refer the reader to (Leadbetter et al.~1983, Chapter 5)
for a detailed survey on Poisson limit theorems and their
applications concerning the extremal properties of random fields.

\begin{cor}
({\it Eigenvalue statistics}). Assume that the conditions of Theorem
6.2 or Theorem 6.9 are fulfilled, with the normalizing constants 
$A\iz{V}>0$ and $B\iz{V}$ specified herein. Then, for fixed $K\dly 1$ and $V\uparrow \ttZ^\nu$, we have the following assertions:

(i) The normalized spectral gaps
$$
(\rrl\iz{1, V}-\rrl\iz{2, V})A\iz{V}, \ldots , (\rrl\iz{K-1,
V}-\rrl\iz{K, V})A\iz{V}, (\rrl\iz{K, V}-B\iz{V})A\iz{V}
$$
are asymptotically mutually independent and
have limiting joint distributions with the density
$$
\exp\lf \{ -t_1-\cdots -(K-1)t_{K-1}-Kt_{K}-\e^{-t_K}\rg \}
$$
for all $t_k\dly 0$ ($1\mly k\mly K-1$) and all $t\iz{K}\in \ttR $.

(ii) The normalized localization centers
$$
\frac {z\iz{\tau(1),V}}{\ve V\ve^{1/\nu}}, \ \frac
{z\iz{\tau(2),V}}{\ve V\ve^{1/\nu}},\ldots,\frac
{z\iz{\tau(K),V}}{\ve V\ve^{1/\nu}}
$$
are asymptotically mutually independent, and each of them is asymptotically uniformly distributed
on $[-1/2;1/2]^\nu$.

(iii)  As a consequence of (ii), the distance between the
localization centers is of order $\ve V\ve ^{1/\nu}$, i.e., for
all $1\mly l<k \mly K$,
$$
\big\ve z\iz{\tau(k),V}-z\iz{\tau(l),V}\big\ve \asymp \ve V\ve
^{1/\nu} \quad \mbox{ in\ probability}.
$$
\end{cor}
\qed

\medskip

We now give three examples of distributions $1-\e^{-Q}$, where
$\log Q(t)=\mro(t)$.

\begin{exmp}
{\rm (Astrauskas 2012).} {\it Polynomially decaying
distributions}. For some $\beta >0$, assume that $f\circ \log \in
RV_{1/\beta}$ or, equivalently, $\e^Q \in RV_{\beta}$. The latter
is the sufficient and necessary condition for the distribution
$1-\e^{-Q}$ to be in the domain of attraction of the max-stable
Fr\'{e}chet law $G_\beta (\cdot)$, or equivalently, the following
limit holds true:
$$
\lim\limits_V|V| \bP \big( \xi(0)
>tf(\log\ve V\ve) \big)=t^{-\beta} \quad \mbox{for\ all} \quad t \in \ttR_{+};
$$
see, e.g., (de Haan and Ferreira 2006, Chapter 1). Since $f \in A\Pi^{0}_{\infty}$ (4.1), from Theorem~4.3(i) with
$p=0$ and Theorem 2.2(i) we see that $\rrl_{K,V}=\xi_{K,V}+
\mro(1)$ in probability, for fixed $K \in \ttN$. Using this limit
and denoting $B_V\equiv 0$ and $a_V=1/f(\log \ve V\ve)$, we obtain
similarly as in the proof of Theorem~6.9(C0) that the point
process $\cN_V^\rrl$ (6.19) converges weakly to the Poisson
process on $[-1/2;1/2]^\nu \times \ttR_+$ with the intensity
measure $\beta \dr x\times t^{-\beta-1} \dr t$.
\end{exmp}

\medskip

\begin{exmp}
{\rm (Grenkova et al.~1990; Astrauskas and Molchanov 1992;
Astrauskas 2008).} {\it Weibull distributions}. Let $Q(t)=t^\rra$
for $t\dly0$, where $\rra >0$. For $\rra \dly 1$, the function $Q$
satisfies conditions (6.7) and (6.8) of Corollary 6.4. For $\rra
< 3$, the inverse function $f(s):=Q^\la(s)=s^{1/\rra}$ ($s\dly 0$)
satisfies conditions (6.3) and (C2) of Theorem~6.9. Therefore, by
Theorems 6.2, 6.9 and Corollary 6.4,
$$
\lim_V |V| \bP \Big( \xi(0)
>b_V+\frac{t}{a_V}\Big) =\lim_V  |V|
\bP \Big(\wt \rrl(0) > B\iz{V} +\frac{t}{a_V} \Big) = \e^{-t}\ \
(t \in \ttR).
$$
Consequently, the point processes $\cN_V^\xi$ (6.2) and
$\cN_V^\rrl$ (6.5) converge weakly to the same Poisson process
$\cN$ as in Theorem 6.2, where the normalizing constants can be
chosen as follows: $b_V=(\log\br{V})^{1/\rra}$, $A_V
=a_V=Q^{\prime}(b_V)=\rra b_V^{\rra-1}$ and

(a) $B\iz{V}= b\iz{V}$ if $\rra < 2$,

(b) $B\iz{V}=b_V+2\nu \rrk ^2b_V^{-1}$ if $2 \mly \rra <3$ \\*
and, as $\br{V}\to \infty $,

(c) $B_V= b_V+ 2\nu\rrk^2 b_V^{-1}+ \rO( b_V^{-\frac {\rra
+1}{\rra -1}})$ if $\rra \dly 3$. \\* For $\rra \dly 3$,
asymptotic equations for $B\iz{V}$ are given in (Astrauskas 2008,
Section 6).

In the case $\rra < 1$, we obtain the following almost sure
asymptotic bounds for the eigenvalues and their spacings for any
fixed $K \in \ttN$ and $m \in \ttN \backslash \{1\}$:
$$
\lsup _V \frac{\vert\rrl \iz{K, V}-b\iz{V}\vert a\iz{V}}
{\log_2\br{V}}=\frac{1}{K},
$$
$$
\linf_V\bigg (\log \Big((\rrl \iz{K, V}\!-\!\rrl \iz{K+1, V}
)a\iz{V}\Big)+\sum _{i=2}^{m-1}\log_i \br{V} \bigg)\Big
/\log_m\br{V}\!=\!-1
$$
and
$$
\lsup _V \bigg ((\rrl \iz{K, V}\!-\!\rrl \iz{K+1, V}
)a\iz{V}-\frac{1}{K}\sum _{i=2}^{m-1}\log_i \br{V} \bigg)\Big
/\log_m\br{V}\!=\frac{1}{K},
$$
with $a_V$ and $b_V$ as above; here $\log_m:=\log_{m-1}(\log)$ for
$m \dly 2$. For any $\rra
>0$, these strong limits for $\xi$ instead of $\rrl$ are proved in
(Astrauskas 2006, Section 3). Therefore, the case of $\rrl$ and
$\rra <1$ is derived by the same arguments as in the proof of
Theorem~6.9(C1) above, where one explores Theorem 4.4(i) ($p=1$)
instead of Theorem 4.3(i).
\end{exmp}

\medskip

\begin{exmp}
{\rm (Astrauskas 2013; 2016).} {\it Distributions with
fractional double exponential tails}. Let $Q(t)=\e^{ t^{\rrg}}$
for $t \dly t_0$, where $0< \rrg <1$. Obviously, $Q$ satisfies
conditions (6.7) and (6.8) of Corollary 6.4, therefore,
$$
\lim_V |V| \bP \Big( \xi(0)
>b_V+\frac{t}{a_V}\Big) =\lim_V  |V|
\bP \Big(\wt \rrl(0) > B\iz{V} +\frac{t}{a_V} \Big) = \e^{-t}\ \
(t \in \ttR).
$$
Consequently, the point processes $\cN_V^\xi$ (6.2) and
$\cN_V^\rrl$ (6.5) converge weakly to the same Poisson process
$\cN$ as in Theorem 6.2; here
$$
b_V=(\log \log \br{V})^{1/\rrg},\ \  A_V =a_V=Q^{\prime}(b_V)=\rrg
b_V^{\rrg-1}\log \ve V\ve
$$
and
$$
B_V=b_V + c_1\frac {b_V ^{\rrg-1}}{\log b_V} + c_2 \frac {b_V
^{\rrg-1}\log \log b_V}{(\log b_V )^{2}}+ c_3 \frac {b_V
^{\rrg-1}}{(\log b_V )^{2}} \big(1+\mro(1) \big)
$$
as $|V|\to\infty$, where $c_1:=\nu \rrk^2 \rrg (1-\rrg)^{-1}$,
$c_2:=c_1(\rrg-1)^{-1}$ and $c_3:=c_2\log
\frac{2(1-\rrg)\sqrt{\e}}{\rrk\rrg}$. The last formula and the
asymptotic equations for $B_V$ are derived in (Astrauskas 2016);
see also (Astrauskas 2013, Section 3).
\end{exmp}

\medskip

\subsection{\normalsize{\textsf {Distribution tails lighter than the double exponential function} }}

Throughout this subsection, we assume that the upper tails $\e^{-Q}$ 
are lighter than the double exponential function. This class of 
distributions includes fractional double exponential tails (1.8) with $\rrg > 1$ 
and bounded tails ($t_Q < \infty$).

We start with the second order asymptotic expansion formula for
the largest eigenvalues $\rrl\iz{K,V}$ of $\cH_V=\rrk \Delta_V
+\xi_V$.

\begin{thm}  Let $\e^f \in RV_0$, and fix
$K \in \ttN$. Then with probability~1
$$
\lim_V (\rrl\iz{K, V}-f\lf (\log \vert V\vert\rg)) =2\nu \rrk.
$$
\end{thm}

\begin{proof} We apply part (i) of Theorem 3.1, where $R \in \ttN$ is
fixed, $\theta_R(\cdot)\equiv \theta$ is a constant, and $m:=\ve
\ttB_R\ve$. Thus, with probability one there is a ball
$\ttB_R(z_V)$ such that $\xi(\cdot) \dly L_{V, \theta}$ in
$\ttB_R(z_V)$ for some $\theta \in \big (\frac {m-1}m; \frac
m{m+1}\big )$ and each $V\supset V_0(\omega; R)$. Therefore, by
Theorem 4.4(iii) and Lemma A.3(ii) with $\rho=0$, almost surely
$$
\xi(\cdot)-\xi_{1,V} \dly L_{V, \theta}-L_{V, 0}+ \mro(1)=\mro(1)
\ \mbox{uniformly\ in}\ \ \ttB_R (z_V),
$$
as $\ve V\ve \to \infty$, for any $R \in \ttN$. The latter means
that with probability one the samples $\xi_V$ satisfy the
condition of Theorem 2.6, therefore, $\rrl\iz{K, V}=L_{V,0}+2\nu
\rrk+\mro(1)$, as claimed. Theorem 6.14 is proved.\qed
\end{proof}

\medskip

From the proof of Theorems 2.6 and 6.14 we see that the top
eigenvalue $\rrl\iz{K, V}$ of $\cH_V$ is approximated by the
corresponding eigenvalue of the Hamiltonian restricted to the (random)
relevant regions $\ttB_{\rm{opt}}:=\ttB_{R_V} (z_V)\subset V$
where $\xi(\cdotp)$ is close to $\xi_{1,V}$ and the diameter of
which tends to infinity as $\ve V\ve \to \infty$.

Bishop and Wehr~(2012) have obtained more accurate asymptotic
bounds for the principal eigenvalue $\rrl\iz{1,V}$ of the
one-dimensional Schr\" odinger operators in $V\subset \ttZ$, with the
Bernoulli i.i.d.~potential. In particular, their results imply
the following

\medskip

\begin{thm} {\rm (Bishop and Wehr 2012).} Let $\nu
=1$, and suppose that the random sequence $\xi(\cdot)$ has a common
Bernoulli distribution: $a=\bP(\xi(0)=1)$ and $1-a=\bP(\xi(0)=0)$, so that $t_Q=1$.
Then with probability one
$$
\rrl\iz{1,V}=t_Q+ 2\rrk+(\log \ve V\ve)^{-2} (-D +\mro(1))\quad
\hbox{as}\quad V\uparrow \ttZ,
$$
where $D=D(\rrk,a) >0$ is the universal constant depending on
$\rrk$ and $a$.
\end{thm}
\qed

\medskip

In the proof of this theorem, the authors have established that
almost surely the relevant region $\ttB_{\rm{opt}} \subset V$ is
the longest consecutive sequence of sites in $V$ with $\xi(\cdot)$
equal to 1, so that the size of $ \ttB_{\rm{opt}}$ is of order
$\log \ve V\ve$. See, e.g., the review paper by Binswanger and
Embrechts (1994) for the strong and weak limit theorems for the
length $\ve \ttB_{\rm{opt}} \ve$ as $V\uparrow\ttZ$.

Recently, Germinet and Klopp (2013) have proved the Poisson limit
theorem for the top eigenvalues under nonlinear renormalization,
i.e., for the so-called unfolded eigenvalues. Write, as above,
$N(\rrl)$ ($\rrl \in \ttR$) for the integrated density of states,
i.e., the nonrandom distribution function of eigenvalues defined
as the almost sure limit of the empirical distribution function
$N_V(\rrl):=\# \{k \colon \ \rrl \iz{k,V} \mly \rrl \}/ \ve V\ve$
as $\ve V\ve \to \infty$ (Kirsch 2008).

\medskip

\begin{thm} {\rm (Germinet and Klopp 2013, Theorem 2.3).} Assume
that $\cH_V$ is the one-dimensional Anderson Hamiltonian ($\nu
=1$), where
the potential is bounded ($\ve \xi(0)\ve \mly \const$) with the
distribution density $p(\cdot):=\e^{-Q(\cdot)}Q^{\prime}(\cdot)$.
Assume, in addition, that the density $p(\cdot)$ is bounded and
does not decay too fast at $t_Q$ (say, $p(t)=\e^{-(t_Q
-t)^{\mo(1)}}$ or $p(t)=\e^{-(t_Q -t)^{-\vartheta}}$ as $t\uparrow
t_Q$, for $0 < \vartheta < 1/2$). Define the point process
$\cM_V^\rrl$ on the positive half-axis $\ttR_+$ by
$$
 \cM_V^\rrl
:=\sum_{k=1}^{\ve V\ve} \delta_{\ve V \ve(1-N(\rrl \iz{k, V}) )}.
$$
Then $\cM_V^\rrl$ converges weakly to the Poisson process on
$\ttR_+$ with the intensity measure $\dr t$, the Lebesgue measure.
\end{thm}
\qed

\medskip

For $\nu \dly 2$, this Poisson limit theorem was shown to hold if
the Laplacian $\Delta $ is replaced by some translation invariant
operator $\cT\psi(x)=\sum _{y\in \ttZ^\nu}T(y)\psi(x-y)$, where
$T(\cdotp )$ is  a real nonrandom function decaying exponentially at infinity.

The proof of Theorem
6.16 relies on the improved versions of Wegner and Minami
estimates that control the structure of eigenvalues  $\rrl \iz{k,V}$ in
a small neighborhood $I$ of the upper spectral edge, so the amount
$N(I)$ is allowed to be exponentially small in $\ve I\ve^{-1}$;
see Section 1.3 for more explanations. See also (Minami 2007) for
a detailed background of Poisson convergence results for unfolded
spectral values. It is important for applications that this
convergence result is given in terms of the integrated density of
states, the main quantity in the theory of random Schr\" odinger
operators. However, in the proof of Theorem 6.16, neither extreme
value theory, nor links to the asymptotic geometric properties of
random potential are explored.

\medskip

\begin{cor}
Assume
that $\cH_V$ is the one-dimensional Anderson Hamiltonian ($\nu
=1$), where the potential has the distribution density
$p(\cdot)$ satisfying the conditions of Theorem 6.16. Then, for fixed
$K\dly 1$ and $V\uparrow \ttZ$, we have the following asymptotic formulas in probability:

(i) If $p(t)=\e^{-(t_Q -t)^{\mo(1)}}$ as $t\uparrow t_Q$, then
$$
\rrl\iz{K,V}=t_Q+ 2\rrk-(\log \ve V\ve)^{-2+\mo(1)}.
$$

(ii) If there is $0 < \vartheta < 1/2$ such that $p(t)=\e^{-(t_Q
-t)^{-\vartheta}}$ as $t\uparrow t_Q$, then
$$
\rrl\iz{K,V}=t_Q+ 2\rrk-(\log \ve
V\ve)^{-2/(1+2\vartheta)+\mo(1)}.
$$
\end{cor}
\qed

\medskip

Corollary 6.17 can be proved by using the limit theorems for the
unfolded eigenvalues and asymptotic expansion formulas for the
tails $1-N(\cdot)$ at the upper spectral edge derived, e.g., in
(Klopp 1998; Biskup and K\"onig 2001). See also (Klopp 2000) and
Section 3.5 of (Kirsch and Metzger 2007) for a detailed discussion
on the edge asymptotics of the integrated density of states. Thus,
the bifurcations in the asymptotic behavior of the top eigenvalues
are caught by those in the tail behavior of the integrated density
of states, i.e.~\lk Lifshits tails\rk. Notice that the asymptotics
in Corollary 6.17(i) agree with the results in the Bernoulli case
(Theorem 6.15). Meanwhile, if the tails $\e^{-Q}$ decays faster at
$t_Q$, then the fluctuations of $\rrl\iz{K,V}$ are much sharper
(Corollary 6.17(ii)) than those in the Bernoulli case (Theorem
6.15).

For the tails lighter than the double exponential function including the case $t_Q< \infty$, some
heuristics on the asymptotic formulas for $\rrl \iz{1,V}$ ($\nu
\dly 1$) and their relations to the long-time asymptotic formulas
for the parabolic Anderson model have earlier been discussed by
Biskup and K\"onig (2001), and van der Hofstad et al.~(2006). Their
assumptions are given in terms of scaling and regularity
properties of the cumulant generating function $\log\bE
\e^{t\xi(0)}$ as $t \to \infty$.

Recently, Biskup et al.~(2014) have proved the \lk homogenized\rk\
versions of limit theorems for the largest eigenvalues of the
(scaled) finite-volume discrete  Schr\"odinger operators
$\cH^{(\rre)}$ with a bounded random potential. Their results assert
that, as the scale parameter $\rre$ tends to zero, then 1) the
largest eigenvalues of $\cH^{(\rre)}$ converge in probability to the corresponding eigenvalues of the limiting (nonrandom)
finite-volume continuous Schr\"odinger operator, and 2) the
fluctuations of the largest eigenvalues centered by their means
are Gaussian in limit.

We notice that the conditions of the above statements imply that
$t^{-1} \log Q(t) \to \infty$ as $t\uparrow t_Q
> 0$ (see Lemma~A.3 with $\rrr =0$ below).

\medskip

We end this subsection with a discussion on the following
important model of spatially continuous random Schr\" odinger
operators:

\begin{exmp}
({\it Schr\" odinger operators in $\ttR^\nu$ with a bounded Poisson
potential of obstacles}). Let $\Delta^{{\rm cont}}$ be the
$\nu$-dimensional continuum Laplacian. Define the random potential
$\xi(\cdot)$ by \beq \xi(x)=-\sum_iW(x+x_i)\quad (x \in \ttR^\nu);
\eeq here $\{x_i\}$ is a Poisson point process in $\ttR^\nu$ with the
constant intensity $\mu
>0$; $W(\cdot)$ is a fixed nonnegative compactly supported,
bounded measurable function, $W(\cdot)$ is non-identically zero Lebesgue-a.e. The
potential $\xi(\cdot)$ is known as a {\it Poisson field of \lk
soft\rk \ obstacles}. Denote by $V:=[-s;s]^\nu \subset \ttR^\nu$ the cubes of the
volume $\ve V\ve$ such that $V \uparrow \ttR^\nu$.
Let us consider the principal Dirichlet eigenvalue
$\rrl\iz{1,V}^{{\rm cont}}<0$ of the operator $\Delta^{{\rm cont}}
+\xi(\cdot)$ in $V$. The eigenvalue $\rrl\iz{1,V}^{{\rm cont}}$
satisfies the following asymptotic formula (e.g., Sznitman 1998):
As $V \uparrow \ttR^\nu$, almost surely \beq \rrl\iz{1,V}^{{\rm
cont}}=(\log \ve V\ve)^{-2/\nu} (-C(\nu,\mu)+\mro(1)),
\eeq where $C(\nu,\mu) >0$ is the universal constant depending on
$\nu$ and $\mu$; see below.

\medskip

Let us {\it sketch a derivation of the lower bound} for
$\rrl\iz{1,V}^{{\rm cont}}$. By the monotonicity property of
eigenvalues, we have the bound $\rrl\iz{1,V}^{{\rm cont}} \dly
\rrl_{1,\ttA}^0$, where $\rrl_{1,\ttA}^0$ is the local
principal Dirichlet eigenvalue of the operator $\Delta^{{\rm
cont}} +\xi(\cdot)$ restricted to the obstacle-free connected open region $\ttA
\subset V$. We now maximize the eigenvalue $\rrl_{1,\ttA}^0$
over such $\ttA$. First, since the probability of region
$\ttA$ to have no points $\{x_i\}$ is equal to $\e^{-\mu\ve \ttA\ve}$
and the number of disjoint shifts of $\ttA$ is of order $\ve V\ve$, we obtain
from the Borel-Cantelli lemma that the volume $\ve \ttA\ve$ should
be approximately equal to $\mu^{-1}\log\ve V\ve$; cf.~the proof of Theorem~3.1 above. On the other
hand, we have from the Faber-Krahn inequality that the principal Dirichlet eigenvalues $\rrl_{1,\ttA}^0$ of
the Laplacian $\Delta^{{\rm cont}}$ in regions $\ttA$ of the
constant volume achieve their maximum at the ball. Thus,
$\rrl\iz{1,V}^{{\rm cont}} \dly \rrl_{1,\ttB\iz{{\rm opt}}}^0(1+\mro(1))$, where
$\rrl_{1,\ttB\iz{{\rm opt}}}^0$ is the principal Dirichlet eigenvalue of
the operator $\Delta^{{\rm cont}}$ in the ball $\ttB_{{\rm opt}}:=\ttB_{R_V}(z_V)$
centered at some random $z_V \in V$ with the radius $R_V:=\lf (\ve
\ttB_1\ve^{-1}\mu^{-1} \log \ve V\ve \rg )^{1/\nu}$, and $\ve
\ttB_1\ve$ is the volume of the unit ball $\ttB_1 \subset
\ttR^\nu$. Consequently, as $V\uparrow \ttR^\nu$, almost surely
$$
\rrl\iz{1,V}^{{\rm cont}} \dly \rrl_{1,\ttB\iz{{\rm opt}}}^0(1+\mro(1))=R_V^{-2}\lf
(\rrl_{1,\ttB_1}^0 +\mro(1)\rg )=(\log \ve V\ve)^{-2/\nu}
(-C(\nu,\mu)+\mro(1)),
$$
where $C(\nu,\mu):=-\rrl_{1,\ttB_1}^0 (\ve \ttB_1\ve \mu
)^{2/\nu}>0$. This lower bound can be shown to be equal to the
upper bound for $\rrl\iz{1,V}^{{\rm cont}}$, concluding the proof
of (6.21). Notice also that the rough upper bound for
$\rrl\iz{1,V}^{{\rm cont}}$ can be derived by using the spatially
continuous version of Lemma 2.8 above.

Summarizing, we conclude that the principal Dirichlet eigenvalue
$\rrl\iz{1,V}^{{\rm cont}}$ is
approximated, as $V \uparrow \ttR^\nu$, by the local principal Dirichlet eigenvalue
$\rrl_{1;\ttB\iz{{\rm opt}}}^0$ of the operator restricted to the relevant region $\ttB_{{\rm opt}}:=\ttB_{R_V}(z_V)\subset V$, so that $\rrl\iz{1,V}^{{\rm
cont}}\leftrightarrow \ttB_{{\rm opt}}$. These observations
and formulas agree with the corresponding formulas for the
discrete Anderson models in $\ttZ$ with the Bernoulli i.i.d.~potential
(Theorem 6.15).

As already mentioned, formula (6.21) and the more explicit
asymptotic bounds for the principal Dirichlet eigenvalue
$\rrl\iz{1,V}^{{\rm cont}}$ ($V \uparrow \ttR^\nu$) were proved
by Sznitman (1998) exploring his original method of enlargement of
obstacles. By this method, the geometry of the spatial regions
where $\xi(\cdot)> 0$ and $\xi(\cdot)\equiv 0$ is reduced to the
simpler geometry of regions associated with the modified potential
in the spectral problems, without changing the eigenvalues very
much. Finally, notice that the asymptotic bounds for the principal
eigenvalues are crucial for study of the intermittent behavior of
a Brownian motion in a Poisson field of obstacles; cf.~Section 7
below.
\end{exmp}
\qed

\medskip

\subsection{\normalsize{\textsf {The double exponential tails} }}

In the double exponential case (1.14), G\" artner and Molchanov (1998,
Theorem 2.16) have obtained the second order expansion formula for
the principal eigenvalue $\rrl_{1,V}$ of $\cH_V=\rrk \Delta_V
+\xi_V$ by claiming a continuity of $Q$. We now provide their
result with the continuity condition removed.

\begin{thm}  If $\e^f \in RV_\rho$ for some
$0 < \rho < \infty$, then with probability~1
$$
\lim_V (\rrl_{1, V}-f\lf (\log \vert V\vert\rg)) =2\nu \rrk
q(\rho / \rrk),
$$
where the nonrandom function $q$ is defined in Section 2.5.
\end{thm}

\begin{proof} We check the conditions of
Theorem 2.7. First, by Theorem 4.7, almost surely $\xi_V$
satisfies condition (2.29). To prove (2.30), we fix constants $R
\in \ttN$, $\delta >0$, and write
$$
\theta(y):=\theta_R(y):= 1- \exp \lf\{\lf
(h^{\ttB_R}_{\rm{opt}}(y)-\delta\rg)/\rho \rg\} \ \ (y \in
\ttB_R),
$$
where the nonrandom function $h^{\ttB_R}_{\rm{opt}}(\cdot)$ is
defined in Section 2.5. Consequently, $\theta(\cdot)$ satisfies
the assumptions of Theorem 3.1(i). Combining the statements of
Theorem~3.1(i), Lemma A.3(ii) and Theorem 4.4(iii), we obtain the
following assertion with probability one: for any $V \supset
V_0(\omega ; \delta, R)$ there is $z_V \in V$ such that
$$
\xi(y) \dly L _{V, \theta(y-z\iz{V})}  \dly \xi _{1,V}
+h^{\ttB_R}_{\rm{opt}}(y-z_V) -2\delta \ \ \mbox{for\ all} \ \ y \in
\ttB_R(z_V).
$$
Since $R \in \ttN$ and $\delta >0$ are arbitrary constants, this
estimate concludes the proof of the almost sure limit~(2.30). Now,
Theorems 2.7 and 4.4(iii) imply the assertion of Theorem 6.19.

\qed
\end{proof}

\medskip

From the proof of Theorems 2.7 and 6.19 we see that almost surely
the eigenvalue $\rrl_{1, V}$ approaches (as $\ve V\ve \to
\infty$) the local principal eigenvalue in the random region,
where $\xi(\cdotp)\approx\xi _{1,V}
+h^{\ttB_R}_{\rm{opt}}(\cdotp)$ for $R$ arbitrarily large, so that
$\rrl_{1, V}$ is associated with the (random) \lk relevant
island\rk \ of high $\xi_{V}$-values of optimal shape, the
diameter of which is asymptotically bounded. From Theorem
3.1(iv), Remark 3.5 and the last assertion of Lemma A.3,
it follows that the \lk islands\rk \ of
$\xi_{V}$-extremes are located asymptotically far away from each
other. Moreover, if the constant $\rho / \rrk$ is large enough,
these \lk islands\rk \ are located in the neighborhood of single
extremely high $\xi_V$-peaks; see (Astrauskas 2008, Theorem 4.4
and Corollary 4.5) and (Astrauskas~2013, Theorem 2.1(iii)).

For arbitrary $0 < \rho < \infty$, Poisson limit theorems for the
largest eigenvalues and the corresponding localization centers are
proved by Biskup and K\"onig (2016); see also the survey by
K\"onig (2016) on these limit theorems and related
topics. To formulate their result, we again define the sites
$z_{\tau(k),V} \in V$ by (6.18), i.e., the localization centers of
the $k$th eigenfunctions $ \psi(\cdotp; \rrl \iz{k,V})$ ($1 \mly k
\mly |V|$) of the Hamiltonian $\cH_V=\rrk \Delta_V +\xi_V$.

\medskip

\begin{thm} {\rm (Biskup and K\"onig 2016, Theorem 1.2).} Let $t_Q = \infty$, and assume that $Q$ is a continuously differentiable function such that
$$
\lim_{t \to \infty}\frac{Q'(t)}{Q(t)}=\frac{1}{\rrr} \ \
\mbox{for\ some} \ \ 0 < \rrr < \infty.
$$
Then the following assertions (I) and (II) hold true:

{\rm (I) (Poisson limit theorem)} There are constants $ B_{V}=
f(\log \ve V\ve)+2\nu \rrk q(\rho / \rrk)+\mro(1)$ and $A_{V}=\rrr^{-1} \log \ve
V\ve$ such that the
point process $\cN_V^\rrl$ (6.5) converges weakly to the Poisson process $\cN$ on
$[-1/2;1/2]^\nu \times \ttR$ with the intensity measure $\dr
x\times \e^{-t} \dr t$.

{\rm (II) (Exponential localization)}  As $V\uparrow \ttZ^\nu$ and $K \dly 1$ fixed,
we have with probability $1+\mro(1)$  that
there exist non-random constants $C>0$, $M >0$, $C^\prime>0$ and $M^\prime >0$ such that
$$
\vert \psi(x; \rrl
\iz{K,V})\vert \mly C \exp \{ -M\vert x-z\iz{\tau (K), V}\vert
\}\quad \mbox{for\ all}\quad x\in V,
$$
and
$$
\vert \psi(x; \rrl \iz{K,V})\vert \mly C^\prime \exp \{ -M^\prime
(\log \log \ve V\ve) \big\vert x-z\iz{\tau (K), V}\big\vert
\}\quad \mbox{for}\quad \ve x-z\iz{\tau(K),V}\ve \dly \log\ve
V\ve,
$$
i.e., the $K$th eigenfunction is highly concentrated in the
neighborhood of its localization center.
\end{thm}

\medskip

For sufficiently large $\rrr$, i.e. $\rrr >\rrr_0$, Poisson limit theorems and localization theorems for
the top eigenvalues were earlier proved by Astrauskas (2007; 2008;
2013). For $\rrr >\rrr_0$, the corresponding localization
properties present an interesting intermediate case between the
single site concentration property, i.e. $\rrl_{K,
V}\leftrightarrow z\iz{\tau(K), V}$, in the case $\rrr=\infty$
(Theorem 6.2(II)) and the non-single site concentration property,
i.e. $\rrl_{K, V}\leftrightarrow\ttB ^{K}_{\rm opt}$, for $0
\mly \rrr < \rrr_0$ (Section 6.2 and Theorem~6.20).

From Remark A.5(i) and Lemma A.3, we notice that the conditions of
Theorem~6.20 imply $\e^f \in RV_{\rrr}$, i.e., the assumption of
Theorem 6.19. Moreover, from Remark~A.5(ii), Lemma A.1 and Lemma
6.1 with $a(\cdot)\equiv 1/Q^{\prime}(\cdot)$ we also see that the
conditions of Theorem 6.20 yield the limit (6.1) with $b_{V}=
f(\log \ve V\ve)$ and $a_{V}=A_{V}=\rrr^{-1} \log \ve V\ve$. Consequently, the distribution $1-\e^{-Q}$ is in the domain of
attraction of the max-stable Gumbel law $G_{\exp}(\cdot)$. We
finally notice that the conditions of Theorem~6.19 or Theorem~6.20
imply the limit $f(s)=(\rrr +\mro(1))\log s$ as $s \to \infty$,
which in turn is equivalent to $\log Q(t) = \rrr^{-1}t+\mro(t)$
as $t \to \infty$; see Lemma A.3.

\medskip

\subsection{\normalsize{\textsf {Some comments on the proofs} }}

In this section, we briefly comment and compare the proof of
Theorems~2.2, 2.3, 6.2, and 6.9 by Astrauskas and Molchanov (1992)
and Astrauskas (2007; 2008; 2012; 2013) (\lk relevant single
peak\rk approximation) and the proof of Theorem 6.20 by Biskup and
K\"onig (2016) (\lk relevant island\rk approximation).

\textbf{(RSP)} {\it \lk Relevant single peak\rk approximation.} As
already mentioned in Section~1.2, the proof of Theorems~2.2, 2.3,
6.2, and 6.9 is based on the finite-rank perturbation arguments and
the analysis of Green functions involving the cluster expansion over
paths. To be more precise, fix $Z:=z_{\tau(K),V}\in V$, the
localization center of the $K$th eigenfunction, and denote by $\rrl(Z)$ the extreme eigenvalue of $\cH_V=\rrk
\Delta_V +\xi_V$ associated with the site $Z$. Let $\cG_{V}^{(Z)}(\rrl; \cdotp ,\cdotp)$ be the Green function of the
Hamiltonian $\rrk \Delta_V +(1-\delta _{Z})\xi_V$ on $l^2(V)$.
Under the conditions of Theorem~2.2 or 2.3 (i.e.~sparseness and
difference in height of $\xi_V$-peaks as $V \uparrow \ttZ^\nu$),
the eigenvalue $\rrl(Z)$ is a solution to
the dispersion equation \beq \cG_{V}^{(Z)}(\rrl; Z,Z)=\frac{1}{\xi(Z)}
\eeq and the corresponding eigenfunction is $\cG^{(Z)}_V(\rrl(Z);
\cdotp ,Z)$. By expanding the Green function $\cG_{V}^{(Z)}(\rrl;
\cdotp ,\cdotp)$ over paths, one proves that equation~(6.22) is
approximated by the corresponding equation $\wdt \cG_{V}(\rrl;
Z,Z)=1/\xi(Z)$ for the principal eigenvalue of the \lk single
peak\rk Hamiltonian $\rrk \Delta _{V}+\wdt
\xi_{V}+\xi(Z)\delta _Z$; here $\wdt \cG_{V}(\rrl;
\cdotp,\cdotp)$ stands for the Green function of the operator
$\rrk \Delta _{V}+\wdt \xi_{V}$. Again expanding $\wdt
\cG_{V}(\rrl; \cdotp,\cdotp)$ over paths, one finds that the
eigenvalue $\rrl(Z)$ of $\cH_{V}$ is approximated by a certain
(nonlinear) function on $\wdt \xi _{V}$ and $\xi(Z)$; cf.
(2.20)--(2.22). Moreover, because of the sparseness of
$\xi_V$-peaks, the extreme eigenvalues $\rrl(Z)$ become
asymptotically independent, so that they obey asymptotic Poisson
behavior as $V\uparrow \ttZ^\nu$ (see Theorems~6.2 and 6.9).

We notice that the analysis of the Green functions combined with
the finite-rank perturbation theory is essential to study the
largest eigenvalues of the finite-volume operators $\cH_{V}$ in
the \lk relevant single peak\rk approximation. Recall that these
techniques also play a crucial role in the proof of the Anderson
localization for the infinite-volume Hamiltonian $\cH$ (Kirsch
2008; Stolz 2011); see also Section 1.3 above.

\textbf{(RI)} {\it \lk Relevant island \rk approximation.}
Recently, Biskup and  K\"onig (2016) have developed novel
arguments to prove Poisson limit theorems for the largest
eigenvalues in the case of double exponential tails (see Theorem
6.20 above). As in the single-peak approximation, the analysis of
the extreme eigenvalues is here based on controlling the
dependence of an eigenvalue on the geometric properties of
$\xi\iz{V}$-peaks and the associated regions in $V$. This enables
to identify the \lk relevant\rk regions $\ttB ^{K}_{\rm
opt}:=\ttB_{R_V}(z_V^K) \subset V$ (where $\xi (\cdotp)$ is high
and of the optimal shape) such that the $K$th largest eigenvalue
$\rrl_{K,V}$ of $\rrk \Delta_V +\xi_V$ is approximated by the
local principal eigenvalue $\rrl_{1,\ttB^{K}_{\rm opt} }$ of the
Hamiltonian restricted to $l^2(\ttB ^{K}_{\rm opt})$ (cf.~also Theorems~2.7
and 6.19 and their proofs in the present survey). In other words,
the eigenvalues associated with a block of \lk relevant islands
\rk of high $\xi_V$-values can be determined by the local
principal eigenvalues associated with separate \lk relevant
islands \rk. It is worth noticing that the conditions of Theorem 6.20
imply that the islands of high $\xi_V$-values are located
extremely far from each other as $V\uparrow \ttZ^\nu$. (See also
Theorem~3.1 for the related limits under the continuity assumption
(3.1)). The proof of Theorem 6.20 involves the following
procedures on a simplification of potential configurations: 1)
those regions, where the potential possesses the lower values, are
deleted from $V$ ({\it domain truncation and component trimming});
2) for the radius $R_V$ tending to infinity slowly, the analysis
of the local principal eigenvalues in all balls $\ttB_{R_V}(z)
\subset V$ is reduced to the consideration of independent
identically distributed local principal eigenvalues in disjoint
balls in $V$ ({\it coupling to i.i.d.~variables}); 3) the local
principal eigenvalue in the region $\ttB ^{K}_{\rm opt}$ is
separated from other local eigenvalues in $\ttB ^{K}_{\rm opt}$
({\it reduction to one eigenvalue per component}); 4) extremal
type limit theorems for the local principal eigenvalues in
$\ttB_{R_V}(z)$ are comparable to each other for the different
increase rate of $R_V \to \infty$, with the same normalizing
constants $A_ {V}$ and $B_ {V}$ ({\it stability with respect
to partition side}), and so on.

Summarizing, the main idea of the proof of Theorem 6.20 explores
the straightforward geometric arguments controlling the dependence
of eigenvalues on potential configurations, rather than the
techniques of resolvents or Green functions. This is in contrast
to the relevant single peak approximation in (RSP), where the Green
functions are the main object of analysis. On the other hand,
although most of the proof of Theo\-rem~6.20 is based on
deterministic arguments, we are not able to reformulate this
assertion in terms of $\xi_V$-extremes (like in Theorems 2.2--2.7
above), except for the case of sufficiently large $\rrr$ considered in (Astrauskas 2008, Theorem B.3).

\medskip

\setcounter{equation}{0} \addtocounter{sct}{1}


\section{\normalsize{\textsf {APPLICATIONS TO THE PARABOLIC ANDERSON MODEL} }}

\setcounter{prop}{1}

\subsection{\normalsize{\textsf {The parabolic Anderson model} }}

The parabolic Anderson model (PAM) is the Cauchy problem for the
following heat equation with random potential: \beq
\begin{split}
&\frac{\partial u(s, x)}{\partial s}\!=\!\rrk\!\sum _{\vert y\vert
=1}\!\lf (u(s, x\!+\!y)\!-\!u(s,x)\rg )\!+\!\xi(x)u(s, x),\ \ s
\dly 0,\ x\in
\ttZ^\nu,\\
&u(0,x)=\delta_0(x),\ \ x\in \ttZ^\nu; \mbox{\hspace{2.5cm}}
\end{split}
\eeq here, as above, $\xi(\cdot)$ is an i.i.d.~random field
(potential) with distribution~$\bP(\xi(0)>~t)=\e^{-Q(t)}$;
$\delta_0$ is the Kronecker delta function at the origin (i.e.,
the localized initial datum of the problem); the variable $s\dly
0$ is referred to as a time. The equation has almost surely a
unique nonnegative solution, provided $\xi(0)\vee 0$ has a finite
moment of order $>\nu$ (G\"artner and Molchanov 1990).

The PAM appears in the context of population dynamics, chemical
kinetics, magnetism and turbulence, etc. (e.g., G\"artner and
Molchanov 1990; Molchanov 1994). The following interpretation of
the solution $u$ is well-known in the mathematical literature
(e.g., Molchanov 1994): Let $(X(s)\colon\ s \dly 0)$ be a
continuous-time random walk in $\ttZ^\nu$ with a generator $\rrk
\Delta_{\rm dif}$, where $\Delta_{\rm dif}\psi(x):=\sum_{\ve
y-x\ve=1}(\psi(y)-\psi(x))$. Let $\xi^+(\cdot)\dly 0$ and
$\xi^-(\cdot)\dly 0$ be independent random i.i.d.~fields on
$\ttZ^\nu$, and write $\xi(\cdot):=\xi^+(\cdot)-\xi^-(\cdot)$. For
a fixed realization $(\xi^+(\cdot);\xi^-(\cdot))$, consider a
system of particles which obey the following diffusion and
branching mechanism:

1) at time $s = 0$, there is a single particle at the origin;

2) particles move independently of each other according to the
random walk $X(\cdot)$;

3) at the site $x$ a particle disappears with intensity $\xi^-(x)$ and
splits into two new 

particles with intensity $\xi^+(x)$, which
further move according to $X(\cdot)$.
\\
Then, for a fixed realization $(\xi^+(\cdot);\xi^-(\cdot))$, the
solution $u(s,x)$ to (7.1) is the expected number of particles at the
site $x$ at time $s$, where the expectation is taken over a
branching mechanism and diffusion (but not over random medium
$\xi^+(\cdot),\xi^-(\cdot)$). Thus, the sum $U(s):=\sum_{x\in
\ttZ^\nu}u(s,x)$ is the expected total mass of particles at time
$s$. We see from the Feynman-Kac formula (7.2) below that $U(s)$ is
equal to $U(s,0)$, where $U(s,\cdot)$ is the solution to equation
(7.1) with the homogeneous initial datum $U(0,\cdot)\equiv 1$ instead
of the localized one.

For i.i.d.~random potentials, the PAM exhibits an {\it
intermittency effect}:  As $s\to \infty$, the overwhelming
contribution to the total mass $U(s)=\sum_x u(s,x)$ of the
solution $u$ to (7.1) comes from a small number of spatially
separated and relatively small islands of large
$u(s,\cdot)$-values, i.e., {\it intermittent islands}. This is in
contrast to the case of constant potential $\xi(\cdot)\equiv
\const$, for which the solution $u(s,\cdot)$ is spread over
the spatial ball of radius $\rO(\sqrt s)$ as $s \to \infty$, i.e.,
{\it diffusion effect}.

Various aspects of long-time intermittent behavior of the PAM
(asymptotic expansion formulas for the total mass $U(s)$ and its
statistical moments, concentration properties for the solutions
$u(s,\cdot)$, etc.) have been intensively studied, during the last
two decades, by mathematicians Molchanov, G\"artner, Sznitman,
K\"onig, Biskup, M\"orters, den Hollander, Sidorova, van der Hofstad, and their
colleagues. See (K\"onig 2016) for a recent survey on the subject
and references therein. The main technical tools of intermittency
theory are a spectral representation and the Feynman-Kac formula for the
solutions $u$, $U$. Recall the latter formula:
\begin{eqnarray}
&&u(s,x)=\textbf{E}_x \lf [\exp \lf \{\int^s_0\xi(X(a)) \dd
a \rg \}u(0,X(s))\rg ] \n \\[6pt]
&&\hspace{.7cm} =\textbf{E}_x \lf [\exp \lf
\{\int^s_0\xi(X(a)) \dd a \rg \}\delta_0\lf (X(s)\rg )\rg ],\
s\dly 0,\ x\in \ttZ^\nu;
\end{eqnarray}
here the expectation $\textbf{E}_x$ is taken with respect to the
random walk $X(\cdot)$ in $\ttZ^\nu$ as above, conditioned by
$X(0)=x$. Looking at (7.2), we see that the intermittent behavior
is determined by competition between two factors: (i) extremely large
exponential factor in (7.2) associated with a portion of the
trajectories $X[0;s]:=(X(a)\colon\ 0\mly a \mly s)$ spending much
time at spatial regions where $\xi(\cdot)$ is high, and (ii) very small
probabilities of such trajectories in (7.2). In the model, the
potential necessitates concentration properties of $u$; meanwhile,
the Laplacian forces these properties to be less expressed. It
turns out that, as the upper tails of potential distribution get
heavier, the $\xi_V$-extremes get more pronounced as
$V=V(s)\uparrow \ttZ^\nu$ (cf.~Section 1.2); therefore, the
long-time intermittent properties (in particular, mass
concentration properties) of the PAM become stronger. In
(G\"artner and Molchanov 1998; G\"artner and den Hollander 1999;
G\"artner et al.~2007), the emphasize has been made on the double
exponential tails (1.14). Such tails indicate the critical
situation between formation of widely-spaced single peaks of $u(s,
\cdotp )$ in the case of tails heavier than in (1.14) (e.g.,
G\"artner et al.~2007; K\"onig et al.~2009; Sidorova and Twarowski
2014; Fiodorov and Muirhead 2014), and formation of widely-spaced
extremely large \lk islands\rk\ of higher values in the behavior of
$u(s, \cdotp )$ for the tails lighter than in (1.14) (Biskup and
K\"onig 2001; van der Hofstad et al.~2006). The results of these papers
suggest that the optimal strategy of particles is to move quickly
to the spatial region where the potential values are high and of
preferred shape, and to stay here for the remaining time.

In this section, we will focus on the representation of the
solutions $u$ and $U$ in the spectral terms of the Anderson
Hamiltonian $\cH_V=\rrk\Delta_V +\xi_V$ where $V=V(s) \uparrow
\ttZ^\nu$. In view of this representation, we will discuss some
techniques of the extreme value theory for eigenvalues, which can be
applied to study intermittent properties of the PAM; cf.~Theorems
7.1--7.2 below.

\medskip

\subsection{\normalsize{\textsf {Asymptotic expansion formulas for the total mass} }}

The first result in this direction was obtained by G\"artner and
Molchanov (1998), who particularly derived the second-order
asymptotic formula for the logarithm of the total mass $U(s)$
($s\to \infty$)  with probability one (and, as a by-product, the
corresponding result for the principal eigenvalue of the operator
$\cH_V$), provided the potential distribution satisfies mild RV
conditions and has all positive exponential moments finite; 
cf.~the assumptions of Theorem 7.1 below. Let us
{\it sketch the proof of the almost sure asymptotic formula} for
the PAM following the arguments of their paper. The first key
observation is that the solutions $u$ and $U$ are approximated, as
$s\to \infty$, by their finite-volume analogues $u\iz{V(s)}$ and
$U\iz{V(s)}$, respectively. I.e., $u\iz{V(s)}$ and $U\iz{V(s)}$
are solutions to the corresponding equations in $V(s)$ with the
Dirichlet boundary condition; here $V(s)\subset \ttZ^\nu$ denote
cubes centered at the origin, whose size length is of order
$s(\log s)^{c}$ for some constant $c>1$. In particular, almost
surely \beq U(s)=U\iz{V(s)}(s)+\mro (1)\ \ \mbox{and}\ \ u(s,
x)=u\iz{V(s)}(s, x)+\mro (1)\ \ \mbox{as}\ \ s\to \infty \eeq
uniformly in $x$, by using the standard cut-off procedure for the
solutions $u$ and $U$, based on the following facts: 1) the overwhelming asymptotic
contribution to the Feynman-Kac representation of $u$ and $U$ is
given by trajectories $X(\cdot)$ which stay inside the box $V(s)$ during
the whole time interval $[0; s]$ i.e. $X[0;s]\subset V(s)$; and 2) the contribution from trajectories $X(\cdot)$ visiting the complement of $V(s)$ during
the time interval $[0; s]$ is much smaller. On
the other hand, the solution $u\iz{V}$ with $V:=V(s)$ admits the
spectral representation \beq u\iz{V}(s, x)=\sum _{k=1}^{\vert
V\vert }\e^{\rrl\iz{k, V}s-2\nu \rrk s }\psi (0; \rrl\iz{k, V} )
\psi (x; \rrl\iz{k, V} ) \quad \hbox {for \ each} \ s\dly 0 \
\hbox{and} \  x\in V,
\eeq and, therefore, the total mass $U\iz{V}(s)=\sum_{x\in
V}u\iz{V}(s, x)$ has the following representation \beq
U\iz{V}(s)=\sum _{k=1}^{\vert V\vert }\e^{\rrl\iz{k, V}s-2\nu \rrk
s } \big (\psi (\cdotp ;\rrl\iz{k, V}),1\big )_{V} \psi (0;
\rrl\iz{k, V} ) \quad \hbox {for \ each} \ s\dly 0;
\eeq here $\big (\psi ,\varphi\big )_{V}$ denotes the inner
product of the functions $\psi$ and $\varphi$ in $l^2(V)$, and 1
stands for the function taking everywhere value~1. Recall that
$\rrl\iz{k, V}$ and $\psi (\cdotp ; \rrl\iz{k, V})$ are the $k$th
eigenvalue and the corresponding eigenfunction of the Hamiltonian
$\cH_{V}=\rrk\Delta_V +\xi_V$. Moreover, the eigenfunctions are
chosen to form an orthonormal basis of $l^2(V)$, and the principal
eigenfunction to be strictly positive in $V$. We may also assume,
without loss of generality, that $\psi (0 ; \rrl\iz{k, V})\dly 0$
for all $k$.

Let us prove that the main asymptotic contribution to the
logarithm of $U\iz{V}(s)$ in (7.5) comes from the first term
associated with the principal eigenvalue. Looking at (7.5) and
applying the Cauchy-Schwarz inequality and Parseval's
identity, we easily obtain the upper bound \beq U_V(s)\mly
\e^{\rrl\iz{1, V}s-2\nu \rrk s }\sqrt{\ve V \ve}
\eeq for each $s\dly 0$ and each $V$. To derive the lower bound
for $U(s)$, one needs more sophisticated arguments: Let
$\bun{V}=\bun{V}(s)$ denote centered cubes of side 
length of order $s(\log s)^{-c}$ for some $c>1$, 
and the site $z_1\in \bun{V}$ as a localization
center of the principal eigenfunction of the operator
$\cH_{\bun{V}}=\rrk\Delta_{\bun{V}} +\xi_{\bun{V}}$. Recall the Feynman-Kac
representation for $U(s)$:
$$
U(s)=\textbf{E}_0 \lf [\exp \lf \{\int^s_0\xi(X(a)) \dd a
\rg \}\rg ],\ \ s\dly 0,
$$
so that $U(s) \dly \underline{U}(s)$, where $\underline{U}(s)$ is
the same expectation $\textbf{E}_0[\cdot]$ when restricted to the particle
trajectories $X[0;s]$, which initially move from the
origin to the site $z_1$ until time 1, then stay in $\bun{V}(s)$ time
interval of length $s-1$, at the end of which the particles return
to $z_1$ and the remaining time move freely. Assuming that $\xi(\cdot)$ is bounded from below (so
percolation effects of very low values of $\xi(\cdot)$ are
neglected), we obtain from the strong Markov property that
$\underline{U}(s)=v(s-1,z_1)\e^{\mo(s)}$ almost surely, where
$v(s,x)$ is the expectation over trajectories in $\bun{V}(s)$ starting
from $x \in \bun{V}(s)$ and ending at $z_1$ during the time interval $[0;s]$,
viz.
$$
v(s,x)=\textbf{E}_x \lf [\exp \lf \{\int^s_0\xi(X(a)) \dd a
\rg \}\vien_{\{X[0;s] \subset \bun{V}(s)\}} \delta\iz{z_1}(X(s))\rg ].
$$
Thus, $v(s,\cdot)$ is the solution to equation (7.1) in $\bun{V}(s)$
with the initial datum $\delta\iz{z_1}$ instead of $\delta\iz{0}$.
Using the spectral representation for $v(s-1,z_1)$ (where all
terms are nonnegative!) combined with the previous
estimates, we finally obtain the lower bound for $U(s)$: As $s \to
\infty$, almost surely
\begin{eqnarray*}
&&U(s) \dly v(s-1,z_1)\e^{\mo(s)} \dly \e^{\rrl_{1,\bun{V}}s-2\nu\rrk s
+\mo(s)}\psi (z_1;
\rrl_{1, \bun{V}} )^2 \\
&&\dly  \e^{\rrl_{1,\bun{V}}s-2\nu\rrk s +\mo(s)}\ve
\bun{V}(s)\ve^{-1}=\e^{\rrl_{1,\bun{V}}s-2\nu\rrk s +\mo(s)},
\end{eqnarray*}
since $\psi (z_1;
\rrl_{1, \bun{V}} )^2 \dly  \ve
\bun{V}\ve^{-1}$ by the definition of the site $z_1\in \bun{V}$. From this formula and (7.6), we get that almost surely
$$
\rrl \iz{1,\bun{V}(s)}-2\nu\rrk +\mro(1) \mly s^{-1}\log U(s)\mly \rrl \iz{1,V(s)}-2\nu\rrk +\mro(1) \ \ \hbox{as} \
\ s\to \infty.
$$
Applying the almost sure asymptotic formulas for the principal
eigenvalue $\rrl \iz{1,V}$ in Corollary 6.6 ($\rrr=\infty$),
Theorem  6.19 ($0<\rrr< \infty$) and Theorem 6.14 ($\rrr=0$) of
the present paper, we can now derive the corresponding asymptotics
for the total mass $U(s)$:

\begin{thm}{\rm (G\"artner and
Molchanov 1998; Section 2.1)}{\bf .} Let ${\it essinf}\
\xi(0)>-\infty$. Assume that  $f:=Q^{\leftarrow}$ satisfies the
following RV conditions: there is a constant $0 \mly \rrr \mly
\infty$ such that $\e^f \in RV_{\rrr}$ and, additionally,
$f(a+\log a)-f(a)\to 0$ as $a \to \infty$. Then with
proba\-bi\-li\-ty~1 \beq \frac{\log U(s)}{s}=f(\nu \log s)-2\nu
\rrk(1-q(\rrr/\rrk))+\mro (1)\ \ \mbox{as}\ \ s\to \infty .
\eeq Here the nonrandom constants $q(\rrr)$ are specified in
Section 2.5; in particular, $q(\rrr)$ are strictly decreasing in
$\rrr$; $q(0)=1$ and $q(\infty)=0$.
\end{thm}

\medskip

Recall that the condition $\e^f \in RV_{\rrr}$ with  $0 < \rrr <
\infty$ (resp.,  $\rrr =\infty$ and $\rrr =0$) ensures the double
exponential upper tails (1.14) (resp., heavier and lighter upper
tails than the double exponential) of the potential distribution
$1-\e^{-Q}$; see Lemma~A.3 of Appendix A. The additional RV
condition of the theorem is to exclude heavy-tailed distributions
of potential. Therefore, the first term on the right-hand side of
(7.7) is equal (with the accuracy $\mro(1)$) to the largest values of
the potential in $V(s)$; see Theorem 4.4(iii). The second term
describes the shape of the potential in the neighborhood of its
maxima and is specified by deterministic variational principles;
see Section 2.5. From Sections 1.2, 2 and 6, we know that the principal 
eigenvalue $\rrl\iz{1,V(s)}$ of the operator $\cH_{V(s)}=
\rrk\Delta_{V(s)} +\xi_{V(s)}$ is approximated (as $s \to \infty$) 
by the local principal eigenvalue in the connected region 
$\ttA_{\rm opt}(s)\subset V(s)$ where the potential $\xi_{V(s)}$ 
possesses high values of a particular preferred shape. The 
logarithmic asymptotics of $U(s)$ is
therefore fully specified by these high values of the potential.

G\"artner and Molchanov (1998) derived also the second-order
expansion formulas for statistical moments of the total mass
$U(s)$ as $s\to \infty $. For the upper distributional tails of
$\xi(0)$ lighter than the double exponential, Biskup and K\"onig
(2001), van der Hofstad et al.~(2006) obtained more accurate expansion
formulas for statistical moments and almost sure behavior of
$U(s)$. The spatial correlation structure for the PAM was
investigated by G\"artner and den Hollander (1999). See also
(Molchanov and Zhang 2012) for the Anderson parabolic model with
$\Delta_{{\rm dif}}$ replaced by the fractional Laplacian 
$-(-\Delta_{{\rm dif}})^\theta$, $0 < \theta <1$, where the
potential has Weibull type tails. In these papers, refined
variational arguments were involved to obtain additional
information on intermittent islands of solutions $u(s,\cdot)$ and the
asymptotic structure of related $\xi_{V(s)}$-extremes (their size,
optimal shape, etc., as $s \to \infty$). See also (K\"onig 2016)
for a recent survey on the subject.

Van der Hofstad et al.~(2008) considered the case of i.i.d.~potentials
with heavy upper tails, i.e., polynomially decaying (Pareto)
distributions and Weibull distributions (1.5) with $\rra < 1$.
Thus, all positive exponential moments of the potential are infinite,
in contrary to the assumptions of Theorem 7.1. For such classes of
potentials, they proved extremal type limit theorems and almost
sure asymptotic bounds for the logarithm of the total mass $U(s)$.

\medskip

\subsection{\normalsize{\textsf {Asymptotic concentration formulas} }}

However, the rough asymptotic expansion formulas for the PAM (like
in Theorem~7.1) provide only appropriate information on the
geometric structure of intermittent islands of the solutions $u(s,
\cdot)$ to (7.1). Recall that the intermittent islands are formed
by those highly concentrated $u(s,\cdot)$-values which give the
main contribution to the total mass $U(s)=\sum_x u(s,x)$, and the
contribution from the complement of these islands is negligible as
$s \to \infty$. Recently, there has been a considerable attention
to the following mathematical problems regarding a geometric
characterization of intermittency effect:

1) description of the shape and location of intermittent islands;

2) description of the shape of potential values which generate
intermittent islands;

3) specification of the minimal number of these islands, etc.

Thus, taking into account (7.3), one needs to prove the following
concentration formula for the total mass $U(s)$: As $s \to
\infty$, \beq U(s)\sim U\iz{V(s)}(s) \sim \sum_{k=1}^{n(s)}
\sum_{x \in \ttA_{\rm opt}^k(s)} u(s, x) \eeq
in the sense of asymptotic equivalence almost surely or in
probability, where $\ttA_{\rm opt}^k(s)$ are believed to present
random connected regions (i.e.~intermittent islands) in $V(s)$ at a large distance
from each other, such that the diameter of each $\ttA_{\rm
opt}^k(s)$ is much smaller than this distance; $n(s)$ are relatively
small numbers, and $V(s) \uparrow \ttZ^\nu$ are centered cubes as
above.

Let us give a heuristic explanation of formula (7.8) in the
spectral terms of the operator $\cH\iz{V}=\rrk\Delta_V +\xi_V$ in
$V=V(s)$, provided the conditions of Theorem 7.1 are fulfilled.
To this end, we need more careful inspection of spectral
representation formulas (7.4)--(7.5), by applying extreme value
theory for eigenvalues including Poisson limit theorems and
localization properties discussed in Sections 1.2, 2, and 6. First,
notice that the exponents of the top eigenvalues are essentially
larger than those associated with lower eigenvalues. Therefore, it suffices to
consider the sum of a few first terms of (7.4)--(7.5) associated
with the largest eigenvalues $\rrl\iz{k,V}$, $1\mly k \mly n:=
n(s)$; the other terms in (7.4)--(7.5) associated with the lower
eigenvalues are asymptotically negligible. Thus, the random field $u(s,\cdot)$
($s\to \infty$) may be interpreted as a superposition of a few
wave functions $\psi (\cdotp ; \rrl\iz{k, V})$ for $1\mly k \mly
n$. The general theory of Anderson localization suggests that the
$k$th eigenfunction $\psi (\cdotp ; \rrl\iz{k, V})$ is
exponentially localized at the $s$-dependent center $z\iz{\tau(k),
V}\in V$, and it is highly concentrated in a $s$-dependent
neighborhood $\ttA_{\rm opt}^k\subset V$ of the site
$z\iz{\tau(k), V}$. The regions $\ttA_{\rm opt}^k\subset V$
possess a relatively small size and are asymptotically far from
each other according to Poissonian behavior of the localization
centers $z\iz{\tau(k), V}\ve V\ve^{-1/\nu}$; $k\dly 1$. Moreover,
$\xi(\cdot)$ possesses in $\ttA_{\rm opt}^k$ the deterministic
optimal shape specified by the variational principles; cf.~Section
2.5. From these observations when applied to (7.4) as $s\to
\infty$, we see that the function
$$
\e^{\rrl\iz{k, V}s-2\nu \rrk s }\psi (0; \rrl\iz{k, V} ) \psi
(\cdot; \rrl\iz{k, V} )
$$
is a very good approximation for $u_V(s,\cdot)$ in the region
$\ttA_{\rm opt}^k$ for each $1\mly k \mly n$. This in turn
suggests that the mass concentration formula (7.8) holds true,
where $V=V(s)\uparrow \ttZ^{\nu}$ as above, $n=n(s)$ are
relatively small numbers, and $\ttA_{\rm opt}^k \subset V$ are the
relevant regions defined above. Thus, the intermittent islands
$\ttA_{\rm opt}^k(s)$ have relatively small size 
and are far away from each other as $s
\to \infty$. This concludes the heuristic explanation of formula (7.8).

G\"artner et al.~(2007) proved the concentration formula (7.8) 
almost surely, for potential
distributions satisfying the conditions of Theorem 7.1 with $0<
\rrr \mly \infty$, i.e., the double exponential upper tails and
heavier than the double exponential. They showed that almost surely
$n(s)=s^{\mo(1)}$, and the connected regions $\ttA_{\rm
opt}^k(s)\subset V(s)$ are asymptotically bounded 
only when defined properly (see above), 
and the distance between them is of order $s^{1-\mo(1)}$. 
For $\rrr=\infty$, the regions $\ttA_{\rm
opt}^k(s)$ shrink to
singletons. Moreover, the shape of the solutions $u(s,\cdot)$ and 
potential values in $\ttA_{\rm opt}^k(s)$ are specified
(via the variational formulas) by the local principal
eigenfunction and the principal eigenvalue in $\ttA_{\rm opt}^k(s)$.
This agrees with the heuristics given above; cf.~also Sections~1.2,
2, and 6 treating the extreme value theory for eigenvalues of the
Hamiltonian $\cH_V$. However, the proof of the asymptotic
concentration becomes complicated by applying straightforwardly the
asymptotic results for the top spectrum of $\cH_V$, because of the
possibly different signs of the $k$th eigenfunctions in
(7.4)--(7.5) where the factor $\psi (0; \rrl\iz{k, V} )$
should be also taken into account. Instead, the authors explore the
Feynman-Kac formulas for $u(s,\cdot)$, $U(s)$ as well as for the
principal eigenfunctions with a slightly modified potential. To
prove the exponential decay of the principal eigenfunctions, they
apply a decomposition technique for the trajectories of the random
walk in the corresponding Feynman-Kac representations.

Sznitman (1998) earlier proved similar mass concentration results
for a Brownian motion in $\ttR^{\nu}$ among Poisson obstacles;
here the potential $\xi(\cdot)$ is given by formula (6.20). In
particular, the spatial regions $\ttA_{\rm opt}^k(s)\subset
\ttR^{\nu}$ in (7.8) were shown to have no obstacles and
unboundedly increase almost surely as $s\to \infty$. The optimal
strategy of the Brownian particle during the time period $[0;s]$ is to
move quickly to one of the obstacle-free regions $\ttA_{\rm
opt}^k(s)$ of the optimal shape, i.e.~the ball of radius $\const (\log
s)^{1/\nu}$, and to stay here for the remaining time. Notice that
the intermittent behavior of this model is rather similar to that of the
spatially discrete PAM with the potential bounded from above. The
related asymptotic results for the principal Dirichlet eigenvalues
are discussed in Example 6.18 of the present survey.

However, in (G\"artner et al.~2007) and (Sznitman 1998), the
problem of the minimal number of intermittent islands was not
considered. This problem in precise setting was solved by several
mathematicians for Pareto distributions \beq \bP(\xi(0)>
t)=\e^{-Q(t)}=t^{-\beta} \quad (t\dly1)
\eeq with $\beta >\nu$, and Weibull distributions (1.5) with
arbitrary $\rra > 0$. Recall that the choice of $\beta >\nu$ in (7.9)
is to guarantee the existence and uniqueness of the solution
$u(s,\cdot)$ to equation (7.1). Write, as above,
$U(s):=\sum_{x}u(s,x)$ for the total mass.

\begin{thm}{\rm (e.g., K\"onig et al.~2009; Sidorova and Twarowski
2014; Fiodorov and Muirhead 2014)}{\bf .} Assume that the
potential has either Pareto distribution with $\beta >\nu$, or
Weibull distribution {\rm (1.5)} with arbitrary $\rra > 0$ . Then
there exists a random process $Z_{{\rm opt}}(s)$ ($s \dly 0$) with
values in $\ttZ^\nu$ such that \beq \lim_{s\to\infty}\frac{u \lf (
s,Z_{{\rm opt}}(s)\rg )}{U(s)}=1\ \ \mbox{in \ probability.}
\eeq
\end{thm}

\medskip

This theorem states the complete localization property for
$u(s,\cdot)$ as $s\to \infty$, which is the strongest case of mass
concentration formula (7.8) in probability  with $n(s)\equiv 1$
and $\ttA_{\rm opt}^1(s)\equiv\{Z_{{\rm opt}}(s)\}$, a singleton.

For the Weibull distribution with $0<\rra <2$, asymptotic formula
(7.10) was proved by Sidorova and Twarowski (2014). This result
was extended by Fiodorov and Muirhead (2014) to an arbitrary
$\rra>0$. See also (Muirhead  and Pymar 2014) for the proof of the
single-site concentration property for a random walk in a random
environment (instead of the standard random walk like in the PAM)
for the Weibull-distributed i.i.d.~potential.

Earlier, K\"onig et al.~(2009) proved (7.10) for the
Pareto-distributed potential. They also established a two-site
concentration property for the PAM: Almost surely $U(s)\sim u(s,
Z^{(1)}(s))+u(s, Z^{(2)}(s))$ as $s\to \infty$, where $Z^{(k)}(s)$
are distinct random processes with values in $\ttZ^\nu$. I.e., one
obtains (7.8)  with probability one, where $n(s)\equiv~2$,
$\ttA_{\rm opt}^1(s)$ and $\ttA_{\rm opt}^2(s)$ are two singletons.
This is the strongest almost sure version of the localization property, since for i.i.d.~random
potentials, the random field $u(s,\cdot)$ is asymptotically
concentrated on at least two distinct sites in $\ttZ^\nu$. The
reason of this fact lays on the observation that the localization
sites are changing infinitely often as $s \to \infty$; thus, at
very rare time moments when particles move from the previous
localization site to a new one, the mass of particles should be
concentrated on at least two different sites.

We recall from Sections 3--5 of the present survey that, for
potential distributions as in Theorem 7.2, the $\xi_V$-peaks are
strongly pronounced as $V\uparrow\ttZ^\nu$; therefore, the
single-site concentration formula, no surprise, holds true
according to the general picture of intermittency based on the
Feynman-Kac formulas. On the other hand, for such potential
distributions, there is a very precise extreme value theory for
eigenvalues of the operators $\cH_V=\rrk\Delta_V +\xi_V$, which
can (and does) provide powerful techniques for the investigation
of single-site concentration properties; cf.~Sections 1.2,
2.2--2.3 and 6.1. See also (K\"onig 2016) for a recent survey on
the subject.

\medskip

Let us {\it sketch the proof of Theorem 7.2 for the
Weibull-distributed potentials} by applying the extreme value theory
for eigenvalues. We follow the terminology and ideas of (Fiodorov
and Muirhead 2014). First, one obtains finite-volume approximation
formulas (7.3) with \lk good\rk accuracy, where $V=V(s)$ are
centered cubes in $\ttZ^\nu$ with the size length of order $s(\log
s)^{1/\rra}$. Then, let us look at the spectral representation
formulas (7.4) and (7.5) for $u_V(s,\cdot)$ and $U_V(s)$.
Because of the factor $\big (\psi (\cdotp ;\rrl\iz{k,
V}),u(0,\cdot)\big )_{V}= \psi (0; \rrl\iz{k, V} )\dly 0$ ($1\mly
k \mly \ve V\ve$) in (7.4), we need to study {\it the penalised
spectrum}
$$
\Psi(s;k):=\rrl_{k,V}+s^{-1} \log \psi(0;\rrl\iz{k,V})-2\nu\rrk
\quad (1\mly k \mly \ve V\ve)
$$
instead of the usual $\Spect (\cH_V)=\lf \{ \rrl\iz{k,V} \colon\ 1
\mly k \mly \ve V\ve \rg \}$. Let us rewrite $u_V$ (7.4) in the
terms of the penalised spectrum $\Psi(s):=\lf \{ \Psi(s; k)
\colon\ 1 \mly k \mly \ve V\ve \rg \}$: \beq u\iz{V}(s, x)=\sum
_{k=1}^{\vert V\vert }\e^{s \Psi(s;k)}\psi (x; \rrl\iz{k, V}
)\quad (x\in V).
\eeq For $1\mly l \mly \ve V\ve$, denote by $\Psi \iz{l,V}(s)$ the
$l$th largest value in the sample $\Psi(s)$. It will turn out that
the gap $\Psi \iz{1,V}(s)-\Psi \iz{2,V}(s)$ between the the first
largest $\Psi \iz{1,V}(s)$  and the second largest $\Psi
\iz{2,V}(s)$ in $\Psi(s)$ is sufficiently large; and moreover, the
eigenfunctions $ \psi(\cdot;\rrl\iz{k,V})$ decay exponentially (so that $ \psi(0;\rrl\iz{k,V})>0$),
for each $k\mly \ve V \ve^\rre$ with $\rre > 0$ small enough; cf.
also Theorem 6.2, Corollary 6.10 and Example 6.12 of the present
paper. We shall prove that the right-hand side of (7.4)--(7.5)
is dominated by just one term associated with $\Psi \iz{1,V}(s)$;
therefore, the localization center of the corresponding
eigenfunction should be the concentration site for the random
field $u_V(s,\cdot)$ as $s \to \infty$. To be more precise, let $\rrl_{{\rm
opt}}:=\rrl\iz{k^\ast,V}$ and $\psi_{{\rm opt}}(\cdot)$ denote the
eigenvalue and eigenfunction  of the operator $\cH_V$ associated
with the first largest value $\Psi \iz{1,V}(s)$ among the
penalised spectrum $\Psi(s)$, i.e.
$$
\Psi_{1,V}(s)=\rrl_{k^\ast,V}+s^{-1} \log \psi_{{\rm opt}}(0)-2\nu\rrk;
$$
here $k^\ast < \ve V\ve^\rre$ with
$\rre > 0$ small enough. Let $Z_{{\rm opt}}(s)\in V$ stand for the
localization center of the eigenfunction $\psi_{{\rm
opt}}(\cdot)$. Also, for random processes $Y(s)$ and $W(s)$, we
write $Y(s)\approx W(s)$ as $s\to \infty$, if the difference
$Y(s)-W(s)$ tends to zero sufficiently fast in probability. Using
this abbreviation and applying Cauchy-Schwarz inequality for
eigenfunctions in (7.11), we obtain that
$$
\max_x \lf \ve \frac {u_V(s,x)} {\e^{s\Psi \iz{1,V}(s)}}-\psi_{\rm
opt}(x) \rg \ve \mly \ve V\ve \exp \lf \{-s \lf (\Psi
_{1,V}(s)-\Psi _{2,V}(s)\rg)\rg\} \approx 0,
$$
therefore, for the total mass $U_V(s)$ we have that
$$
\lf \ve \frac {U_V(s)} {\e^{s\Psi \iz{1,V}(s)}}-\lf (\psi_{\rm opt},1\rg
)_{V} \rg \ve \mly \ve V\ve^2 \exp \lf \{-s \lf (\Psi
_{1,V}(s)-\Psi _{2,V}(s)\rg)\rg\} \approx 0,
$$
as $s\to \infty$. Consequently,
$$
\frac {u_V(s,Z_{{\rm opt}}(s))} {U_V(s)}\approx \frac {\psi_{\rm
opt}(Z_{{\rm opt}}(s))} {\lf (\psi_{\rm opt},1\rg )_{V}}\approx 1
$$
as $s\to \infty$, because of the sharp exponential decay of
$\psi_{{\rm opt}}(\cdot)$. Since $u(s,\cdot)\approx u_V(s,\cdot)$
in $V$ and $U(s)\approx U_V(s)$, the last formula implies (7.10),
as claimed.

However, the penalised spectral values are too complicated to
handle. In order to study the spacings of the largest values in
$\Psi(s)$ as well as further properties of the concentration site
$Z_{{\rm opt}}(s)$, one needs a good approximation for $\Psi_{k,V}(s)$
by a simpler function on potential configurations. To this end, we introduce the following auxiliary quantities: Write
$J:= [(\rra-1)/2]$ (= the integer part) for $\rra \dly 1$, and
$J=0$ otherwise. Given $z\in V$, let $\rrl^{(J)}(z)$ denote the
principal eigenvalue of the single-peak Hamiltonian
$$
\rrk \Delta\iz{V} +\sum_{y\colon 1 \mly |y-z| \mly J} \wt
\xi(y)\delta _y+\xi(z)\delta_z \quad \hbox{on}\quad  l^2(V).
$$
I.e., $\rrl^{(J)}(z)$ is the local principal eigenvalue on the
lattice ball of radius $J$ with a single $\xi_V$-peak at $z$
surrounded by the island of lower $\xi_V$-values. Let
$\rrl^{(J)}_{k,V}:=\rrl^{(J)}(z_{\tau(k),V})$ denote the $k$th largest value of the sample $\{ \rrl^{(J)}(z)\colon\ z\in
V\}$. We now observe from Section~6.1 of the present survey that
the eigenvalue $\rrl_{k,V}$ is approximately equal to the local
eigenvalue $\rrl^{(J)}_{k,V}$, i.e., $\rrl_{k,V}\approx
\rrl^{(J)}_{k,V}$ as $V=V(s)\uparrow \ttZ^\nu$  and $k < \ve
V\ve^\rre$ with $\rre
> 0$ small enough. This observation and more careful inspection of
the exponential decay of the $k$th eigenfunction
$\psi(\cdot;\rrl\iz{k,V})$ at the origin (cf.~Astrauskas 2008;
Section 6) suggest that, with high probability, the $K$th largest
value $\Psi \iz{K,V}(s)$ of the penalised spectrum $\Psi(s)$ is
approximately equal to the $K$th largest value $\Upsilon^{(J)}
_{K,V}(s)$ of the {\it penalisation functional} \beq
\Upsilon^{(J)} (s,z):=\rrl^{(J)}(z)-\frac{\ve z\ve}{s}
\cdot\frac{\log\log s}{\rra} -2\nu\rrk \quad (z \in V).
\eeq Comparing the penalised spectral values to (7.12), we observe, in addition, that the quantity $\rra^{-1}\log\log s =\log b_{V(s)}+\rO(1)$
is the nonrandom rate of the exponential decay of the $k$th
eigenfunction $\psi(\cdot;\rrl\iz{k,V})$ described in Proposition
1.3; moreover, $\log \psi(0;\rrl\iz{k,V})$ asymptotically behaves
like $-|z\iz{\tau(k),V}|\log b_V$ as $V\uparrow \ttZ^\nu$. Recall also
that, for each $z\in V$, the eigenvalue $\rrl^{(J)}(z)$ is a
certain (nonlinear) function of the sample $\lf \{ \xi(z+x)
\colon\ 0 \mly \ve x\ve \mly J \rg \}$; cf.~(2.20)--(2.22). Now
the concentration site $Z_{{\rm opt}}(s)$ can be defined as the
maximizer of the random field $\Upsilon^{(J)} (s,\cdot)$ in $V$;
here $\Upsilon^{(J)} (s,z)$ is a certain function of both the sample
$\lf \{ \xi(z+x) \colon\ 0 \mly \ve x\ve \mly J \rg \}$ and the
site $z$. Similarly as in Theorem 6.2 above, one obtains Poisson
limit theorems for the (normalized) penalisation functionals
$\Upsilon^{(J)} (s,\cdot)$ and their locations. This limit theorem
implies the limiting distributions for the normalized concentration site
$Z_{{\rm opt}}(s)$ as well as for the spacings $\Upsilon^{(J)}
_{K,V}(s)-\Upsilon^{(J)} _{K+1,V}(s)$ of the random field
$\Upsilon^{(J)} (s,\cdot)$ in $V$. The latter in turn implies the
existence of a sufficiently large gap between the largest values in
the penalised spectrum $\Psi(s)$, as claimed. This assertion concludes the heuristic proof of the single-site concentration
property (7.10) of the solution $u(s,\cdot)$. \qed

\medskip

As stated in (Fiodorov and Muirhead 2014), the above
considerations are a starting point in obtaining more information
on the asymptotic behavior of the concentration site $Z_{{\rm
opt}}(s)$ and the shape of potential in its neighborhood as $s\to
\infty$, provided the random potential has Weibull distributions.
In particular, Poisson limit theorems for the penalisation
functional (7.12) imply that the site $Z_{{\rm opt}}(s)\in V(s)$
is of order $s$ up to logarithmic corrections and has the limiting
distribution as a product of univariate Laplace distributions.
Recall also that the random field $\Upsilon^{(J)} (s,\cdot)$ has
the finite range (=$J$) of dependency; see (Fiodorov and Muirhead
2014) where $J$ is called as the {\it radius of influence}. For
sufficiently heavy tails (Weibull distributions with $\rra <2$),
the eigenvalue $\rrl^{(J)}(z)$ in (7.12) may be replaced by
$\xi(z)$, so that $\Upsilon^{(J)} (s,z)$ ($z\in V$) are
independent non-identically distributed random variables. (Recall
that the case $\rra <2$ was studied in (Sidorova and Twarowski
2014) by exploring the Feynman-Kac representations). These
observations are crucial for describing the shape of $\xi(\cdot)$
in the neighborhood of the concentration site $Z_{{\rm opt}}(s)$:
As $s \to \infty$, with probability $1+\mro(1)$ the single peak
$\xi \lf (Z_{{\rm opt}}(s)\rg )$ is extremely high:
$$
\xi \lf (Z_{{\rm opt}}(s)\rg )=b\iz{V(s)}(1+\mro(1)),
$$
where $b\iz{V}:=(\log \ve V\ve)^{1/\rra}$; cf.~Example 6.12 above.
Meanwhile, the neighboring values $\xi(x)$ ($1 \mly \lf \ve
x-Z_{{\rm opt}}(s)\rg \ve \mly J$) are essentially lower: there
exists a strictly decreasing nonrandom function $d(\cdot)\colon
[1;J]\mapsto [0; 1)$ such that
$$
\xi \lf (x\rg )\asymp b\iz{V(s)}^{d(\lf \ve x-Z_{{\rm opt}}(s)\rg
\ve)} \quad \hbox{for} \ 1 \mly \lf \ve
x-Z_{{\rm opt}}(s)\rg \ve \mly J.
$$

This characterization of the concentration site $Z_{{\rm opt}}(s)$
agrees with the asymptotic results for the largest eigenvalues
$\rrl\iz{K, V}$ and eigenfunctions $\psi (\cdotp ; \rrl\iz{K, V})$
of the Hamiltonian $\cH\iz{V}=\rrk\Delta_V +\xi_V$ as $V \uparrow
\ttZ^\nu$ and $K \dly 1$ fixed, provided $\xi(0)$ has Weibull
distribution; cf.~Sections 1.2, 2.2--2.3, and 6.1 of the present
paper. In particular, let us look at the asymptotic expansion
formulas for $\rrl\iz{K, V}$ (Example 6.12) to observe that the
leading term $b_V$ comes from an isolated high peak
$\xi(z\iz{\tau(K), V})$; meanwhile, the further terms of order
$\mro(1)$ come from the neighboring values $\xi(x)$ ($1 \mly \lf
\ve x-z\iz{\tau(K), V}\rg \ve \mly J$) with the same influence
radius $J$ as in the PAM above; cf.~also (Astrauskas 2008; Section
6). Moreover, the eigenfunction $\psi (\cdotp ; \rrl\iz{K, V})$ is
highly concentrated at the site $z\iz{\tau(K), V}$, the
localization center (cf.~Theorem 6.2(II)); and the
neighboring values $\xi(x)$ ($1 \mly \lf \ve x-z\iz{\tau(K), V}\rg
\ve \mly J$) have significant influence on the asymptotic behavior
of the localization index $\tau(K)=\tau_V(K)$ and on the concentration
degree of $\psi (\cdotp ; \rrl\iz{K, V})$ in the neighborhood of
$z\iz{\tau(K), V}$; see (Astrauskas 2013). These observations
suggest the following conclusion:  the lighter are the tails of
potential, the larger is the influence radius in both models; thus, the
weaker are the localization properties of both models, PAM
and (time-independent) Anderson Hamiltonian.

It is worth mentioning, at the heuristic level of rigor, that
Section 2 (resp., Section 6) of the present paper exhibits all
classes of \lk typical\rk configurations of the potential (resp., RV
classes of potential distributions) which are thought to guarantee
certain concentration properties for the solutions $u$ to equation
(7.1). For instance, the results of Sections 2.3 and 6.1 (i.e.,
\lk relevant single-peak\rk approximation) should be preliminaries
to establish a single-site concentration property in probability
for solutions $u$, provided the upper tails of potential are
heavier than the double exponential. In particular, the results of
Section 2.2 and Theorem 6.9 (\lk sharp single-peak \rk
approximation) are related to the single-site concentration
property for $u$ with zero influence radius $J=0$, provided the
distributional tails are heavier than Weibull's tails with
$\rra=3$. In view of the results of Sections 2.5 and 6.3 (\lk
relevant island\rk approximation in the double exponential case),
it can be conjectured that, with high probability, the solution
$u(s,\cdot)$ exhibits entire concentration on a single island,
the diameter of which is asymptotically bounded 
as $s\to \infty$. See (K\"onig 2016) for a heuristic explanation
of this conjecture.

If the upper tails of potential are lighter than those of double
exponential including bounded tails ($t_Q < \infty$), no results
on concentration properties for the PAM are known, with the exception of the
very special (but important) spatially continuous model of
Brownian motion in a Poisson field of obstacles studied by
Sznitman (1998).

\medskip

\smallskip

\newenvironment{lemm}{\smallskip {\textbf{\textsf{Lemma}}}\
{\bf {\textsf{\Alph{section}.\arabic{prop}.}}}\ \it }{\smallskip
\addtocounter{prop}{1}}
 \setcounter{equation}{0}

 \setcounter{equation}{0} \addtocounter{sct}{1}

\renewcommand{\thesection}{\Alph{section}}
\renewcommand{\thesubsection}{\normalsize{\textsf
{\Alph{section}.\arabic{subsection}.}}}

\newcounter{lemm}
\newcounter{skaitl}

\newenvironment{thmm}{\smallskip {\textbf{\textsf{Theorem}}}\
{\bf {\textsf{\Alph{section}.\arabic{prop}.}}}\ \it }{\smallskip
\addtocounter{prop}{1}}

\renewcommand{\theequation}{\Alph{section}.\arabic{equation}}

\setcounter{lemm}{0}

 \setcounter{section}{0}

\setcounter{equation}{0} \addtocounter{sct}{1}

\setcounter{prop}{1}




\newenvironment{exm}{\smallskip {\bf {\textsf{Example}}}\
{\bf {\textsf{\Alph{section}.\arabic{prop}.}}}\ \rm }{\smallskip \addtocounter{prop}{1}}


\appendix
\section{\normalsize{\textsf {APPENDIX A: REGULAR VARIATION}}}

In this section, we study the classes of functions
$f:=Q^\leftarrow$ (the left-continuous inverse of the cumulative
hazard function) introduced in Sections 4--6. These classes are
characterized in terms of $Q$. The tail behavior of $Q(t)$ as
$t\uparrow t_Q$ is also treated. In Section A.1, we recall the
classical results on the domain of attraction of max-stable Gumbel 
law and regular variation $RV_\rho$. The classes $A\Pi^p_\infty$
(4.1), $A\Pi^p_0$ (4.2) and $\rO A\Pi^p$ (4.3) are studied in
Section~A.2, and $PI_{<2}$ (4.5) in Section A.3. Finally, examples
and counterexamples are given in Section~A.4.




\subsection{\normalsize \textsf {The domain of attraction of max-stable Gumbel 
distribution and regular variation} }

We now give the well-known characterization statements for the
distribution function to be in the domain of attraction of the
max-stable Gumbel law $G_{\exp}(\cdot)$.



\begin{lemm} {\rm (Resnick 1987; de Haan and Ferreira 2006).} The
following assertions {\rm (i)}--{\rm (iii)} are equivalent:

{\rm (i)} $f\in A\Pi$ {\rm (6.3)} with an auxiliary function
$a(\cdotp)>0 $;

{\rm (ii)} there exists another auxiliary function $a_1:(-\infty;
t_Q) \rightarrow \ttR_+$ such that
$$
Q(t+ca_1(t))-Q(t) \rightarrow c \ \ \mbox{as} \ \ t\uparrow t_Q, \
\ \mbox{for any} \ \ c\in \ttR;
$$

{\rm (iii)} there exist functions $b:(-\infty; t_Q) \rightarrow
\ttR $ and $a_2:(-\infty; t_Q) \rightarrow \ttR_+$ such that
$$
Q(t)=b(t)+\int^t_{t_0}1/a_2(s)\dd s \quad \ (t<t_Q),
$$
here $b(t) \to \overline{b} \in \ttR$ ($t \uparrow t_Q$), the
function $a_2$ is locally absolutely continuous with the density
$a^{\prime}_2(t) \to 0$ ($t\uparrow t_Q$) and, for $t_Q < \infty$,
$a_2(t) \to 0$ ($t\uparrow t_Q$).

In this case, $a \circ f(s)=a_1 \circ f(s)(1+\mro(1))= a_2 \circ
f(s)(1+\mro(1))$ as $s \to \infty$. Moreover, the limit in (ii)
with an auxiliary function $a_1(\cdotp)> 0$ implies that $f\in
A\Pi$ {\rm (6.3)} with the same auxiliary function
$a(\cdotp)\equiv a_1(\cdotp)$.

\end{lemm}

\medskip

\begin{exm} For $t_Q=\infty$, $p \dly 0$ and $B>0$,
consider the subclass $A\Pi^p_B \subset A\Pi$ associated with the
auxiliary function $a_2(s):=B(p+1)^{-1}s^{-p}$ in Lemma A.1(iii).
In this case, $Q(t)=B^{-1}t^{p+1}+\const +\mro(1)$, i.e.,
$1-\e^{-Q}$ are Weibull type distributions. In the next section,
we extend the subclass $A\Pi^p_B$ to the boundary cases
$B=\infty$, $B=0$ and $\rO$-type asymptotics.
\end{exm}

\medskip

Let us discuss the class $RV_\rho$ of (nondecreasing) functions,
which are regularly varying at infinity with index $\rho$. Recall
that, for $0 < \beta < \infty$, the condition $f\circ \log \in
RV_{1/\beta}$ is sufficient and necessary for the distribution
$1-\e^{-Q}$ to be in the domain of attraction of max-stable
Fr\'{e}chet law $G_\beta(t):=\exp \lf \{ -t^{-\beta} \rg \}$
($t>0$); cf.~Example 6.11. We now explore the class $RV_\rho$ to
characterize the double exponential type distributions.

\begin{lemm} For $t_Q= \infty$ and $0\mly \rho \mly \infty$, the following assertions are
equivalent:

{\rm (i)} $\e^f \in RV_\rho$;

{\rm (ii)} $f(s)-f(\delta s) \to -\rho \log \delta$ as $s \to
\infty$, for any $0< \delta <1$;

{\rm (iii)} $Q(t+C)/Q(t)\to \e^{C/\rho}$ as $t\to \infty$, for any
$C>0$.

For any $0\mly \rho \mly \infty$, either of {\rm (i)--(iii)} implies that $\lim_{t\to \infty} t^{-1}
\log Q(t)=\rho^{-1}$. Finally, for any $0< \rho \mly \infty$, either of {\rm (i)--(iii)} yields that $Q(t-)/Q(t) \to 1$ as $t\to \infty$, i.e., the continuity condition {\rm (3.1)}.
\end{lemm}

\begin{proof}
The equivalence of (i)--(iii) follows from Theorems 1.5.12, 2.4.7
and Propositions 2.4.4(iv) and 1.3.6(i) in (Bingham et al.~1987)
combined with the observation in (Resnick 1987, Sect.~0.2) that
$Q(t-\rre)\mly f^\leftarrow(t) \mly Q(t)$ for all $t \in \ttR$ and
all $\rre > 0$. For $\rrr=\infty$ and $\rrr=0$, the equivalence of
(i)--(iii) is also proved, respectively, in Lemma~A.6 ($p=0$) and
Lemma~A.7 ($p=0$) of the present paper adapted for the
argument-additive functions; see also (de Haan and Ferreira 2006, 
Chapter 1 and Appendix B.1) for the case $0 < \rrr < \infty$.\qed
\end{proof}

\medskip

The following lemma is compounded of Lemma A.1 and Lemma A.3 with
$\rrr =\infty$, provided that there exists the density of the
distribution $1-\e^{-Q}$.

\begin{lemm}{\rm (Cf.~Corollary 6.4)}. Let $t_Q= \infty$. For some large $t_0$, assume
that $Q:[t_0; \infty) \rightarrow \ttR_+$ is (locally) absolutely
continuous with the positive density $Q^{\prime}:[t_0; \infty)
\rightarrow \ttR_+$ obeying the following conditions: \beq
\lim_{t\to \infty }\frac {Q^{\prime}(t+C)} {Q^{\prime}(t)}= 1\quad
\mbox{for\ any} \quad C>0,
\eeq and \beq \linf_{t\to \infty }Q^{\prime}(t)>0.
\eeq

Then the following limits {\rm (I)}--{\rm (IV)} hold true:

\vspace{6pt} {\rm (I)} $\lim_{t\to \infty } Q(t+u)/Q(t)=
\lim_{t\to \infty }Q^{\prime}(t+u)/Q^{\prime}(t)= 1$ uniformly in
compact sets of $u \in \ttR$;

\vspace{6pt} {\rm (II)} $\linf_{t\to \infty }Q(t)/t>0$;

\vspace{6pt} {\rm (III)} $\lim_{t\to \infty }\big
(Q(t+va_1(t))-Q(t)\big )=v$ uniformly in compact sets of $v \in
\ttR$, with $a_1(\cdot)\equiv 1/Q^{\prime}(\cdot)$ in $[t_0;
\infty)$;

\vspace{6pt} {\rm (IV)} with $p(t):=\e^{-Q(t)}Q^{\prime}(t)$ ($t
\dly t_0$) as the distribution density and $a_1(\cdot)$ as in part
{\rm (III)},
$$
\lim_{t\to \infty }\frac {p(t+u+va_1(t))} {p(t+u)}=\e^{-v}
$$
uniformly in compact sets of $v, u \in \ttR$.
\end{lemm}

\begin{proof}
(I) By L'H\^opital's rule, we obtain the first limit for any $u
\in \ttR$. The uniform convergence follows from Theorem 1.2.1 in
(Bingham et al.~1987) adapted for the argument-additive functions.

(II) The assertion follows from (A.2).

(III) Writing
$$
Q(t+va_1(t))-Q(t)-v=v\int^1_0 \Big ( a_1(t)Q^{\prime}(t+\theta
va_1(t))-1\Big)\dd \theta \quad \mbox{for} \quad t \dly t_0
$$
and applying assertion (I) and condition (A.2), we easily obtain
the claimed limit.

(IV) Let us rewrite the ratio under the limit in the form:
\begin{eqnarray}
&&\frac {p(t+u+va_1(t))} {p(t+u)}=\exp \Big\{ -\big
(Q(t+u+va_1(t))-Q(t+u)\big) \Big \}\times \n \\
&&\hspace{2.9cm}\frac {Q^{\prime}(t+u+va_1(t))}
{Q^{\prime}(t)}\bigg (\frac {Q^{\prime}(t+u)} {Q^{\prime}(t)}\bigg
)^{-1};
\end{eqnarray}
$t \dly t_0$. Since $a_1(\cdot)$ is a bounded function, assertion
(I) implies that the last two ratios on the right-hand side of
(A.3) converge to 1 locally uniformly in $u,v \in \ttR$. It
remains to prove the uniform convergence of the exponent on the
right-hand side of (A.3). By the theorem of continuous convergence
(see, e.g., p. 2 in (Resnick 1987)), it suffices to check that,
for arbitrary functions $u(t)\rightarrow u$ and $v(t)\rightarrow
v$, the following limit holds true:
$$
Q\big(t+u(t)+v(t)a_1(t)\big)-Q\big(t+u(t)\big) -v(t)\rightarrow 0\
\ \mbox{as} \ \ t\uparrow \infty.
$$
This is shown similarly as in part (III), so we omit the details.
Lemma A.4 is proved.~\qed
\end{proof}

\medskip

\newenvironment{remm}{\smallskip {\bf {\textsf{Remark}}}\
{\bf {\textsf{\Alph{section}.\arabic{prop}.}}}\ \rm }{\smallskip
\addtocounter{prop}{1}}

\begin{remm}(Cf.~Theorem 6.20).
Let $t_Q= \infty$. For some $t_0$, assume that $Q:[t_0; \infty)
\rightarrow \ttR_+$ is (locally) absolutely continuous with the
positive density $Q^{\prime}$ satisfying  the following condition:
$(\log Q)^{\prime}(t) \to \rrr^{-1}$ as $t \to \infty$, for some
$0 < \rrr < \infty$. Then the following limits hold true:

(i) $\lim_{t\to \infty } Q(t+C)/Q(t)= \lim_{t\to \infty
}Q^{\prime}(t+C)/Q^{\prime}(t)= \e^{C/\rho}$ uniformly in compact
sets of $C \in \ttR$;

(ii) with $a(\cdot)\equiv 1/Q^{\prime}(\cdot)$ in $[t_0; \infty)$,
$$
\lim_{t\to \infty }\big (Q(t+Ca(t))-Q(t)\big )=C \quad \mbox{for \
any} \quad  C \in \ttR.
$$.
\end{remm}

\medskip

{\it The proof of the assertions of Remark A.5}. (i) Write $\log
Q$ in the form:
$$
\log Q(t)=\const +\frac{t}{\rrr} +\int^t_{t_0}\rre(s)\dd s \ \ (t
\dly t_0),
$$
where $\rre(t):=(\log Q)^{\prime}(t) -\rrr^{-1} \to 0$ as $t \to
\infty$. Using this representation and the conditions of Remark
A.5, we obtain the claimed limits for any $C \in \ttR$. The
uniform convergence follows from Theorem 1.5.2 in (Bingham et al.~1987) adapted for the argument-additive functions.

(ii) Since $a(t)=\mro(1)$, the claimed limit is derived similarly
as in the proof of Lemma~A.4(III). \qed

\medskip

\subsection{\normalsize \textsf {Classes $A\Pi^p_\infty$, $A\Pi^p_0$ and $\rO A\Pi^p$} }

Recall that, for $p \dly 0$, the classes $A\Pi^{p}_{\infty}$,
$A\Pi^{p}_{0}$ and $\rO A\Pi^p$ consist of functions $f:=Q^{\la}$
satisfying, respectively, $f(s)^p\big ( f(s+c)-f(s) \big )\to
\infty $, $\to 0$ and $\asymp 1$ as $s \to\infty$, for any $c>0$;
cf.~(4.1)--(4.3). We first formulate the results of this section.
To avoid trivialities, we restrict ourselves to the case $t_Q =
\infty$.

\begin{lemm} For any $p\dly 0$, $f \in A\Pi^{p}_{\infty}$ if and only if
\beq \lim_{t\to \infty }\big(Q(t+ct^{-p})-Q(t)\big)= 0\quad
\mbox{for\ any} \quad c>0.
\eeq In this case, \beq Q(t)=\mro(t^{p+1}) \quad \mbox{as} \quad
t\to\infty.
\eeq
\end{lemm}

\medskip

\begin{lemm}
For any $p\dly 0$, $f \in A\Pi^{p}_{0}$ if and only if \beq
\lim_{t\to \infty }\big(Q(t+ct^{-p})-Q(t)\big)= \infty\quad
\mbox{for\ any} \quad c>0.
\eeq In this case, \beq Q(t)t^{-p-1} \to \infty\quad \mbox{as}
\quad t\to\infty.
\eeq

\end{lemm}

\medskip

\begin{lemm} For any $p\dly 0$ and $f \in \rO A\Pi^p$, the following assertions hold
true:

{\rm (i)} there is a constant $c>0$ such that $Q(t+ct^{-p})-Q(t)
\asymp 1$ as $t\to \infty$;

{\rm (ii)} $Q(t) \asymp t^{p+1}$ as $t \to \infty$;

{\rm (iii)} if a function $a: \ttR_+ \rightarrow \ttR_+$ is chosen
to satisfy $\linf_{s\to \infty }a(s)\dly c_1>0$ and $\linf _{s\to
\infty }(s-a(s))\dly c_2>0$, then
\begin{eqnarray*}
&&\const(s-a(s))s^{-p/(p+1)}\mly f(s)-f(a(s)) \\[6pt]
&&\mly \const \pr \lf
(s^{1/(p+1)}-a(s)^{1/(p+1)}+a(s)^{-p/(p+1)}\rg )
\end{eqnarray*} for any $s\dly s_0$ and for some $\const\pr \dly \const >0$.

\end{lemm}

\medskip

Before proving Lemmas A.6--A.8, we provide an example of $Q$
satisfying assertion~(i) of Lemma A.8 such that
$f:=Q^{\leftarrow}$ does not belong to $\rO A\Pi^p$ for $p \dly
0$, and further on, two technical lemmas for later use.

\medskip

\begin{exm} For $p \dly 0$, write $Q(t):=t^{p+1}+[t]$, $t \dly 0$.
Note that, for each $c > 0$,
$$
Q(t+ct^{-p})-Q(t)=c(p+1)+g(t) +\mro(1)\quad \mbox{as} \quad
t\to\infty,
$$
where $0 \mly g(t):=[t+ct^{-p}] -[t] =\rO (1)$. I.e., $Q$
satisfies the assertion of Lemma~A.8(i) for any $c >0$. However,
for each $t:=n \in \ttN$, we get $Q(n)-Q(n-)=1$, therefore,
$f:=Q^{\leftarrow} \notin \rO A \Pi^p$ according to Lemma~A.10
below.
\end{exm}

\medskip

\begin{lemm}
If $\linf_{s \to \infty} f(s)^p\big ( f(s+c)-f(s) \big ) >0$ for
each $c >0$ and for some $p \dly 0$, then
$$
\lim_{t\to \infty }\big(Q(t)-Q(t-)\big)= 0,
$$
i.e., $Q$ is continuous at infinity.
\end{lemm}

\begin{proof}
Assume for a moment that there exists a sequence $t_n\to \infty $
such that $Q (t_n)-Q (t_n-)\to c^0>0$. This limit implies that
$s_n:=Q (t_n)\to \infty $ and, in addition, that $f(s_n-c)=f(s_n)$
for any $0<c<c^0$ and any $n\dly n_0(c)$, contradicting the
assumption of the lemma. This completes the proof of the claimed
assertion.\qed
\end{proof}

\medskip

\begin{lemm} {\rm (Resnick 1987, pp. 4).}
For all $s \in \ttR_+$ and $t \in (-\infty; t_Q)$, the following
assertions hold true:

{\rm (i)} $f (s) \mly t$ if and only if $s \mly Q(t)$;

{\rm (ii)} $f (s) > t$ if and only if $s > Q(t)$;

{\rm (iii)} $Q(f(s)-) \mly s \mly Q(f(s))$.

\end{lemm}

\medskip

We now are in a position to prove Lemmas A.6--A.8. To simplify the
proceedings, we need the following abbreviations:
$$
f_p(s;c):= f(s)^p \big(f(s+c)-f(s) \big) \quad \mbox{and} \quad
Q_p(t;c):=Q(t+ct^{-p})-Q(t).
$$

\smallskip

{\it Proof of Lemma A.6.} Assume first that (A.4) holds true.
I.e., for each $\rre >0$ there is $s_0=s_0(\rre) >0$ such that
\beq Q_p(f(s);\rre^{-1}) < \rre \quad \mbox{for\ all} \quad s \dly
s_0. \eeq
Since (A.4) implies a continuity of $Q$ at infinity, from Lemma
A.11(iii) we have that $Q(f(s)) \mly s+\rre$ for each $s \dly
s_0$. This and (A.8) yield that
$Q\big(f(s)+\rre^{-1}f(s)^{-p}\big) < s+2 \rre$ for each $s \dly
s_0$. Inverting $Q$ (see Lemma A.11(ii)), we get that $f_p(s; 2
\rre) > 1/ \rre$ for each $\rre > 0$ and each $s \dly s_0(\rre)$.
I.e., $f \in A\Pi^{p}_{\infty}$.

To prove the inverse implication, assume for a moment  that there
are a sequence of reals $t_n \to \infty$ and constants $c > 0$,
$\rre >0$ such that $Q(t_n+ct_n^{-p}) \dly Q(t_n) +\rre$ for all
$n \dly n_0(\rre, c)$. Inverting $Q$ (see Lemma A.11(i)), we get
that $f \big( Q(t_n) +\rre \big) \mly t_n+ct_n^{-p}$, which
combined with $t_n < f \big( Q(t_n) +\rre /2 \big)$ (see Lemma
A.11(ii)) gives that, for each $n \dly n_0(\rre, c)$,
$$
f_p(s_n; \rre /2) \mly c \quad \mbox{with} \quad s_n:= Q(t_n)+
\rre /2.
$$
Since $s_n \to \infty$, the latter violates the assumption $f \in
A\Pi^{p}_{\infty}$, concluding the proof of the first part of
Lemma A.6.

To prove (A.5), we note that, for any natural $M \dly 2$ and any
$s\dly 2M$,
$$
f(s+1)^{p+1}-f(s)^{p+1}\dly f_p(s; 1)\dly I(M):=\inf _{s\dly
M}f_p(s; 1),
$$
and, therefore,
$$
f(s)^{p+1}-f(M)^{p+1}\dly (s-M-1)I(M).
$$
Hence $\linf _{s\to \infty }f(s)s^{-1/(p+1)}\dly I(M)^{1/(p+1)}$.
Since $I(M)\to \infty $ (as $M\to \infty $) by the assumption, the
latter implies that $f(s)s^{-1/(p+1)}\to \infty $ as $s\to \infty
$, which in turn yields (A.5). Lemma A.6 is proved.\qed

\medskip

{\it Proof of Lemma A.7.} Assume first that (A.6) holds true,
i.e., for each $\rre > 0$ there  is $t_0=t_0(\rre)$ such that
$Q_p(t; \rre) \dly 1/\rre$ for each $t \dly t_0$. By Lemma
A.11(i), the latter is equivalent to $f \big( Q(t)+1/\rre \big)
\mly t+ \rre t^{-p}$. Substituting $t:=f(s) \to \infty$ into this
inequality and then applying $Q(f(s)) \dly s$ (see Lemma
A.11(iii)), we obtain $f_p(s; 1/\rre) \mly \rre$ for each $\rre >
0$ and each $s \dly s_0(\rre)$, i.e., $f \in A\Pi^{p}_{0}$.

To prove the inverse implication, assume for a moment that there
are a sequence $t_n \to \infty$ and constants $\delta >0$, $c >0$
such that $Q_p(t_n; \delta) < c$ for each $n \dly n_0(c, \delta)$.
Here, inverting $Q$ (see Lemma A.11(i),(ii)) and denoting
$s_n:=Q(t_n) \to \infty$, we obtain that $f_p(s_n; c) > \delta$
for each $n \dly n_0(c, \delta)$, contradicting the assumption $f
\in A\Pi^{p}_{0}$. This completes the proof of the first part of
Lemma A.7.

We will prove (A.7) under the weaker condition by assuming (A.6)
for some $c > 0$. (The forthcoming arguments are applied to prove
assertion (ii) of Lemma A.8 as well). Write $Q^{(c)}(t):=
Q(c^{1/(p+1)} t)$ and observe that
$$
\lim_{t \to \infty}Q^{(c)}_p(t; 1)=\lim_{t \to \infty}Q_p(t; c)=
\infty,
$$
i.e., limit (A.6) for $c
> 0$ is reduced to that for $c=1$. With the abbreviation $R(k):=[k^p]$, we obtain that, for
fixed natural $M\dly 3$ and any natural $t\dly 2M$,
\begin{eqnarray*}
&&\hspace{-.7cm}Q^{(c)} (t)-Q^{(c)} (M)\\
&&\hspace{-.7cm}\quad {}=\sum _{k=M}^{t-1}\sum _{l=0}^{R(k)-1}\lf
(Q^{(c)} \Big
(k+\frac{l+1}{R(k)}\Big )-Q^{(c)} \Big (k+\frac{l}{R(k)}\Big )\rg ) \\
&&\hspace{-.7cm}\quad {}\!\dly\!\sum _{k=M}^{t-1}\sum
_{l=0}^{R(k)-1}Q^{(c)}_p \Big (k+\frac{l}{R(k)}; 1\Big
)\!\dly\!\inf_{\tau \dly M}Q^{(c)}_p  (\tau; 1)\sum
_{k=M}^{t-1}R(k)\\
&&\hspace{-.7cm}\quad {}=\inf_{\tau \dly M}Q^{(c)}_p  (\tau;
1)\frac{t^{p+1}}{p+1}\lf (1+\mro (1)\rg )
\end{eqnarray*}
as $t\to \infty $. Here, by (A.6), the infimum tends to infinity
as $M\to \infty $, therefore, (A.7) is fulfilled. This completes
the proof of Lemma A.7.\qed

\medskip

{\it Proof of Lemma A.8.} (i) We first prove that if, for each $c
> 0$, the function $f_p(s; c) $ is asymptotically bounded away from
zero as $s \to \infty$, then $Q_p(t; \delta) =\rO (1)$ as $t \to
\infty$, for some $\delta
>0$. Assume otherwise that, for each $\delta > 0$, there exists a
sequence $t_n \to \infty$ such that $Q_p(t_n; \delta) \dly 2M$ for
any $M > 0$ and any $n \dly n_0(M)$. Here, inverting $Q$ similarly
as in the proof of Lemma A.6, we obtain that $f_p(s_n; M) \mly
\delta$ with $s_n :=Q(t_n) +M \to \infty$ as $n \to \infty$.
Therefore, $\linf_{s \to \infty}f_p(s; M) \mly \delta$. Since
$\delta > 0$ is arbitrary, we obtain the contradiction proving the
desired implication. We next observe that, if
$Q_p(t;\delta)=\rO(1)$ for some $\delta >0$, then \beq Q_p (t;
k\delta)=\rO(1) \quad \mbox {as} \quad t \to \infty,  \quad \mbox
{for\ any} \quad k \in \ttN.
\eeq This implication is easily proved by induction in $k$. We
omit the details.

With the abbreviation
$$
M:= 2 \lsup_{s \to \infty}f_p(s; 1) > 0,
$$
we finally show that the function $Q_p(t; M)$ is asymptotically
bounded away from zero as $t \to \infty$. For this, fix an
arbitrary sequence $t_n\to \infty $, and define a sequence
$\{s_n\}$ by $f(s_n+)\dly t_n\dly f(s_n)$ ($n\in \ttN $). Write
$\tau_n:=f(s_n)$. Combining Lemmas A.10 and~A.11(iii), we have
that $Q (t_n)-Q (\tau_n)=\mro (1)$ and, consequently, \beq Q_p
(t_n; M)\dly Q_p (\tau_n; M)+\mro (1) \quad \mbox {as} \quad n \to
\infty.
\eeq On the other hand, from the definition of $M$, it follows
that $f(s_n)+M f(s_n)^{-p}>f(s_n+1)$ for any $n\dly n_0$. Applying
$Q$ to both sides of this inequality and then using Lemmas A.10
and A.11(iii), we obtain that $Q_p (\tau_n; M)\dly 1/2$ for $n
\dly n_0$. The latter combined with (A.10) implies that the
sequence $Q_p(t_n; M)$ is asymptotically bounded away from zero,
as claimed. This and (A.9) conclude the proof of part (i).

(ii) The assertion is shown by the same arguments as in the proof
of limits (A.5) and (A.7). We omit the details.

(iii) If $a(s)$ or $s-a(s)$ are bounded from above for any large
$s$, then the bounds in (iii) simply follow from part (ii) and
condition (4.3).

For simplicity we abbreviate $a:=a(s)$, and assume that both $a$
and $s-a$ tend to infinity as $s \to\infty$. By combining
assumption (4.3) and the limit $f(s)\asymp s^{1/(p+1)}$, we obtain
that, for $s\dly s_0$,
\begin{align*}
f(s)-f(a)&\mly \sum _{0\mly k\mly s-a}\lf (f(a+k+1)-f(a+k)\rg ) \\
&\mly \const \sum _{0\mly k\mly s-a}(a+k)^{-p/(p+1)} \\
&\mly \const
a^{-p/(p+1)}+\const \int_0^{s-a}(a+k)^{-p/(p+1)}\dd k \\
\noalign{\noindent\mbox{and}}
f(s)-f(a)&\dly \sum _{0\mly k\mly s-a-1}\lf (f(a+k+1)-f(a+k)\rg ) \\
&\dly \const \sum _{0\mly k\mly s-a-1}(a+k)^{-p/(p+1)} \\
&\dly  \const (s-a-1)s^{-p/(p+1)},
\end{align*}
as claimed. Lemma A.8 is proved.\qed

\medskip

\begin{remm} (A relationship with classical regular variation).
(i) Consider the case $p=0$. Obviously, for $\beta=\infty$ or
$\beta= 0$, $f$ is in $A\Pi^0_\beta $ if and only if $g:=\exp
\circ f \circ \log \in RV_\beta$. Therefore, for $p=0$, Lemmas A.6
and A.7 follow from the well-known results for the class
$RV_\beta$ with $\beta =\infty$ and $\beta=0$, respectively
(Bingham et al.~1987).

The class $\rO A \Pi^0$ links to the exponential type
distributions $1-\e^{-Q}$, with $Q(t) \asymp t$ as $t \to \infty$.
Moreover, if $f \in \rO A\Pi^0$, then $g:=\exp \circ f \circ \log$
is in $\rO RV$, the class of $\rO$-regularly varying functions
studied, e.g., in (Bingham et al.~1987, Section 2).

(ii) In the case of $p \dly 0$, if $f \in \rO A\Pi^p$, then $f
\circ \log$ is asymptotically balanced or, equivalently, the
maximum $\xi _{1, V}$ of i.i.d.~sample $\xi _{V}$ is
stochastically compact (Bingham et al.~1987, Sections 3.11 and
8.13.12).
\end{remm}

\medskip

\subsection{\normalsize \textsf {Class $PI_{<2}$} }

Recall that the class $PI_{<2}$ consists of functions $f:=Q^{\la}$
such that
$$
\linf_{s \to \infty}\frac{f(c s)}{f(s)}> 1 \quad \mbox {for\ some}
\quad 1 < c <2;
$$
cf.~(4.5). This is the subclass of positive increase class $PI$
considered, e.g., in (Bingham et al.~1987, Section 2.1.2).

\begin{lemm}
$f \in PI_{<2}$ if and only if \beq \lsup_{t \to \infty} Q(C_0
t)/Q(t) <2 \quad \mbox{for\ some} \quad C_0>1.
\eeq In this case, \beq \lsup_{t \to \infty} \frac{\log Q(t)}{\log
t} \mly \const :=\frac {\log 2}{\log C_0}.
\eeq

\end{lemm}

\medskip

\begin{proof}
Assume in contrary to (A.11) that, for any $C > 1$, there is a
sequence $t_n =t_n(C) \to \infty$ such that $Q(Ct_n) \dly (2-\rre)
Q(t_n)$ for any $\rre > 0$ and any $n \dly n_0(\rre)$. Similarly
as in the proof of Lemma A.6, inverting $Q$ (see Lemma
A.11(i),(ii)) and writing $s_n:=Q(t_n) +\rre \to \infty$, we
obtain that
$$
\lsup_{n} f(\delta s_n)/f(s_n) \mly C \quad \mbox {for\ any} \quad
1 < \delta <2 \quad \mbox{and\ any} \quad C>1
$$
or, equivalently, $ \lim_{n} f(\delta s_n)/f(s_n)=1$,
contradicting the assumption $f \in PI_{<2}$. According to these
arguments, the inclusion $f \in PI_{<2}$ implies (A.11).

To show the inverse implication, we suppose otherwise that
$f:=Q^{\la} \notin PI_{<2}$, i.e., for each $1< c<2$, there exists
a sequence $s_n =s_n(c) \to \infty$ such that $f(cs_n) \mly (1+
\rre)f(s_n)$ for any $\rre
> 0$ and any $n \dly n_0(\rre)$.
In this inequality, we invert $f$ (see Lemma~A.11(i)) to obtain
$cs_n \mly Q\big((1+ \rre)f(s_n) \big)$. On the other hand, by
Lemma~A.11(iii),
$$
cs_n \dly cQ\big(f(s_n)- \big)\dly cQ\big((1- \rre)f(s_n) \big).
$$
Summarizing these estimates and using the abbreviations $t_n:=(1-
\rre)f(s_n)$ and $\delta :=(1+\rre)/(1-\rre)$, we have that $Q(
\delta t_n) \dly cQ(t_n)$. Since $1<c<2$ is an arbitrary constant
but close to 2, the latter implies the limit $\lsup_{t}Q( \delta t
)/Q(t) \dly 2$ for each $\delta > 1$, contradicting assumption
(A.11). This concludes the first part of the lemma.

Let us show (A.12). By (A.11), there exist numbers $C > 1$ and
$t_0=t_0(C)$ such that $Q(Ct)/ Q(t) \mly 2$ for all $t \dly t_0$.
Applying this estimate, we obtain that, for any $n \in \ttN $ and
any $t \in \big[ C^n t_0; C^{n+1}t_0 \big)$,
\begin{eqnarray*}
&&Q(t) =\frac{Q(t)}{Q(C^nt_0)} \cdot \frac{Q(C^n
t_0)}{Q(C^{n-1}t_0)} \cdots \frac{Q(Ct_0)}{Q(t_0)} \cdot Q(t_0)\\[6pt]
&&\hspace{.7cm}\mly Q(t_0)2^{n+1} \mly \const t^{(\log 2 )/\log
C},
\end{eqnarray*}
i.e., (A.12) is done. Lemma A.13 is proved.\qed

\end{proof}

\subsection{\normalsize \textsf {Comparison of the classes
$A\Pi^p_\infty$, $A\Pi$ and $PI_{<2}$. Examples} }

In view of limit theorems for eigenvalues (see Theorems 6.2 and
6.9), we need to compare the classes $A\Pi^{p}_{\infty}$ (4.1), $S
A \Pi^{2}_{\infty}$ (5.4), $A\Pi$ (6.3) and $PI_{<2}$ (4.5) of
functions $f:=Q^{\la}$.

\begin{lemm}
{\rm (i)} For any $p \dly 0$, there exist examples $f_1 \in
PI_{<2}\backslash A\Pi^{p}_{\infty}$ and $f_2 \in
A\Pi^{p}_{\infty} \backslash PI_{<2}$. Consequently, there is $f_2
\in S A\Pi^{2}_{\infty} \backslash PI_{<2}$.

{\rm (ii)} There exist examples $f_3 \in PI_{<2}\backslash A\Pi$
and $f_4 \in A\Pi \backslash PI_{<2}$ with an auxiliary function
$a_4
 \dly 1$.

{\rm (iii)} For any $p>0$, there exists an example $f_5 \in \big(
A\Pi^{p}_{\infty} \cap PI_{<2} \big) \backslash A\Pi$ and,
therefore, there is $f_5 \in \big(S A \Pi^{2}_{\infty} \cap
PI_{<2} \big) \backslash A\Pi$.

{\rm (iv)} For $p \dly 0$, if $f \in A\Pi$ with an auxiliary
function $a: \ttR_+ \rightarrow \ttR_+$ such that $s^p a(s) \to
\infty$ as $s \to \infty$, then $f \in A \Pi^{p}_{\infty}$.

{\rm (v)} For $0 \mly \rho < \infty$, if $\e^f$ is in $RV_{\rho}$,
then $f(s+ \log s)- f(s) \to 0$ as $s \to \infty$.

\end{lemm}

\medskip

\begin{proof}
(i) Write $Q(t):=[t]$, $t \dly 0$; i.e., $1-\e ^{-Q}$ is the
geometric distribution. Let us show that $f_1:=Q^{\la} \in
PI_{<2}\backslash A\Pi^{p}_{\infty}$. Indeed, since $Q(n)-Q(n-)=1$
for all $n \in \ttN$, Lemma A.10 implies that $f_1 \notin
A\Pi^{p}_{\infty}$ for any $p \dly 0$. However, since $Q$
satisfies (A.11), we have that $f_1 \in PI_{<2}$, as claimed.

Consider the function $f_2(s):=\int_1^sb(t)\dd t$, where $b(t):=n$
if $2^{2n}<t\mly 2^{2n+1}$, and $b(t):=2^n$ if $2^{2n+1}<t\mly
2^{2n+2}$ for all $n \in \ttN \cup \{0\}$. Let us show that $f_2
\in A\Pi^{0}_{\infty} \backslash PI_{<2}$. Obviously $f_2 \in
A\Pi^{0}_{\infty}$. With $s_n:= 2^{2n+1}$ and $1/2 < \delta <1$,
we see that
$$
\frac {f_2(\delta s_n)}{f_2(s_n)}=1- \frac{\int_{\delta s_n}^{s_n}
b(t) \dd t}{\int_{1}^{s_n} b(t) \dd t},
$$
where $\int_{\delta s_n}^{s_n} b(t) \dd t=\const 4^n n$ and
\begin{eqnarray*}
&&\int_{1}^{s_n} b(t) \dd t =\sum_{l=0}^{n-1} \bigg( \int
^{s_l}_{s_l/2}l \dd t +\int ^{2s_l}_{s_l}2^l \dd t \bigg)+\int
^{s_n}_{s_n/2}n \dd t \\[6pt]
&&\hspace{.7cm} =\sum_{l=0}^{n-1} \Big( l\cdot 4^l +2\cdot 8^l
\Big)+n\cdot 4^n =\frac{2}{7}\cdot 8^n (1+\mro(1))
\end{eqnarray*}
as $n \to \infty$. Summarizing, we find that $f_2(\delta s_n)/
f_2( s_n) \to 1$ (as $n \to \infty$) for each $1/2 < \delta <1$,
i.e., $f_2 \notin PI_{<2}$. Since $A\Pi_{\infty}^1 \subset S
A\Pi_{\infty}^2$, we also obtain that $f_2 \in S A\Pi^{2}_{\infty}
\backslash PI_{<2}$.

(ii) As in part (i) above, let $f_3$ be the inverse of the
cumulative hazard function of the geometric distribution. Since
$f_3$ is not in $A\Pi$ (Resnick 1987, Corollary 1.6), we obtain
that $f_3 \in PI_{<2}\backslash A\Pi$.

To prove the existence of $f_4 \in A\Pi \backslash PI_{<2}$, it
suffices (via Lemmas A.1 and A.13) to find a continuous function
$a: [1; \infty) \rightarrow [1; \infty)$, with derivative
$a^{\prime}(t) \to 0$ as $t\uparrow \infty$, such that the
function $Q(t):=\int^t_{1}1/a(s)\dd s $ does not satisfy (A.11).
For this, we abbreviate $t_n := (\log n)^n$, $m_n:= (\log n)^2$,
$$
\rre_n :=\frac{1}{m_n} \lf (\frac{t_{n+1}}{t_n}-1 \rg
)-\frac{1}{t_n} \quad \mbox {and} \quad b_n:= 1+t_{n+1} \rre_n ,
$$
and consider the functions \beq a(t):=\begin{cases} {\displaystyle
\ \ t/t_n} &
\hbox{if\ \ $t_n < t \mly t_{n+1}- m_n,$} \\
-\rre_n t +b_n &\hbox{if\ \ $t_{n+1}- m_n < t \mly t_{n+1}$;\ \ \
for $n \dly 2$,}
\end{cases}
\eeq
and $Q(t):=\int^t_{1}1/a(s)\dd s $. Obviously, $a \dly 1$ is
continuous in $[1; \infty)$ and $a^{\prime}(t) \to 0$ as
$t\uparrow \infty$. Therefore, $f_4:=Q^{\la} \in A\Pi$ with the
auxiliary function $a$. Let us show that, for each $c > 1$, \beq
\frac{Q(ct_n)}{Q(t_n)}=1+ \frac{\int^{ct_n}_{t_n}1/a(t) \dd
t}{\int^{t_n}_{1}1/a(t) \dd t} \to \infty \quad \mbox {as} \quad n
\to \infty,
\eeq so that $f_4 \notin PI_{<2}$. Indeed, from (A.13) we see
that, for $n \dly n_0(c)$, \beq \int^{ct_n}_{t_n} \frac{ \dd
t}{a(t)}=\int^{ct_n}_{t_n} \frac{ t_n}{t} \dd t=t_n \log c .
\eeq To estimate the integral $\int^{t_n}_1 $ in (A.14), we again
use (A.13) and the bound $a \dly 1$. Thus,
\begin{eqnarray}
&&\int^{t_n}_{t_2} \frac{ \dd t}{a(t)} \mly \sum^{n-1}_{k=2} \bigg
( t_k \int^{t_{k+1}-m_k}_{t_k}\frac{\dd t}{t} +
 \int^{t_{k+1}}_{t_{k+1}-m_k}1 \dd t \bigg) \n \\
&&\quad {}\mly \sum^{n-1}_{k=2} \big( t_k \log(t_{k+1}/t_k) +m_k \big
) \mly t_n (\log n)^{-1/2}
\end{eqnarray}
for any $n \dly n_0$. Now (A.15) and (A.16) imply (A.14), as
claimed.

(iii) Consider the example $Q(t) =t +\sin t$ ($t \dly 0$) given by
Von Mises. Obviously, $Q$ satisfies (A.4) and (A.11),
consequently, $f_5:=Q^{\la}$ is in $ A\Pi^{p}_{\infty} \cap
PI_{<2} $ for any $p
> 0$. However, $f_5 \notin A\Pi$. We observe that the function
$f_6(s):=s+ \sin s$ ($s \dly 0$) is also in $\big(
A\Pi^{p}_{\infty} \cap PI_{<2} \big) \backslash A\Pi$ for any $p
> 0$. (This is verified by straightforward calculations).
Consequently, since $ A\Pi^{1}_{\infty}\subset S
A\Pi^{2}_{\infty}$, the functions $f_5$ and $f_6$ are in $\big(S
A\Pi^{2}_{\infty} \cap PI_{<2} \big) \backslash A\Pi$.

(iv) The assertion follows from the definition of $A\Pi$ (6.3) and
$A\Pi^p_{\infty}$ (4.1).

(v) The assertion follows from Lemma A.3(ii). Lemma A.14 is
proved.\qed
\end{proof}

We finally provide two examples of distributions which represent
RV classes considered in Sections 3--6.
\medskip

\begin{exm} For $\rra > 0$,
let $Q_{\rra}(t)= t^{\rra}$ ($t \dly 0$), i.e., $1-\e^{-Q_{\rra}}$
is Weibull distribution. Clearly $
f_\rra(s):=Q_{\rra}^{\leftarrow}(s)=s^{1/\rra }$ for $s \dly 0$.
By straightforward calculations, we obtain that if $\rra < p+1$
(resp., $\rra > p+1$ or $\rra = p+1$), then $f_{\rra} \in
A\Pi^p_{\infty}$ (resp., $f_{\rra} \in A\Pi^p_{0}$ or $f_{\rra}
\in \rO A\Pi^p$). Also, for $\rra <3$, $f_{\rra} \in S
A\Pi^2_{\infty}$. Finally, for any $\rra > 0$, $f_{\rra} \in
PI_{<2}$, $\exp \circ f_{\rra} \in RV_{\infty}$ and $f_{\rra}$ is
in $A\Pi $ with the auxiliary function $a(t):=\rra ^{-1}
t^{1-\rra}$. The latter means that the distribution
$1-\e^{-Q_{\rra}}$ is in the domain of attraction of the
max-stable Gumbel law; cf.~Section 6.1.
\end{exm}

\medskip

\begin{exm} Given $\rrg > 0$ and $\rho >0$, let
$Q_{\rrg, \rho}(t)=\e^{\rho^{-1} t^{\rrg}}$ ($t \dly t_0$), i.e.,
the fractional double exponential distribution. Then $ f_{\rrg,
\rho}(s)=(\rho \log s)^{1/\rrg }$ for $s \dly s_0$. Obviously, if
$0 <\rrg <1$ (resp., $\rrg > 1$ or $\rrg =1$ ), then $\exp \circ
f_{\rrg, \rho}$ is in $RV_{\infty}$ (resp., $RV_{0}$ or
$RV_{\rho}$); cf.~Section 6.
\end{exm}

\medskip

{\bf Acknowledgments}. The author would like to thank the referee who reads the first version of this paper for his insightful remarks, critical comments and suggestions.

\medskip

\newcommand{\ms}{\vspace{2pt}}

\section*{\normalsize \textsf {REFERENCES}}

\footnotesize


\begin{list}{}{\topsep 2pt\parsep 0pt\itemsep 0pt\leftmargin 18pt}

\item[1.] Adler, R. J., Taylor, J. E.: {\it Random Fields and
Geometry}. Springer, New York (2007)

\item[2.] Aizenman, M., Molchanov, S.: Localization at large
disorder and at extreme energies: an elementary derivation. {\it
Commun. Math. Phys.} {\bf 157}, 245--278 (1993)

\item[3.] Aizenman, M., Schenker, J. H., Friedrich, R. M.,
Hundertmark, D.: Finite-volume fractional-moment criteria for
Anderson localization. {\it Commun. Math. Phys.} {\bf 224},
219--253 (2001)

\item[4.]  Anderson, G. W.,  Guionnet, A., Zeitouni, O.: {\it An
introduction to random matrices}. Cambridge Studies in Advanced
Mathematics, vol. 118. Cambridge University Press, Cambridge
(2010)

\item[5.]  Anderson, P. W.: Absence of diffusion in certain random
lattices. {\it Phys. Rev.}  {\bf 109}, 1492--1505 (1958)

\item[6.] Astrauskas, A.: On high-level exceedances of i.i.d.~random fields. {\it Stat. Probab. Letters} {\bf 52}, 271--277
(2001)

\item[7.] Astrauskas, A.: On high-level exceedances of Gaussian
fields and the spectrum of random Hamiltonians. {\it Acta Appl.
Math.} {\bf 78}, 35--42 (2003)

\item[8.] Astrauskas, A.: Strong laws for exponential order
statistics and spacings. {\it Lithuanian Math. J.} {\bf 46},
385--397 (2006)

\item[9.] Astrauskas, A.: Poisson-type limit theorems for
eigenvalues of finite-volume Anderson Hamiltonian. {\it Acta Appl.
Math.} {\bf 96}, 3--15 (2007)

\item[10.] Astrauskas, A.:  Extremal theory for spectrum of random
discrete Schr\"odinger operator. I. Asymptotic expansion formulas.
{\it J. Stat. Phys.} {\bf 131}, 867--916 (2008)

\item[11.] Astrauskas, A.:  Extremal theory for spectrum of random
discrete Schr\"odinger operator. II. Distributions with heavy
tails. {\it J. Stat. Phys.} {\bf 146}, 98--117 (2012)

\item[12.] Astrauskas, A.:  Extremal theory for spectrum of random
discrete Schr\"odinger operator. III. Localization properties.
{\it J. Stat. Phys.} {\bf 150}, 889--907 (2013)

\item[13.] Astrauskas, A.: Asymptotic expansion formulas for the
largest eigenvalues of finite-volume Anderson Hamiltonians with
fractional double exponential tails. In preparation (2016)

\item[14.] Astrauskas, A., Molchanov, S. A.: Limit theorems for the
ground states of the Anderson model. {\it Funkts. Anal.
Prilozhen.} {\bf 26}:4, 92--95 (1992); English transl.: {\it
Funct. Anal. Appl.} {\bf 26}, 305--307 (1992)

\item[15.] Auffinger, A., Ben Arous, G., P\'ech\'e, S.: Poisson
convergence for the largest eigenvalues of heavy tailed random
matrices. {\it Ann. Inst. H. Poincar\'e Probab. Statist.} {\bf
45}, 589--610 (2009)

\item[16.] Bai, Z. D., Yin, Y. Q.: Necessary and sufficient
conditions for almost sure convergence of the largest eigenvalue
of a Wigner matrix.  {\it Ann.
Probab.} {\bf 16}, 1729--1741 (1988)

\item[17.] Benaych-Georges, F., P\'ech\'e, S.: Localization and delocalization for heavy tailed band matrices. {\it Ann. Inst. H. Poincar\'e Probab. Statist.} {\bf
50}, 1385--1403 (2014)

\item[18.] Bingham, N. H., Goldie, C. M., Teugels, J. L.: {\it
Regular Variation}. Cambridge University Press, Cambridge (1987)

\item[19.] Binswanger, K., Embrechts, P.: Longest runs in coin
tossing. {\it Insurance Math. Econom.} {\bf 15}, 139--149 (1994)

\item[20.] Biroli, G., Bouchaud, J.-P., Potters, M.: On the top
eigenvalue of heavy-tailed random matrices. {\it Europhys. Lett.
EPL} {\bf 78(1)}, Art 10001, 5 pp (2007)

\item[21.] Bishop, M., Wehr, J.: Ground state energy of the
one-dimensional discrete random Schr\"odinger operator with
Bernoulli potential. {\it J. Stat. Phys.} {\bf 147}, 529--541
(2012)

\item[22.] Biskup, M., Fukushima, R., K\"onig, W.:  Eigenvalue
fluctuations for lattice Anderson Hamiltonians. Preprint {\it
arXiv:1406.5268 [math.PR]} (2014)

\item[23.] Biskup, M., K\"onig, W.:  Long-time tails in the
parabolic Anderson model with bounded potential. {\it Ann.
Probab.} {\bf 29}, 636--682 (2001)

\item[24.] Biskup, M., K\"onig, W.: Eigenvalue order statistics
for random Schr\"odinger operators with doubly-exponential tails. {\it Commun. Math. Phys.} {\bf
341}, 179--218 (2016)

\item[25.] Bourgade, P., Erd\H{o}s, L., Yau, H.-T.: Edge
universality of beta ensembles. {\it Commun. Math. Phys.} {\bf
332}, 261--353 (2014)

\item[26.] Carmona, R., Klein, A., Martinelli, F.: Anderson localization for
Bernoulli and other singular potentials. {\it Commun. Math. Phys.} {\bf
108}, 41--66 (1987)

\item[27.] de Haan, L., Ferreira, A.: {\it Extreme Value Theory:
An Introduction}. Springer, New York (2006)

\item[28.] Deheuvels, P.: Strong laws for the $k$th order
statistics when $k\mly c\log_2n$. {\it Probab. Theory Relat.
Fields} {\bf 72}, 133--154 (1986)

\item[29.] Devroye, L.: Upper and lower class sequences for
minimal uniform spacings. {\it Z. Wahrsch. Verw. Gebiete} {\bf
61}, 237--254 (1982)

\item[30.] Elgart, A., Kr\"uger, H., Tautenhahn, M., Veseli\'c,
I.: Discrete Schr\"odinger operators with random alloy-type
potential. In: Benguria, R., Friedman, E., Mantoiu, M. (eds.),
Spectral Analysis of Quantum Hamiltonians, {\it Operator Theory:
Advances and Applications}, vol. 224, pp. 107--131. Springer,
Basel (2012)

\item[31.] Embrechts, P., Kluppelberg, C., Mikosch, T.: {\it
Modelling Extremal Events for Insurance and Finance}. Springer,
Berlin (1997)

\item[32.]  Erd\H{o}s, L., Knowles, A., Yau, H.-T., Yin, J.:
Spectral statistics of Erd\H{o}s-R\'enyi graphs I: Local
semicircle law. {\it Ann.
Probab.} {\bf 41}, 2279--2375 (2013a)

\item[33.]  Erd\H{o}s, L., Knowles, A., Yau, H.-T., Yin, J.:
Delocalization and diffusion profile for random band matrices.
{\it Commun. Math. Phys.} {\bf 323}, 367--416 (2013b)

\item[34.] Fiodorov, A., Muirhead, S.: Complete localisation and
exponential shape of the parabolic Anderson model with Weibull
potential field. {\it Electron. J. Probab.} {\bf 19}, no. 58,
1--27 (2014)

\item[35.] Fr\" ohlich, J., Martinelli, F., Scoppola, E., Spencer,
T.: Constructive proof of localization in the Anderson tight
binding model. {\it Commun. Math. Phys.} {\bf 101}, 21--46 (1985)

\item[36.] Fr\" ohlich, J., Spencer, T.:  Absence of diffusion in the Anderson tight
binding model for large disorder or low energy. {\it Commun. Math.
Phys.} {\bf 88}, 151--184 (1983)

\item[37.] G\" artner, J., den Hollander, F.: Correlation
structure of intermittency in the parabolic Anderson model. {\it
Probab. Theory Relat. Fields} {\bf 114}, 1--54 (1999)

\item[38.] G\" artner, J., K\"onig, W., Molchanov, S. A.: Almost
sure asymptotics for the continuous parabolic Anderson model. {\it
Probab. Theory Relat. Fields} {\bf 118}, 547–-573 (2000)

\item[39.] G\" artner, J., K\"onig, W., Molchanov, S.: Geometric
characterization of intermittency in the parabolic Anderson model.
{\it Ann. Probab.} {\bf 35}, 439--499 (2007)

\item[40.] G\" artner, J., Molchanov, S. A.: Parabolic problems for
the Anderson model. I. Intermittency and related topics. {\it Commun. Math. Phys.} {\bf 132}, 613--655
(1990)

\item[41.] G\" artner, J., Molchanov, S. A.: Parabolic problems for
the Anderson model. II. Second-order asymptotics and structure of
high peaks. {\it Probab. Theory Relat. Fields} {\bf 111}, 17--55
(1998)

\item[42.] Germinet, F., Klopp, F.: Enhanced Wegner and Minami
estimates and eigenvalue statistics of random Anderson models at
spectral edges. {\it Ann. H. Poincar\'e } {\bf 14}, 1263--1285
(2013)

\item[43.] Germinet, F., Klopp, F.: Spectral statistics for random
Schr\"odinger operators in the localized regime. {\it J. Europ.
Math. Soc.} {\bf 16}:9, 1967--2031 (2014)

\item[44.] G\" otze, F., Naumov, A., Tikhomirov, A. N.: Local semicircle law under moment conditions. Part II: Localization and delocalization. Preprint
{\it arXiv:1511.00862v2 [math.PR]} (2015)

\item[45.] Grenkova, L. N., Molchanov, S. A., Sudarev, Yu. N.: On the
basic states of one-dimensional disordered structures. {\it
Commun. Math. Phys.} {\bf 90}, 101--124 (1983).

\item[46.] Grenkova, L. N., Molchanov, S. A., Sudarev, Yu. N.:  The
structure of the edge of the multidimensional Anderson model
spectrum. {\it Teoret. Mat. Fiz.} {\bf 85}:1, 32--40 (1990);
English transl.: {\it Theor. Math. Phys.} {\bf 85}:1, 1033--1039
(1990)

\item[47.] van der Hofstad, R., K\"onig, W., M\"orters, P.: The
universality classes in the parabolic Anderson model. {\it Commun.
Math. Phys.} {\bf 267}, 307--353 (2006)

\item[48.] van der Hofstad, R., M\"orters, P., Sidorova, N.: Weak
and almost sure limits for the parabolic Anderson model with
heavy-tailed potential. {\it Ann. Appl. Prob.} {\bf 18}, 2450-2494
(2008)

\item[49.] Hundertmark, D.: A short introduction to Anderson
localization. In: {\it Analysis and stochastics of growth
processes and interface models}, pp. 194--218. Oxford Univ. Press,
Oxford (2008)

\item[50.] Killip, R., Nakano, F.: Eigenfunction statistics in the
localized Anderson model. {\it Ann. H. Poincar\'e } {\bf 8},
27--36 (2007)

\item[51.] Kirsch, W.: An invitation to random Schr\"odinger
operator. In: {\it Random Schr\"odinger operators, Panor.
Synth\'eses}, vol. 25, pp. 1--119. Soc. Math. France, Paris (2008)

\item[52.] Kirsch, W., Metzger, B.: The integrated density of
states for random Schr\"odinger operators. In: {\it Spectral
theory and mathematical physics: a Festschrift in honor of Barry
Simon's 60th birthday, Proc. Sympos. Pure Math.}, vol.~76, pp.
649--696. Amer. Math. Soc., Providence (2007)

\item[53.] Klopp, F.: Band edge behavior of the integrated density
of states of random Jacobi matrices in dimension 1. {\it J. Stat.
Phys.} {\bf 90}, 927–-947 (1998).

\item[54.] Klopp, F.: Precise high energy asymptotics for the
integrated density of states of an unbounded random Jacobi matrix.
{\it Rev. Math. Phys. }~{\bf 12}(4), 575--620 (2000)

\item[55.] Klopp, F.: Decorrelation estimates for the eigenlevels
of the discrete Anderson model in the localized regime.  {\it
Commun. Math. Phys.} {\bf 303}, 233--260 (2011)

\item[56.] K\"onig, W.: {\it The Parabolic Anderson
Model}. Birkh\" auser, Basel (2016)

\item[57.] K\"onig, W., Lacoin, H., M\"orters, P., Sidorova,
N.: A two cities theorem for the parabolic Anderson model. {\it
Ann. Probab.} {\bf 37}, 347--392 (2009)

\item[58.] Lankaster, P.: {\it Theory of Matrices}. Academic
Press, London (1969)

\item[59.] Leadbetter, M. R., Lindgren, G., Rootz\'en, H.: {\it
Extremes and Related Properties of Random Sequences and
Processes}. Springer, New York (1983)

\item[60.] Lee, J. O., Yin, J.: A necessary and sufficient
condition for edge universality of Wigner matrices. {\it Duke
Math. J. } {\bf 163}(1), 117--173 (2014)

\item[61.] Mehta, M. L.: {\it Random Matrices}, 3rd ed. Elsevier /Academic Press,
Amsterdam (2004)

\item[62.] Minami, N.: Local fluctuation of the spectrum of a
multidimensional Anderson tight binding model. {\it Commun. Math.
Phys.} {\bf 177}, 709--725 (1996)

\item [63.] Minami, N.: Theory of point processes and some basic
notions in energy level statistics. In: {\it Probability and
Mathematical Physics. CRM Proceedings and Lecture Notes}, vol. 42,
pp. 353–-398. Amer. Math. Soc., Providence (2007)

\item[64.] Molchanov, S.: The local structure of the spectrum of
the one-dimensional Schr\"odinger operator. {\it Commun. Math.
Phys.} {\bf 78}, 429--446 (1981)

\item[65.] Molchanov, S. A.: Lectures on random media. In: {\it
Lectures on Probability Theory, Ecole d'Et\'e de Probabilit\'es de
Saint-Flour XXII-1992}. {\it Lect. Notes in Math.}, vol. 1581, pp.
242--411. Springer, Berlin (1994)

\item[66.] Molchanov, S., Vainberg, B.: Scattering on the
system of the sparse bumps: multidimensional case.{\it Applicable
Analysis} {\bf 71}, 167--185 (1998)

\item[67.] Molchanov, S., Vainberg, B.: Spectrum of
multidimensional Schr\"odinger operators with sparse potentials.
In: Santosa, F., Stakgold, I. (eds.) {\it Analytical and
Computational Methods in Scattering and Applied Mathematics}, pp.
231--253. Chapman and Hall/CRC (2000)

\item[68.] Molchanov, S., Zhang, H.: The parabolic Anderson
model with long range basic Hamiltonian and Weibull type random
potential. In: Deuschel, J.-D., Gentz, B., K\"onig, W., von Renesse,
M., Scheutzow, M., Schmock, U. (eds.) {\it Probability in Complex
Physical Systems}, In Honour of Erwin Bolthausen and J\" urgen G\"
artner, Springer Proceedings in Mathematics, vol. 11, pp. 13--31.
Springer, Heidelberg (2012)

\item[69.] Muirhead, S., Pymar, R.: Localization in the
Bouchaud-Anderson model. Preprint {\it arXiv:\ 1411.4032v2
[math.PR]} (2014)

\item[70.] Pastur, L., Figotin, A.: {\it Spectra of Random and
Almost-Periodic Operators}. Springer, Berlin (1992)

\item[71.] Resnick, S. I.: {\it Extreme Values, Regular Variation,
and Point Processes}. Springer, Berlin (1987)

\item[72.] Shorack, G. R., Wellner, J. A.: {\it Empirical
Processes with Applications to Statistics}. Wiley, New York (1986)

\item[73.] Sidorova, N., Twarowski, A.: Localisation and ageing in
the parabolic Anderson model with Weibull potential. {\it Ann.
Probab.} {\bf 42}, 1666--1698 (2014)

\item[74.] Simon, B., Wolff, T.: Singular continuous spectra under
rank one perturbations and localization for random Hamiltonians.
{\it Commun. Pure Appl. Math.} {\bf 39}, 75–-90 (1986)

\item[75.] Sodin, S.: The spectral edge of some random band
matrices. {\it Annals of Mathematics} {\bf 172}, 2223--2251 (2010)

\item[76.] Soshnikov, A.: Universality at the edge of the
spectrum in Wigner random matrices. {\it Commun. Math. Phys.} {\bf
207}, 697--733 (1999)

\item[77.] Soshnikov, A.: Poisson statistics for the largest
eigenvalues of Wigner random matrices with heavy tails. {\it
Elect. Commun. Probab.} {\bf 9}, 82--91 (2004)

\item[78.] Spencer, T.: Random banded and sparse matrices
(Chapter 23). In: Akemann, G., Baik, J., Di Francesco, P. (eds.)
{\it Oxford Handbook on Random Matrix Theory}. Oxford University
Press, Oxford(2011)

\item[79.] Stolz, G.: An introduction to the mathematics of
Anderson localization. {\it Contemp. Math.} {\bf 552}, 71--108
(2011)

\item[80.] Sznitman, A.-S.: {\it Brownian Motion, Obstacles and
Random Media}. Springer, Berlin (1998)

\item[81.] Tao, T., Vu, V.: Random matrices: Universality of local
eigenvalue statistics up to the edge. {\it Commun. Math. Phys.}
{\bf 298}, 549--572 (2010)

\item[82.] Tao, T., Vu, V.: Random matrices: the universality
phenomenon for Wigner ensembles. In: {\it Modern aspects of random
matrix theory, Proc. Sympos. Appl. Math.}, vol. 72, pp. 121--172.
Amer. Math. Soc., Providence (2014)

\item[83.] Tautenhahn, M., Veseli\'c, I.: Discrete alloy-type
models: regularity of distributions and recent results.  {\it
Markov Process. Related Fields} {\bf 21}, 823--846 (2015)

\item[84.] Vu, V., Wang, K.: Random weighted projections, random
quadratic forms and random eigenvectors. {\it  Random Struct. Alg.} {\bf 47}, 792--821 (2015)

\item[85.] Wellner, J. A.: Limit theorems for the ratio of the
empirical distribution function to the true distribution function.
{\it Z. Wahrsch. Verw. Gebiete} {\bf 45}, 73--88 (1978)

\end{list}

 \ed

 \setcounter{equation}{0} \addtocounter{sct}{1}

\renewcommand{\thesection}{\Alph{section}}
\renewcommand{\thesubsection}{\normalsize{\textsf
{\Alph{section}.\arabic{subsection}.}}}

\newcounter{lemm}
\newcounter{skaitl}

\newenvironment{thmm}{\smallskip {\textbf{\textsf{Theorem}}}\
{\bf {\textsf{\Alph{section}.\arabic{prop}.}}}\ \rm }{\smallskip
\addtocounter{prop}{1}}

\renewcommand{\theequation}{\Alph{section}.\arabic{equation}}

\setcounter{lemm}{1} \setcounter{section}{0}

\setcounter{equation}{0} \addtocounter{sct}{1}

